\documentclass[11pt,letterpaper]{article}
\usepackage[lmargin=1.0in,rmargin=1.0in,bottom=1.0in,top=1.0in,twoside=False]{geometry}

\usepackage{fullpage,amssymb,amsmath}
\usepackage{graphicx}
\usepackage{tikz}
\usetikzlibrary{shapes}
\usetikzlibrary{calc,positioning,arrows.meta}
\usetikzlibrary{decorations.pathmorphing}
\tikzset{snake it/.style={decorate, decoration={snake,amplitude=2}}}

\usepackage{skull}
\usepackage[T1]{fontenc}

 \usepackage{xcolor}
 \usepackage{mathtools}
\usepackage{microtype}
\usepackage{amsfonts}
\usepackage{comment}
\usepackage[english]{babel}
\usepackage{mathrsfs}
\usepackage[vlined, ruled, linesnumbered]{algorithm2e}
\DontPrintSemicolon

\usepackage{trimspaces}
\usepackage{nccfoots}
\usepackage{setspace}
\usepackage{inconsolata}
\usepackage{libertine}
\usepackage[absolute]{textpos}
\usepackage{thmtools}
\usepackage{thm-restate}

\usepackage{xspace}
\usepackage{enumerate}

\usepackage{todonotes}

\usepackage{longtable}

\definecolor{blue}{rgb}{0.1,0.2,0.5}
\definecolor{brown}{rgb}{0.6,0.6,0.2}
\usepackage[ocgcolorlinks, linkcolor={blue}, citecolor={brown}]{hyperref}

\usepackage{comment}

\usepackage[amsmath,thmmarks,hyperref]{ntheorem}
\usepackage{cleveref}

\usepackage{standalone}

\crefformat{page}{#2page~#1#3}%
\Crefformat{page}{#2Page~#1#3}%
\crefformat{equation}{#2(#1)#3}%
\Crefformat{equation}{#2(#1)#3}%
\crefformat{figure}{#2Figure~#1#3}%
\Crefformat{figure}{#2Figure~#1#3}%
\crefformat{section}{#2Section~#1#3}
\Crefformat{section}{#2Section~#1#3}
\crefformat{chapter}{#2Chapter~#1#3}
\Crefformat{chapter}{#2Chapter~#1#3}
\crefformat{chapter*}{#2Chapter~#1#3}
\Crefformat{chapter*}{#2Chapter~#1#3}
\crefformat{part}{#2Part~#1#3}
\Crefformat{part}{#2Part~#1#3}
\crefformat{enumi}{#2(#1)#3}
\Crefformat{enumi}{#2(#1)#3}

\usepackage{latexsym}

\theoremnumbering{arabic}
\theoremstyle{plain}
\theoremsymbol{}
\theorembodyfont{\itshape}
\theoremheaderfont{\normalfont\bfseries}
\theoremseparator{.}

\newtheorem{theorem}{Theorem}
\crefformat{theorem}{#2Theorem~#1#3}
\Crefformat{theorem}{#2Theorem~#1#3}

\newcommand{\newtheoremwithcrefformat}[2]{%
  \newtheorem{#1}[theorem]{#2}%
  \crefformat{#1}{##2\MakeUppercase#1~##1##3}%
  \Crefformat{#1}{##2\MakeUppercase#1~##1##3}%
}
\newcommand{\newseptheoremwithcrefformat}[2]{%
  \newtheorem{#1}{#2}%
  \crefformat{#1}{##2\MakeUppercase#1~##1##3}%
  \Crefformat{#1}{##2\MakeUppercase#1~##1##3}%
}
\newcommand{\newnestedtheoremwithcrefformat}[2]{%
  \newtheorem{#1}{#2}[theorem]%
  \crefformat{#1}{##2\MakeUppercase#1~##1##3}%
  \Crefformat{#1}{##2\MakeUppercase#1~##1##3}%
}

\newtheoremwithcrefformat{lemma}{Lemma}
\newtheoremwithcrefformat{proposition}{Proposition}
\newtheoremwithcrefformat{observation}{Observation}
\newtheoremwithcrefformat{conjecture}{Conjecture}
\newtheoremwithcrefformat{corollary}{Corollary}
\newnestedtheoremwithcrefformat{claim}{Claim}

\theorembodyfont{\upshape}
\newtheoremwithcrefformat{example}{Example}
\newtheoremwithcrefformat{remark}{Remark}
\newseptheoremwithcrefformat{definition}{Definition}
\newseptheoremwithcrefformat{question}{Question}

\theoremstyle{nonumberplain}
\theoremheaderfont{\scshape}
\theorembodyfont{\normalfont}
\theoremsymbol{\ensuremath{\square}}
\newtheorem{proof}{Proof}

\theoremsymbol{\ensuremath{\lrcorner}}
\newtheorem{claimproof}{Proof of Claim}

\def\cqedsymbol{\ifmmode$\lrcorner$\else{\unskip\nobreak\hfil
\penalty50\hskip1em\null\nobreak\hfil$\lrcorner$
\parfillskip=0pt\finalhyphendemerits=0\endgraf}\fi} 

\newcommand{\cqed}{\renewcommand{\qed}{\cqedsymbol}}

\tikzset{
    position/.style args={#1:#2 from #3}{
        at=(#3.#1), anchor=#1+180, shift=(#1:#2)
    }
}

\newcommand{\N}{\mathbb{N}}

\newcommand{\Sc}{\mathcal{S}}

\newcommand{\Pc}{\mathcal{P}}
\newcommand{\Cc}{\mathcal{C}}
\newcommand{\Dc}{\mathcal{D}}
\newcommand{\Lc}{\mathcal{L}}
\newcommand{\Ac}{\mathcal{A}}
\newcommand{\Bc}{\mathcal{B}}
\newcommand{\Gc}{\mathcal{G}}

\newcommand{\Tc}{\mathcal{T}}
\newcommand{\Kc}{\mathcal{K}}
\newcommand{\Fc}{\mathcal{F}}

\newcommand{\Qc}{\mathcal{Q}}
\newcommand{\dist}{\mathrm{dist}}
\newcommand{\Root}{\mathsf{Root}}
\newcommand{\Leaves}{\mathsf{Leaves}}
\newcommand{\Path}{\mathsf{Path}}
\newcommand{\Tree}{\mathsf{Tree}}
\newcommand{\QuasiCage}{\mathsf{QuasiCage}}
\newcommand{\RSemiLadder}{\mathsf{RSemiLadder}}
\newcommand{\SGSemiLadder}{\mathsf{SGSemiLadder}}
\newcommand{\QSemiLadder}{\mathsf{QSemiLadder}}
\newcommand{\Cage}{\mathsf{Cage}}
\newcommand{\Split}{\mathsf{Split}}
\newcommand{\Area}{\mathsf{Area}}
\newcommand{\SPath}{\mathsf{SPath}}
\newcommand{\OrderedCage}{\mathsf{OrderedCage}}
\newcommand{\IdOrderedCage}{\mathsf{IdOrderedCage}}
\newcommand{\NeighborCage}{\mathsf{NeighborCage}}
\newcommand{\SeparatingCage}{\mathsf{SeparatingCage}}

\newcommand{\pw}[1]{\mathrm{pw}\left(#1\right)}
\DeclareMathOperator{\lca}{lca}
\DeclareMathOperator{\Int}{Int}
\DeclareMathOperator{\wcol}{wcol}
\DeclareMathOperator{\poly}{poly}
\DeclareMathOperator{\tw}{tw}
\DeclareMathOperator{\rev}{rev}

\newcommand{\len}{\mathrm{len}}

\renewcommand{\epsilon}{\varepsilon}

\let\originalleft\left
\let\originalright\right
\renewcommand{\left}{\mathopen{}\mathclose\bgroup\originalleft}
\renewcommand{\right}{\aftergroup\egroup\originalright}

\renewcommand{\leq}{\leqslant}
\renewcommand{\geq}{\geqslant}

\renewcommand{\setminus}{-}

\renewcommand{\phi}{\varphi}

\usepackage{todonotes}

\newcommand{\prz}[2][]{}
\newcommand{\map}[2][]{}
\newcommand{\mip}[2][]{}

\usepackage{lineno}

\begin{document}

\title{Bounds on half graph orders in powers of sparse graphs\footnote{
The results of this paper have been presented in the master thesis, defended by the author at the University of Warsaw.}}

\author{
Marek Sokołowski\thanks{University of Warsaw, Poland,
  \texttt{marek.sokolowski@mimuw.edu.pl}.
 This work is 
a part of project TOTAL that has received funding from the European Research Council (ERC) 
under the European Union's Horizon 2020 research and innovation programme (grant agreement No.~677651).}
}

\begin{titlepage}
\def\thepage{}
\thispagestyle{empty}
\maketitle

\begin{textblock}{20}(0, 11.7)
\includegraphics[width=40px]{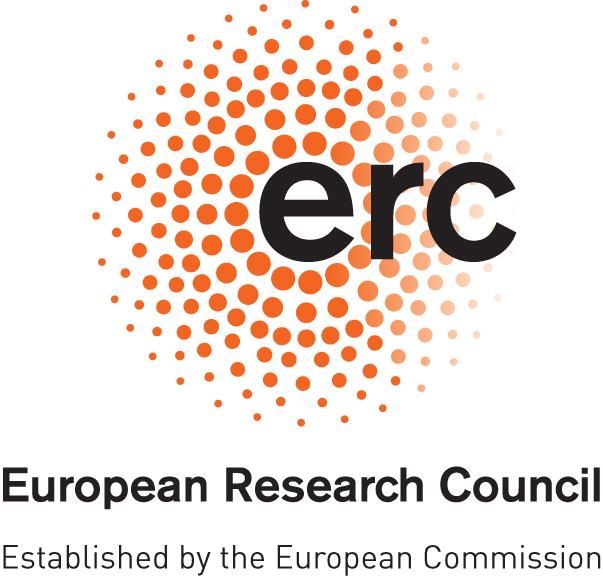}%
\end{textblock}
\begin{textblock}{20}(-0.25, 12.1)
\includegraphics[width=60px]{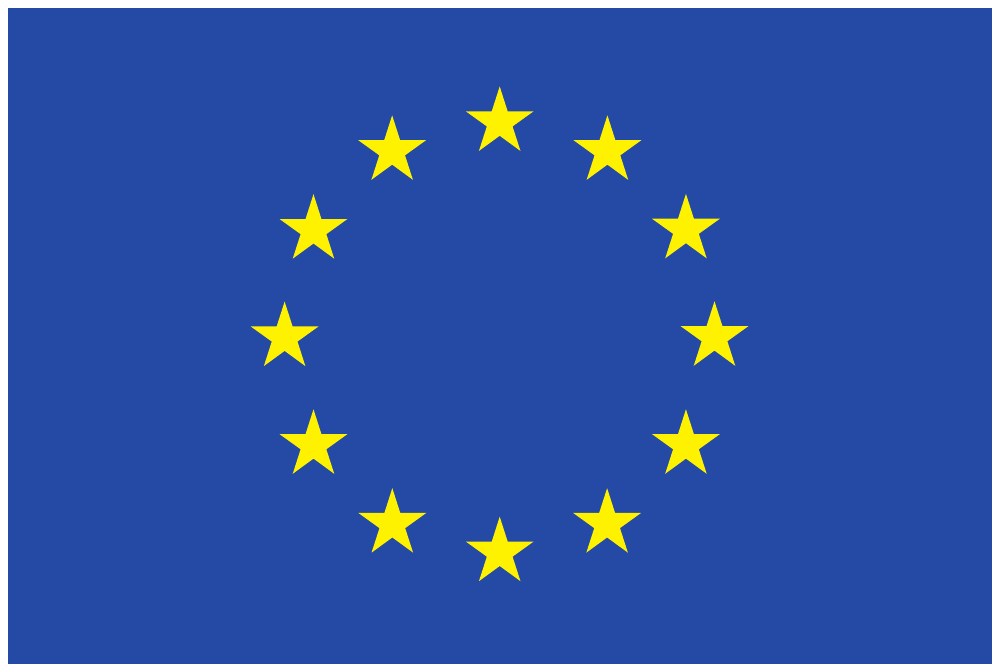}%
\end{textblock}
\begin{abstract}
Half graphs and their variants, such as ladders, semi-ladders and co-matchings, are combinatorial objects that encode total orders in graphs. Works by Adler and Adler~(Eur. J. Comb.; 2014) and Fabiański et~al.~(STACS; 2019) prove that in the powers of sparse graphs, one cannot find arbitrarily large objects of this kind. However, these proofs either are non-constructive, or provide only loose upper bounds on the orders of half graphs and semi-ladders.
In this work we provide nearly tight asymptotic lower and upper bounds on the maximum order of half graphs, parameterized on the distance, in the following classes of sparse graphs: planar graphs, graphs with bounded maximum degree, graphs with bounded pathwidth or treewidth, and graphs excluding a fixed clique as a minor.

The most significant part of our work is the upper bound for planar graphs. Here, we employ techniques of structural graph theory to analyze semi-ladders in planar graphs through the notion of cages, which expose a topological structure in semi-ladders.
As an essential building block of this proof, we also state and prove a new structural result, yielding a fully polynomial bound on the neighborhood complexity in the class of planar graphs.
\end{abstract}
\end{titlepage}

\section{Introduction}\label{introduction-chapter}  
It is widely known that there is a huge array of algorithmic problems deemed to be computationally hard.
One of the ways of circumventing this issue is limiting the set of possible instances of
  a~problem by assuming a more manageable structure.
For example, restricting our attention to graph problems,
  we can exploit the planarity of the graph instances through the planar separator
  theorem \cite{DBLP:journals/siamcomp/LiptonT80} or Baker's
  technique~\cite{DBLP:journals/jacm/Baker94}.
Analogously, if graphs have bounded treewidth or pathwidth, we
  can solve multiple hard problems by means of dynamic programming on tree
  or path decompositions~\cite{ParamAlgo}.
These examples present some of the algorithmic techniques which allow us to utilize the
  structural sparsity of graph instances.

Nešetřil and Ossona de Mendez have proposed two abstract notions of sparsity in
  graphs: bounded expansion \cite{DBLP:journals/endm/NesetrilM05}
  and nowhere denseness \cite{DBLP:journals/jsyml/NesetrilM10, DBLP:journals/ejc/NesetrilM11a}.
Intuitively speaking, a~class $\Cc$ of graphs has bounded expansion if one cannot obtain
  arbitrarily dense graphs by picking a~graph $G \in \Cc$ and contracting pairwise disjoint
  connected subgraphs of $G$ with fixed radius to single vertices.
More generally, $\Cc$ is nowhere dense if one cannot produce arbitrarily large cliques
  as graphs as a result of the same process; see Preliminaries
  (Section \ref{preliminaries-chapter}) for exact definitions.
These combinatorial notions turn out to be practical in the algorithm design.
In fact, they are some of the fundamental concepts of Sparsity --- a research area
  concerning the combinatorial properties and the algorithmic applications of sparse graphs.
The tools of Sparsity allow us to design efficient algorithms in nowhere dense classes of graphs
  and in classes of graphs with bounded expansion for problems such as subgraph isomorphism
  \cite{DBLP:conf/stoc/NesetrilM06, DBLP:journals/ejc/NesetrilM11a}
  or distance-$d$ dominating set \cite{DBLP:journals/endm/Nesetril08}.
In fact, under reasonable complexity assumptions, nowhere denseness exactly
  characterizes the classes of graphs closed under taking subgraphs
  that allow efficient algorithms for all problems definable in first-order logic
  \cite{DBLP:journals/jacm/DvorakKT13, DBLP:journals/jacm/GroheKS17}.
For a comprehensive introduction to Sparsity, we refer to the book by
  Nešetřil and Ossona~de~Mendez~\cite{DBLP:books/daglib/0030491} and
  to the lecture notes from the University of Warsaw \cite{SparsityUWNotes}.

Research in Sparsity provided a plethora of technical tools for structural analysis of sparse graphs,
  many of which are in the form of various graph parameters.
For example, the classes of graphs with bounded expansion equivalently have bounded
  \textbf{generalized coloring numbers} --- the \textbf{weak} and \textbf{strong $d$-coloring
  numbers} are uniformly bounded in a class $\Cc$ of graphs for each fixed $d \in \mathbb{N}$,
  if and only if $\Cc$ has bounded expansion \cite{DBLP:journals/dm/Zhu09}.
Another example is the concept of \textbf{$p$-centered colorings}
  \cite{DBLP:journals/ejc/NesetrilM06}, which have been shown to require
  a~bounded number of colors for each fixed $p \in \N$ exactly for the classes of graphs
  with bounded expansion \cite{DBLP:journals/ejc/NesetrilM08}.
It is interesting and useful to determine the asymptotic growth (depending on the parameter;
  $d$ or $p$ in the examples above) of these measures in well-studied 
  classes of sparse graphs,
  such as planar graphs, graphs with bounded pathwidth or treewidth, proper minor-closed classes
  of graphs etc., as it should lead to a better understanding of these classes.
For instance, the asymptotic behavior of generalized coloring numbers has been studied
  in proper minor-closed classes of graphs and in planar graphs
  \cite{DBLP:journals/siamdm/GroheKRSS18, DBLP:journals/endm/HeuvelMRS15}.
Also, there has been a recent progress on $p$-centered colorings
  in proper minor-closed classes, in planar graphs, and in graphs with bounded maximum degree
  \cite{DBLP:conf/soda/DebskiFMS20, DBLP:journals/corr/abs-2006-04113,
  DBLP:conf/soda/PilipczukS19}.

In this work, we will consider another kind of structural measures that behave nicely in
  sparse graphs, which are related to the concepts of
  half graphs, ladders, semi-ladders, and co-matchings.

\begin{definition}
\label{ladders-def-distone}
In an undirected graph $G=(V,E)$, for an integer $\ell \geq 1$, $2\ell$ different vertices
  $a_1, a_2, \dots, a_\ell$, $b_1, b_2, \dots, b_\ell$ form:
\begin{itemize}
  \item a \textbf{half graph} (or a \textbf{ladder})
    of order $\ell$ if for each pair of indices $i, j$ such that
    $i, j \in [1, \ell]$, we have $(b_i, a_j) \in E$ if and only if $i < j$
    (Figure \ref{ladders-fig}(a));
  \item a \textbf{semi-ladder} of order $\ell$ if we have $(b_i, a_j) \in E$
    for each pair of indices $i, j$ such that $1 \leq i < j \leq \ell$, and $(b_i, a_j) \notin E$
    for each $i \in [1, \ell]$ (Figure \ref{ladders-fig}(b));
  \item a \textbf{co-matching} of order $\ell$ if for each pair of indices $i, j$
    such that $i, j \in [1, \ell]$, we have $(b_i, a_j) \in E$ if and only if $i \neq j$
    (Figure \ref{ladders-fig}(c)).
\end{itemize}
\begin{figure}[h]
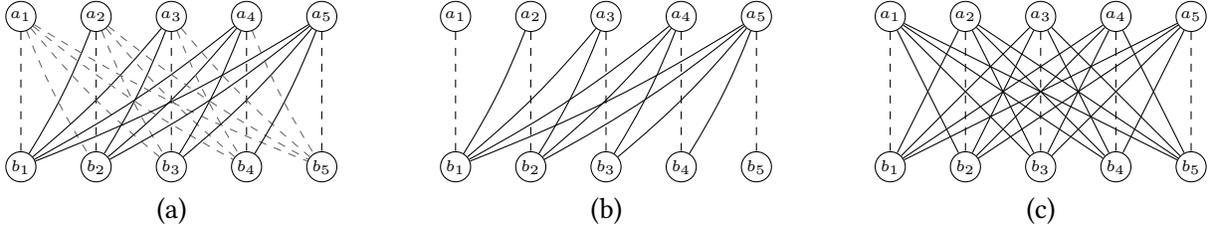

\centering
  \begin{minipage}[t]{0.30\textwidth}
	\centering
	\input{figures/semi-ladder-b.tex}
    (a)
    \end{minipage}\begin{minipage}[t]{0.05\textwidth}\ 
    \end{minipage}\begin{minipage}[t]{0.30\textwidth}
	\centering  
    \input{figures/semi-ladder-a.tex}
    (b)
    \end{minipage}\begin{minipage}[t]{0.05\textwidth}\ 
    \end{minipage}\begin{minipage}[t]{0.30\textwidth}
	\centering    
    \input{figures/semi-ladder-c.tex}
    (c)
    \end{minipage}
  \caption{The objects in Definition \ref{ladders-def-distone}. Solid lines indicate pairs of vertices
    connected by an edge, and dashed edges indicate pairs of vertices
    not connected by an edge. Intuitively, half graphs and co-matchings encode the
    $<$, and $\neq$ relations, respectively, on pairs
    $(a_1, b_1), (a_2, b_2), \dots, (a_\ell, b_\ell)$ of vertices.}
  \label{ladders-fig}
\end{figure}
\end{definition}

Naturally, each half graph, ladder, and co-matching is also a semi-ladder.
Hence, if for a class $\Cc$ of graphs, the orders of semi-ladders occurring in any
  graph in $\Cc$ are bounded from above by some constant $M$,
  then $M$ is also the corresponding upper bound for orders of half graphs, ladders and
  co-matchings.
We also remark that there is a simple proof utilizing Ramsey's theorem
  \cite[Lemma 2]{fabiaski2018progressive} which demonstrates that if an arbitrary class $\Cc$
  of graphs has a uniform upper bound on the orders of half graphs in any graph
  in~$\Cc$, and a separate uniform upper bound on the orders of co-matchings,
  then there also exists a uniform upper bound on the maximum order of
  semi-ladders in~$\Cc$.

In this work, we will consider $\Cc$ to be powers of nowhere dense classes of graphs.
Formally, for an~undirected graph $G$ and an integer $d$, we define the graph $G^d$ as follows:
$$ V(G^d) = V(G)\text{ \ \ and \ \ } E(G^d) = \{(u, v)\,\mid\,\dist_G(u, v) \leq d\}.$$
Then, the $d$-th power of a class $\Cc$ is defined as $\Cc^d = \{G^d\,\mid\,G\in \Cc\}$.

Even though $\Cc^d$ may potentially contain dense graphs, it turns out that the objects
  from Definition~\ref{ladders-def-distone} still behave nicely in $\Cc^d$:

\begin{theorem}[\cite{DBLP:journals/ejc/AdlerA14, fabiaski2018progressive}]
  \label{semi-ladders-bounded-thm}
  For every nowhere dense class $\Cc$ of graphs, there exists a function $f: \N \to \N$ such
  that for every positive integer $d$, the orders of semi-ladders in graphs from $\Cc^d$ are
  uniformly bounded from above by~$f(d)$.
\end{theorem}

Theorem \ref{semi-ladders-bounded-thm} implies that the orders of half graphs,
  ladders and co-matchings found in graphs of $\Cc^d$ are uniformly bounded by $f(d)$ as well.

In fact, a much more general result holds --- each nowhere dense class of graphs $\Cc$
  is \textbf{stable}.
That is, any first-order interpretation of $\Cc$ has a uniform bound on the maximum order of half graphs.
Formally, for a fixed first-order formula $\varphi(\bar{x}, \bar{y})$, where $\bar{x}$ and $\bar{y}$
  are tuples of variables, and a graph $G$, we define the~bipartite graph $G^\varphi = (V_1, V_2, E)$
  as follows:
$$V_1(G^\varphi) = V(G)^{\bar{x}}, \quad V_2(G^\varphi) = V(G)^{\bar{y}}\text{ \ \ and \ \ } 
  E(G^\varphi) = \{(\bar{a}, \bar{b})\,\mid\, G \models \varphi(\bar{a}, \bar{b})\}.$$
Then, the class $\Cc^\varphi = \{G^\varphi\,\mid\, G \in \Cc\}$ of graphs
  has a uniform bound on the maximum order of half graphs, as long as $\Cc$ is nowhere dense.
This has been proved by Adler and Adler~\cite{DBLP:journals/ejc/AdlerA14}, who applied the ideas
  of Podewski and Ziegler \cite{Podewski1978}.
While the proof of Adler and Adler relies on the compactness theorem of first-order logic,
  Pilipczuk, Siebertz, and~Toruńczyk proposed a~direct proof of this fact, giving explicit upper
  bounds on the orders of half graphs \cite{DBLP:conf/lics/PilipczukST18a}.
We stress that the generalization to the first-order logic is only valid for half graphs (it is
  straightforward to find counterexamples for semi-ladders and co-matchings), while
  in the specific case where $|\bar{x}| = |\bar{y}| = 1$ and $\varphi$ describes the distance-$d$
  relation between the pairs of vertices, even semi-ladders are bounded
  (Theorem~\ref{semi-ladders-bounded-thm}).

Stability is one of important properties of structures in mathematical logic, and has been
  studied extensively --- see the book by Pillay \cite{DBLP:books/daglib/0072694} for the general
  exposition of the topic, and the work of Pilipczuk et al. \cite{DBLP:conf/lics/PilipczukST18a}
  for the information on connections of the stability theory with sparse graphs.
As an example use of stability in the graph theory, one can consider stable
  classes of graphs with some additional structural properties.
For instance, stable classes of graphs of bounded
  cliquewidth~\cite{DBLP:journals/dam/CourcelleO00} or linear
  cliquewidth~\cite{DBLP:journals/jct/LozinR07} can, in fact, be constructed
  from the classes of graphs with bounded treewidth and pathwidth, respectively,
  by means of some first-order transformations
  \cite{DBLP:journals/corr/abs-2007-07857, DBLP:journals/ejc/NesetrilMRS21}.
Also, in the classes of graphs with a uniform bound on the orders of half graphs,
  Malliaris and Shelah improved the bounds on the parameters
  in the statement of the Szemer\'{e}di's regularity lemma \cite{10.2307/23813167}.
  
\bigskip

An important motivation for studying semi-ladders is a recent work of Fabiański et
  al.~\cite{fabiaski2018progressive}, which presented a fixed-parameter
  algorithm solving the \textsc{distance-$d$ dominating set} problem in sparse graphs:
  
\smallskip

\textsc{Fixed:} a nowhere dense class $\Cc$ of graphs, distance $d \geq 1$

\textsc{Input:} a graph $G \in \Cc$, a parameter $k \geq 1$

\textsc{Output:} a subset $D \subseteq V(G)$, $|D| \leq k$, such that each
  vertex of $G$ is at distance at most~$d$ from a~vertex of $D$, provided such a set exists.

\smallskip

In any nowhere dense class $\Cc$, the running time of this algorithm is bounded from above by
  $2^{O(k \log k)} \cdot (n + m)$ for each fixed $d \in \mathbb{N}$.
If $\Cc$ has bounded expansion, we can infer a~more concrete upper bound:
  $O\left((\xi \cdot k^L)^k + k^L \cdot (n+m)\right)$, where $L$ is the bound on the order of
  semi-ladders in graphs from $\Cc^d$, and $\xi$~is a constant
  related to the concept of neighborhood complexity of $\Cc$:
\begin{theorem}[\cite{DBLP:journals/ejc/ReidlVS19}]
  \label{neigh-complex-general-thm}
  Let $\Cc$ be a class of graphs with bounded expansion, and let $d \geq 1$ be an integer.
  For any graph $G$ and its vertex $v$, let $N^d[v]$ be the closed $d$-neighborhood of $v$;
    that is, the set of vertices at distance at most $d$ from $v$.
  Then, for every non-empty set $X$ of vertices of $G$, the number of different intersections
    of $X$ with a closed $d$-neighborhood of any vertex of $G$ is linear in the size of $X$:
  \begin{equation}
  \label{neigh-complex-eq}
  \left|\left\{N^d[v] \cap X\,\mid\,v \in V(G)\right\}\right| \leq \xi \cdot |X|
  \end{equation}
  where $\xi$ is a constant depending only on $\Cc$ and $d$.
  The left-hand side of (\ref{neigh-complex-eq}) is called the \textbf{neighborhood complexity}
    of $G$.
\end{theorem}
If $\Cc$ is nowhere dense, the runtime guarantees of this algorithm are slightly looser, but they still rely
  heavily on the value of $L$.
Hence, if we find tight asymptotic bounds on $L$ as a function of $d$ in the $d$-th powers of
  various concrete classes
  of sparse graphs, we should be able to infer concrete bounds on the running time of the algorithm
  of Fabiański et al. in these classes.

As mentioned above, the value of $L$ is bounded in powers of nowhere dense classes of graphs.
However, the proof in \cite{fabiaski2018progressive} relies on the description of nowhere
  denseness through uniform quasi-wideness (see Definition \ref{uniform-quasi-wideness-def}).
Hence, the resulting bounds on the maximum order of semi-ladders are usually
  suboptimal --- for instance, applying them directly to the $d$-th power of the
  class of planar graphs yields a $d^{O(d^6)}$ upper bound;
  compare this to  a $d^{O(d)}$ bound we present in this work (Theorem \ref{planar-upper-bound}).
Therefore, it is natural to pose a question: can we find good estimates on the maximum orders
  of semi-ladders occurring in the powers of well-studied classes of sparse graphs?

\vspace{1em}
\noindent \textbf{Our results.} We state and prove asymptotic lower and upper bounds on the
  orders of semi-ladders in the $d$-th powers of the following classes of sparse graphs:
  graphs with maximum degree bounded by $\Delta$, planar graphs, graphs with pathwidth
  bounded by $p$, graphs with treewidth bounded by $t$, and graphs excluding the
  complete graph $K_t$ as a minor (Figure \ref{results-figure}).
Our results are asymptotically almost tight --- they allow us to understand how quickly
  the maximum orders of semi-ladders grow, as functions of $d$,
  in these classes of graphs.

\begin{figure}[h]
\centering
\begin{tabular}{c c c}
\textbf{\ \ \ \ Class of graphs\ \ \ \ } & \textbf{\ \ Lower bound\ \ } & \textbf{\ \ Upper
  bound\ \ } \\ \hline
    \begin{tabular}{@{}c@{}} \\ \textsc{Degree} $\leq\Delta$ \\ \ \end{tabular} &
    \begin{tabular}{@{}c@{}} $\Delta^{\Omega(d)}$ \\ (Cor. \ref{degree-lower-bound-asym})
      \end{tabular} &
    \begin{tabular}{@{}c@{}} $\Delta^d + 1$ \\ (Thm. \ref{degree-upper-bound})
      \end{tabular} \\ \hline

    \begin{tabular}{@{}c@{}} \\ \textsc{Planar} \\ \ \end{tabular} &
    \begin{tabular}{@{}c@{}} $2^{\left\lceil \frac{d}{2} \right\rceil}$ \\
      (Cor. \ref{planar-lower-bound}) \end{tabular} &
    \begin{tabular}{@{}c@{}} $d^{O(d)}$ \\ (Thm. \ref{planar-upper-bound})
      \end{tabular} \\ \hline
    
    \begin{tabular}{@{}c@{}} \\ \textsc{Pathwidth} $\leq p$ \\ \ \end{tabular} &
    \begin{tabular}{@{}c@{}} $d^{p - O(1)}$ \\ (Cor. \ref{pw-lower-bound-asym})
      \end{tabular} &
    \begin{tabular}{@{}c@{}} $(dp)^{O(p)}$ \\ (Cor. \ref{pw-upper-bound-asym})
      \end{tabular}  \\ \hline
    
    \begin{tabular}{@{}c@{}} \\ \textsc{Treewidth} $\leq t$ \\ \ \end{tabular} &
    \begin{tabular}{@{}c@{}} $2^{\displaystyle d^{\Omega(t)}}$ \\
      (Cor. \ref{tw-lower-bound-asym}) \end{tabular} &
    \begin{tabular}{@{}c@{}} $d^{\displaystyle O(d^{t + 1})}$ \\
      (Cor. \ref{tw-upper-bound-asym}) \end{tabular} \\ \hline
    
    \begin{tabular}{@{}c@{}} \\ $K_t$-\textsc{minor-free} \\ \ \end{tabular} &
    \begin{tabular}{@{}c@{}} $2^{\displaystyle d^{\Omega(t)}}$ \\
      (Cor. \ref{kt-lower-bound-asym}) \end{tabular} &
    \begin{tabular}{@{}c@{}} $d^{\displaystyle O(d^{t - 1})}$ \\
      (Cor. \ref{kt-upper-bound-asym}) \end{tabular} \\
\end{tabular}
\caption{Bounds on the maximum semi-ladders in the $d$-th powers of classes of graphs
  proved in this work.}
\label{results-figure}
\end{figure}

In Section \ref{lower-bounds-chapter}, we describe three constructions proving
  the lower bounds stated in Figure~\ref{results-figure}.
Each of the constructions, for an integer $d$, creates a graph whose $d$-th power contains
  a large half graph.
Hence, these graphs also contain large semi-ladders.
Conversely, the remaining part of this work (Sections \ref{noose-profile-lemma-section},
  \ref{planar-upper-bound-section}, \ref{other-upper-bounds-chapter})
  contains proofs of the upper bounds on the maximum orders
  of semi-ladders in the considered classes of graphs; these proofs also yield the upper bounds
  on the orders of half graphs and co-matchings.

In Section \ref{noose-profile-lemma-section}, we prove that the neighborhood complexity
  in planar graphs has an upper bound which is polynomial both in $|X|$, and in $d$.
We note that Theorem \ref{neigh-complex-general-thm} already proves a bound linear
  in $|X|$; however, the proof of Reidl et al.~\cite{DBLP:journals/ejc/ReidlVS19} requires
  an exponential dependence on $d$, even in the class of planar graphs.
The polynomial dependence on~$d$ asserted by our result is crucial in the proof of the upper
  bound on semi-ladders in the class of planar graphs.

Section \ref{planar-upper-bound-section} is fully devoted to the proof of the upper bound
  on the maximum order of a semi-ladder in the $d$-th power of any planar graph.
The proof is quite involved and should be considered the main contribution of this work.
This proof works as follows:
  we start with a huge semi-ladder in the $d$-th power of a chosen planar graph $G$.
Using this semi-ladder, we find in $G$ objects with more and more structure, which we call:
  quasi-cages, cages, ordered cages, identity ordered cages, neighbor cages,
  and separating cages, in this order.
Each extraction step requires us to forfeit a fraction of the object, but
  we give good estimates on the maximum size loss incurred in the process.
Eventually, we derive a direct upper bound on the maximum order of a separating cage that
  can be found in a planar graph.
Retracing the extraction steps, we find a uniform upper bound on the maximum
  order of any semi-ladder in the $d$-th power of a planar graph.

Section \ref{other-upper-bounds-chapter} presents the proofs of upper bounds in the powers of
  other considered classes of graphs: graphs with bounded degree, graphs with bounded
  pathwidth or treewidth, and $K_t$-minor-free graphs.
In the case of graphs of bounded pathwidth, we generalize the sunflower lemma
  (Lemma \ref{large-sunflower-lemma}), which ultimately helps us
  to extract structural properties from path decompositions of bounded width.
In the proof of the upper bound in the class of graphs excluding $K_t$ as a minor,
  we rely on weak coloring numbers, and upper bounds on them in $K_t$-minor-free graphs
  proved by van den Heuvel et al. \cite{DBLP:journals/endm/HeuvelMRS15}.

\vspace{1em}
\noindent \textbf{Acknowledgements.}
  We would like to thank Michał Pilipczuk for pointing us to this problem, the verification
    of the correctness of the proofs in this work, and proofreading this work.
  
  Also, we want to thank Piotr Micek for reviewing this work and valuable suggestions
    regarding Section~\ref{noose-profile-lemma-section}.

\section{Preliminaries}\label{preliminaries-chapter}
\noindent
\textbf{Distances.} In this work --- unless specified --- we only consider undirected graphs,
  without loops or multiple edges with the same endpoints.
  In a graph $G$, for two of its vertices $u, v$ (not necessarily different), we denote the
  distance between these two vertices in $G$ as $\dist_G(u, v)$; if these vertices are not in the
  same connected component of $G$, we set $\dist_G(u, v) = +\infty$.

  When the graph is clear from the context, we may omit the name of the graph from our
    notation and simply write $\dist(u, v)$.
    
\vspace{1em}
\noindent \textbf{Different definitions of half graphs.}
In some applications (e.g. in the definition of the stability in the Introduction),
  half graphs are allowed to contain multiple occurrences of any single vertex.
However, a straightforward greedy procedure can convert this variant of half graphs into the half graphs from
  Definition~\ref{ladders-def-distone} of order at least half of the original order.
Since we will only consider an asymptotic behavior of the orders of half graphs and other structures
  in various classes of graphs, we can safely assume the variant of half graphs from
  Definition~\ref{ladders-def-distone}.

\vspace{1em}
\noindent \textbf{Distance-$d$ half graphs (ladders, semi-ladders, co-matchings).}
In many parts of this work, it will be more useful to operate on the original graph $G$ instead of
  its $d$-th power $G^d$ --- some of the objects and decompositions defined later are properly
  defined in the original classes of graphs, but not in their powers.
Hence, we define \textbf{distance-$d$} half-graphs, ladders, semi-ladders, and co-matchings
  that exist in $G$ if and only if their counterparts from Definition~\ref{ladders-def-distone}
  exist in $G^d$.

In an undirected graph $G=(V,E)$, for two integers $d, \ell \geq 1$, $2\ell$ different vertices
  $a_1, a_2, \dots, a_\ell$, $b_1, b_2, \dots, b_\ell$ form:
\begin{itemize}
  \item a \textbf{distance-$d$ half graph} (or a \textbf{distance-$d$ ladder})
    of order $\ell$ if for each pair of indices $i, j$ such that
    $i, j \in [1, \ell]$, we have $\dist(b_i, a_j) \leq d$ if and only if $i \leq j$;
  \item a \textbf{distance-$d$ semi-ladder} of order $\ell$ if we have $\dist(b_i, a_j) \leq d$
    for each pair of indices $i, j$ such that $1 \leq i < j \leq \ell$, and $\dist(b_i, a_j) > d$
    for each $i \in [1, \ell]$;
  \item a \textbf{distance-$d$ co-matching} of order $\ell$ if for each pair of indices $i, j$
    such that $i, j \in [1, \ell]$, we have $\dist(b_i, a_j) \leq d$ if and only if $i \neq j$.
\end{itemize}

\vspace{1em}
\noindent \textbf{Pathwidth, path decompositions.}
We define a path decomposition $(W_1, W_2, \dots, W_k)$ of a graph $G$
as any sequence of subsets of $V(G)$ (called \textbf{bags}) with the following properties:
  \begin{itemize}
  \item For each edge $uv \in E(G)$, there exists a bag $W_i$ such that $u, v \in W_i$.
  \item For each vertex $v \in V(G)$, $v$ appears in a contiguous subset of the sequence.
  \end{itemize}
The width of the decomposition is $\max_{i=1}^k |W_i| - 1$. The \textbf{pathwidth} of a graph
  $G$, denoted $\pw{G}$, is the minimum possible width of any path decomposition of $G$.

\vspace{1em}
\noindent \textbf{Treewidth, tree decompositions.}
A tree decomposition of a graph $G$ is an undirected tree $\Tc$ whose each vertex $t$
  is associated with a subset of vertices $X_t$, called the \textbf{bag} of $t$.
The tree has the following properties:
  \begin{itemize}
  \item For each edge $uv \in E(G)$, there exists a vertex $t \in V(\Tc)$ such that $u, v \in X_t$.
  \item For each vertex $v \in V(G)$, the set of vertices $\{t\,\mid\,v \in X_t\}$
    forms a connected subtree of~$\Tc$.
  \end{itemize}
Similarly to the case of path decomposition, the width of a tree decomposition is the maximum
  size of any bag in the decomposition, minus one.
The \textbf{treewidth} of a graph $G$, denoted $\tw(G)$, is the minimum possible width
  of any tree decomposition of $G$.

\vspace{1em}
\noindent \textbf{Graph minors, minor-free classes of graphs.}
A graph $H$ is a \textbf{minor} of $G$ if one can map each vertex $v \in V(H)$ to
  a non-empty connected subgraph $A(v)$ of $G$
  so that every two subgraphs are vertex-disjoint, and for each edge $uv \in E(H)$,
  subgraphs $A(u)$ and $A(v)$ are connected by an edge.
A graph $H$ is a \textbf{depth-$d$ minor} of $G$ if additionally each subgraph $A(v)$ has
  radius not exceeding $d$.

A class $\Cc$ of graphs is $H$-minor-free if no graph $G \in \Cc$ contains $H$ as a minor.
We note that for any class $\Cc$ and graph $H$, if $\Cc$ is $K_{|H|}$-minor-free, then it is also
  $H$-minor-free.
Hence, the maximum order of any semi-ladder in an $H$-minor-free class of graphs is
  bounded from above by the upper bound for $K_{|H|}$-minor-free graphs.

As treewidth is a minor-closed parameter, and $\tw(K_{t+2}) = t + 1$, the following holds:
\begin{theorem}
  \label{tw-minor-free}
  Every graph $G$ such that $\tw(G) \leq t$ is also $K_{t+2}$-minor-free.
\end{theorem}

\vspace{1em}
\noindent \textbf{Bounded expansion.}
We follow the notation used in the lecture notes from the University of
  Warsaw~\cite{SparsityUWNotes}.

For a graph $G$ and an integer $d \in \N$, we define $\nabla_d(G)$ to be the maximum edge
  density of any depth-$d$ minor of $G$:
$$ \nabla_d(G) := \sup\left\{\frac{|E(H)|}{|V(H)|}\,\mid\,H\text{ is a depth-}d\text{ minor of }
  G\right\}. $$
A class $\Cc$ of graphs has \textbf{bounded expansion} if there exists a function $f : \N \to \N$
  such that we have $\nabla_d(G) \leq f(d)$ for each $G \in \Cc$, $d \in \N$.

\vspace{1em}
\noindent \textbf{Nowhere denseness.}
For a graph $G$ and an integer $d \in \N$, we define $\omega_d(G)$ to be the maximum size
  of a clique which is a depth-$d$ minor of $G$:
$$ \omega_d(G) := \sup\left\{t\,\mid\,K_t\text{ is a depth-}d\text{ minor of }G \right\}. $$
A class $\Cc$ of graphs is \textbf{nowhere dense} if there exists a function $f : \N \to \N$
  such that we have $\omega_d(G) \leq f(d)$ for each $g \in \Cc$, $d \in \N$.

It can be easily verified that each class $\Cc$ with bounded expansion is also nowhere dense.
However, the converse statement is not true: there are nowhere dense classes which
  have unbounded expansion.

\vspace{1em}
\noindent \textbf{Oriented paths.}
In an undirected graph $G$, an oriented path is a sequence of vertices
  $P = (v_0, v_1, v_2, \dots, v_\ell)$ such that $(v_{i-1}, v_i) \in E(G)$ for each $i \in \{1, 2, \dots,
  \ell\}$.
The path is simple if all vertices $v_0, v_1, v_2, \dots, v_\ell$ are different.

For an oriented simple path $P$ and two vertices $u, v \in P$ with $u$ occurring
  on the path before $v$, we define $P[u, v]$ as the subpath of $P$ starting in $u$ ending in $v$.
Next, for any two oriented paths $P = (p_0, p_1, \dots, p_x)$, $Q = (q_0, q_1, \dots, q_y)$ such
  that $p_x = q_0$, we denote by $P \cdot Q = (p_0, p_1, \dots, p_x, q_1, q_2, \dots, q_y)$
  the concatenation of these two paths.
Also, for any oriented path $P = (v_0, v_1, \dots, v_\ell)$, we define
  $P^{-1} = (v_\ell, v_{\ell-1}, \dots, v_0)$ as the version of $P$ with flipped orientation,
  and $\len(P) = \ell$ as its length.

\vspace{1em}
\noindent \textbf{Degenerate graphs.}
A $k$-degenerate graph ($k \geq 1$) \cite{lick_white_1970} is a graph in which every
  subgraph contains a vertex of degree at most $k$.
For instance, each planar graph is $5$-degenerate.

In the thesis, will use one property of degenerate graphs --- each $k$-degenerate graph
  is $(k+1)$-colorable~\cite{SZEKERES19681}.
  
We remark that the classes of degenerate graph are not necessarily nowhere dense
  (e.g. each subdivided clique is a $2$-degenerate graph).
Consequently, the class of $2$-degenerate graphs contains arbitrary large semi-ladders since
  subdivided half graphs are $2$-degenerate.

\vspace{1em}
\noindent \textbf{Distance-$d$ profiles.}
For a vertex $v \in V(G)$ and a subset of vertices $S \subseteq V(G)$,
  a {\bf distance\nobreakdash-$d$ profile} of $v$ on $A$ in $G$ is a function $\pi_d[v, S]$
  from $S$ to
  $\{0, 1, 2, \dots, d, +\infty\}$ that assigns each vertex from $S$ either its distance from $v$,
  or $+\infty$ if this distance exceeds $d$.

\section{Lower bounds on the orders of half graphs}\label{lower-bounds-chapter}

In this section, we present a series of constructions which provide us with graphs belonging
  to the considered classes and containing large half graph.
Since each half graph is also a semi-ladder,
  the constructed graphs will also contain large semi-ladders.

The section is divided into several subsections.
In Subsection \ref{degree-section}, we produce a family of graphs, each having bounded degree,
  but containing sizeable half graphs.
It will later turn out that a subfamily of these graphs --- specifically, those graphs in the family
  whose maximum degrees are bounded from above by $4$ --- consists solely
  of planar graphs.
These graphs also witness an exponential lower bound on the maximum half graph size
  in the class of planar graphs.
  
In Subsection \ref{pw-section}, we produce graphs with small pathwidth, yet large half graphs.
Despite the construction being rather simple --- a ``core'' containing $\Theta(p)$ vertices
  where $p$ is the bound on pathwidth, together with a huge number of vertices connected
  only to the vertices of the core --- it gives us a polynomial lower bound on
  the possible co-matching order, with the degree of the polynomial depending linearly on $p$.
  
In Subsection \ref{tw-section}, we construct graphs with small treewidth and large half graphs.
The construction is quite non-trivial --- each graph in the family
  will be formed from multiple copies of smaller graphs in this family.
This dependence allows us to create large half graphs in this class of graphs by means of induction.
As the graphs with bounded treewidth do not contain large cliques as minors (Theorem
  \ref{tw-minor-free}), we will immediately find large half graphs in $K_t$-minor-free graphs.

\subsection{Graphs with bounded degree and planar graphs}
\label{degree-section}

In this subsection, we will prove the following theorem:

\begin{theorem}
\label{degree-lower-bound}
For each $\Delta \geq 4$ and odd $d \geq 1$,
  there exists a graph with maximum degree bounded by $\Delta$
  containing a distance-$d$ half graph of order $\left\lfloor \frac{\Delta}{2}\right\rfloor^{
  \left\lceil \frac{d}{2} \right\rceil}$.
\end{theorem}
  
We fix two integer variables, $k \geq 2$ and $h \geq 1$.
We will construct an undirected graph $H_{k,h}$ with the following properties:
  \begin{itemize}
    \item The degree of each vertex in the graph is bounded by $2k$.
    \item The graph contains vertices $a_1, b_1, a_2, b_2, \dots, a_{k^h}, b_{k^h}$
      forming a distance-$(2h-1)$ half graph of order~$k^h$.
  \end{itemize}

In the construction of $H_{k,h}$, we will create two isomorphic rooted trees
  $\Ac$ and $\Bc$, each having $k^h$ leaves, and connect pairs of vertices of
  these trees using paths of small, positive length.
  
We fix an alphabet $\Sigma = \{1, 2, \dots, k\}$ with the natural ordering of its elements.
This allows us to introduce a lexicographical ordering of all words in $\Sigma^*$,
  denoted $\preccurlyeq$.
For a word $s \in \Sigma^*$, by $|s|$ we mean the length of $s$.
And for a pair of words $x, y \in \Sigma^*$, by $x \cdot y$ we mean the concatenation of
  $x$ and $y$.

For each word $s \in \Sigma^*$, $|s| \in [0, h]$, we create two vertices: $a_s \in V(\Ac)$
  and $b_s \in V(\Bc)$.
The edges in both trees are constructed in the following way:
  for each word $s \in \Sigma^*$, $|s| \leq h-1$, and each character $c \in \Sigma$,
  we add an edge between vertices $a_s$ and $a_{s \cdot c}$ in $\Ac$, and between
  the vertices $b_s$ and $b_{s \cdot c}$ in $\Bc$.

We also connect the trees in the following way:
  for each word $s \in \Sigma^*$, $|s| \leq h-1$, and each pair of characters $c, d \in \Sigma$
  such that $c < d$, we connect the vertices $a_{s \cdot d} \in V(\Ac)$ and $b_{s \cdot c} \in V(\Bc)$ by
  a path of length $2|s| + 1$.
  
\begin{figure}[h]
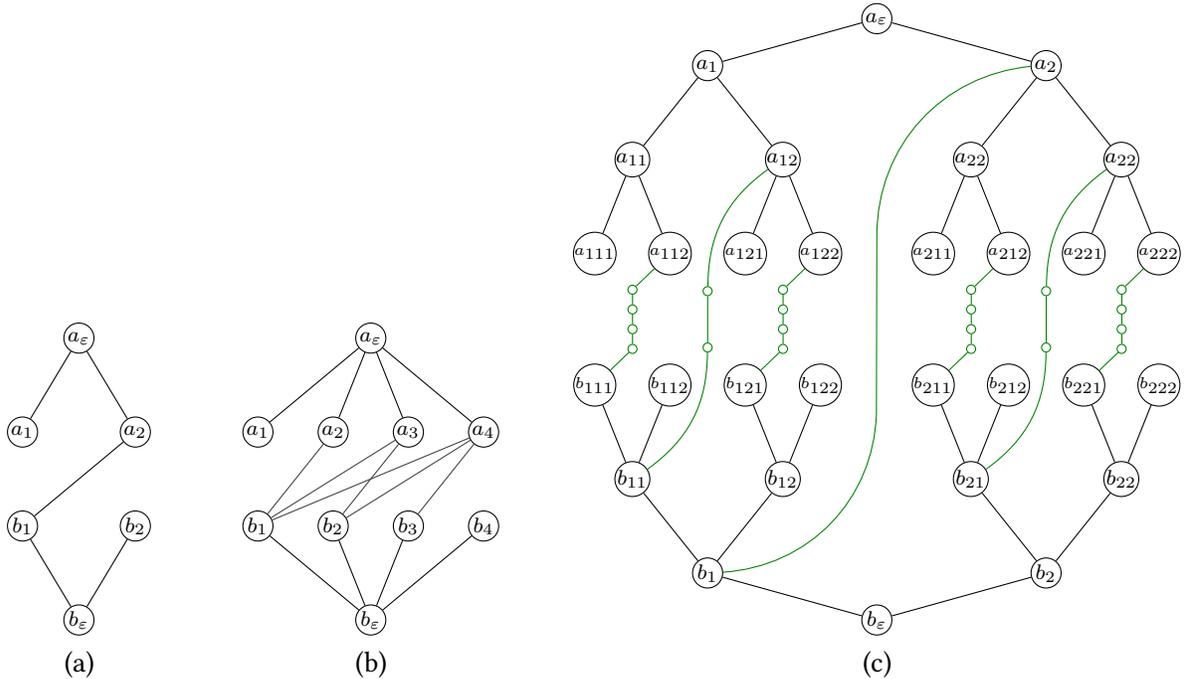

  \centering
  \begin{minipage}[b]{0.185\textwidth}
      \centering
      \input{figures/degree-a.tex}
      (a)
    \end{minipage}\begin{minipage}[b]{0.285\textwidth}
      \centering
      \input{figures/degree-b.tex}
      (b)
    \end{minipage}\begin{minipage}[b]{0.53\textwidth}
      \centering
      \input{figures/degree-c.tex}
      (c)
    \end{minipage}
  \caption{Graphs: (a) --- $H_{2,1}$, (b) --- $H_{4,1}$, (c) --- $H_{2,3}$. $\varepsilon$
    denotes an empty word.}
  \label{degree-examples-fig}
\end{figure}

This finishes the construction of $H_{k,h}$.
Examples of such graphs can be found in Figure~\ref{degree-examples-fig}.

\vspace{0.5em}
We now need to prove the required properties of $H_{k, h}$.
Firstly, we prove that each vertex has at most $2k$ neighbors in the graph.
Vertices belonging to neither $\Ac$ nor $\Bc$ lie on the paths connecting $\Ac$ with $\Bc$,
  so their degrees are equal to $2$.
Let us now fix a vertex belonging to $\Ac$ or $\Bc$.
In its own tree, this vertex has at most $k$ children and one parent.
This vertex is also an endpoint of at most $k-1$ paths connecting $\Ac$ with $\Bc$;
  the bound is satisfied with equality for vertices $a_{s \cdot k}$ and $b_{s \cdot 1}$
  ($1, k \in \Sigma$, $s \in \Sigma^*$, $|s| \leq h-1$).
Hence, the total number of neighbors of each vertex of $\Ac$ and $\Bc$ is bounded by $2k$.
  
Secondly, each tree obviously has $k^h$ leaves, since the leaves are parameterized by a word
  of length $h$ over an alphabet of size $k$.

Finally, we need to find a half graph in this graph. Let $(t_1, t_2, \dots, t_{k^h})$ be the sequence
  of all words of length $h$ over $\Sigma$, ordered lexicographically.
The following lemma proves that the lexicographical ordering of the leaves of $\Ac$ and $\Bc$
  produces a half graph:
\begin{lemma}
  \label{lemma-degree-distance}
  Fix two indices $i, j \in [1, k^h]$.
  The distance between $b_{t_i}$ and $a_{t_j}$ is equal to $2h-1$ if $i < j$, and is greater
    than $2h - 1$ if $i \geq j$.
    
  \begin{proof}
  For simplicity, we can assume that each path of length $\ell$ ($\ell \geq 1$)
    connecting a vertex of $\Ac$ with a vertex of $\Bc$ is, in fact, a weighted edge of length
    $\ell$ connecting these two vertices.

  We introduce the \textbf{level function} $\xi\,:\,V(\Ac) \cup V(\Bc) \to \mathbb{N}$,
    defined in the following way:
    
  $$ \xi(v) = \begin{cases}
    h - |s| & \text{if }v = a_s\text{ for }s \in \Sigma^*, \\
    h + |s| - 1 & \text{if }v = b_s\text{ for }s \in \Sigma^*.
  \end{cases} $$
  
  For a vertex $v$, we say that the value $\xi(v)$ is the \textbf{level} of vertex $v$.
  
  We see that the levels of the vertices of $\Ac$ range between
    $0$ and $h$, while the levels of the vertices of $\Bc$ range between $h-1$ and $2h - 1$
    (Figure \ref{degree-levels-fig}).
    
  \begin{figure}[h]
    \centering
    \input{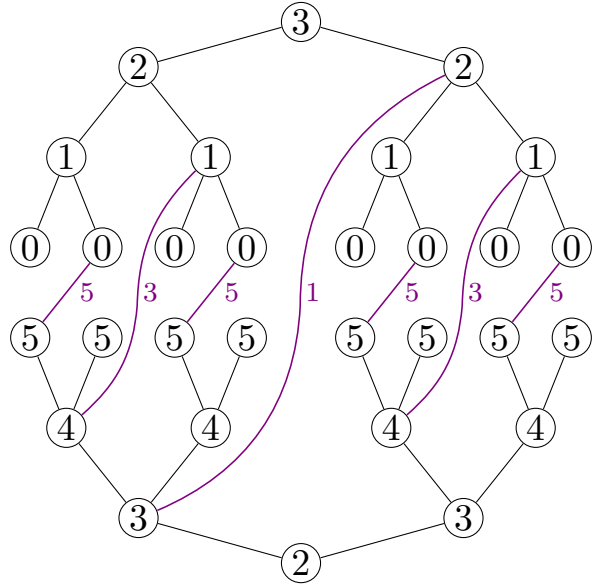}
    \label{degree-levels-fig}
    \caption{The levels of vertices of $H_{2,3}$. Each weighted edge resulting from
      a contracted path is colored violet and is labeled by its length.}
  \end{figure}

  We can now observe that:
  \begin{itemize}
    \item For every two vertices of $\Ac$ connected by an edge --- $a_s$ and $a_{s \cdot c}$
      for some $s \in \Sigma^*$, $c \in \Sigma$ --- we have that $\xi(a_{s}) = \xi(a_{s \cdot c}) + 1$.
    \item For every two vertices of $\Bc$ connected by an edge --- $b_s$ and $b_{s \cdot c}$ ---
      we have that $\xi(b_{s \cdot c}) = \xi(b_s) + 1$.
    \item For every two vertices $a_{s \cdot d} \in V(\Ac)$, $b_{s \cdot c} \in V(\Bc)$
      ($s \in \Sigma^*, c, d \in \Sigma, c < d$), connected by a weighted edge of length
      $2|s| + 1$, we have
      $$\xi(b_{s \cdot c}) - \xi(a_{s \cdot d}) = (h + |s|) - (h - |s| - 1) = 2|s| + 1.$$
  \end{itemize}
  Hence, for each edge of weight $x$ connecting two vertices $u, v$, we have that
    $x = |\xi(u) - \xi(v)|$.
  Therefore, for all pairs of vertices $u, v$, the distance between these vertices
    is bounded from below by $|\xi(u) - \xi(v)|$, and the bound is satisfied with equality only if
    there exists a path between $u$ and $v$ on which the values of $\xi$ change
    monotonically.
  We will call such paths \textbf{monotonous}.
    
  Fix $i, j \in [1, k^h]$, and consider the vertices $b_{t_i}$ and $a_{t_j}$.
  By the definition of $\xi$, we know that $\xi(a_{t_j}) = 0$ and $\xi(b_{t_i}) = 2h - 1$,
    hence the distance between these two vertices is at least $2h - 1$.
  We will now analyze when this distance is equal to $2h - 1$.

  We see that the value of $\xi$ only increases when: going up the tree $\Ac$ (in the direction
    of the root of $\Ac$), taking a weighted edge from a vertex of $\Ac$ to a vertex of
    $\Bc$ (but not the other way around), and going down the tree $\Bc$ (away from the root
    of $\Bc$).
  Hence, each path between $b_{t_i}$ and $a_{t_j}$
    of length $2h-1$ can only have such three segments, in this exact order, as it needs to
    be monotonous.
  
  If $i < j$, then $t_i \prec t_j$, so the words
    $t_i$ and $t_j$ admit the following factorization:
  $$ t_i = p \cdot c \cdot s_i, \qquad t_j = p \cdot d \cdot s_j $$
  for $p, s_i, s_j \in \Sigma^*$, $c, d \in \Sigma$ and $c < d$.

  Therefore, the path originating at $a_{t_j}$, going up $\Ac$ to the vertex $a_{p \cdot d}$, taking
    the weighted edge to $b_{p \cdot c}$ connecting both trees (this is possible since $c < d$),
    and then going down $\Bc$ to the destination $b_{t_i}$, is monotonous.
  Hence, $\dist(b_{t_i}, a_{t_j}) = 2h - 1$.
  
  Conversely, if a path from $a_{t_j}$ to $b_{t_i}$ has length $2h - 1$, then it is monotonous.
  Therefore, it must follow from $a_{t_j}$ to an ancestor of $a_{t_j}$ (we will name it
    $a_x$ where $x$ is a prefix of $t_j$), then take a weighted
    edge to an ancestor of $b_{t_i}$ (we will call it $a_y$ where $y$ is a prefix of $t_i$),
    and then walk down the tree $\Bc$ to $b_{t_i}$.
  By the construction of the tree, $x$ and $y$ must be both of the form
    $x = p \cdot d$, $y = p \cdot c$ where $p \in \Sigma^*$, $c, d \in \Sigma$, $d > c$.
  Since $x$ and $y$ are prefixes of $t_j$ and $t_i$, respectively, this means
    that $t_j \succ t_i$.
  Hence, $j > i$.
  
  Therefore, the path of length $2h - 1$ exists between $a_{t_j}$ and $b_{t_i}$ if and only if
    $i < j$; otherwise the shortest path between these vertices is longer.  
  \end{proof}
\end{lemma}

Lemma \ref{lemma-degree-distance} proves that the graph $H_{k, h}$ contains
  a half graph of order $k^h$.
Since $k \geq 2$ and $h \geq 1$ were two arbitrary variables, this statement is satisfied
  for each choice of these variable.
  
\vspace{0.5em}
Let us now finish the proof of Theorem \ref{degree-lower-bound}.
For each $\Delta \geq 4$ and odd $d \geq 1$ we can see that the graph 
  $H_{\left\lfloor\frac{\Delta}{2}\right\rfloor, \frac{d+1}{2}}$ has its maximum degree
  bounded by $2\left\lfloor\frac{\Delta}{2}\right\rfloor \leq \Delta$,
  and contains a distance-$d$ half graph of order $\left\lfloor\frac{\Delta}{2}\right\rfloor^{
  \left\lceil\frac{d}{2}\right\rceil}$.
Hence, the proof of the theorem is complete.

\medskip

We note the following immediate conclusion from Theorem \ref{degree-lower-bound}:
\begin{corollary}
  \label{degree-lower-bound-asym}
  Let $\Cc_\Delta$ be the class of graphs with maximum degree bounded by $\Delta \geq 4$.
  Then the maximum order of a distance-$d$ half graph ($d \geq 1$) in $\Cc_\Delta$ is bounded
    from below by $\Delta^{\Omega(d)}$.
\end{corollary}

It turns out that the presented family of examples also witness as an exponential lower
  bound on the maximum order of a distance-$d$ half graph in the class of planar graphs.
This is proved explained by following lemma:

\begin{lemma}
  \label{degree-four-is-planar-lemma}
  For each $h \geq 1$, the graph $H_{2, h}$ is planar.
  
  \begin{proof}
  For each $h \geq 1$, let $H'_{2, h}$ be the graph obtained from $H_{2, h}$ by contracting
    each path connecting two trees to a single edge.
  Obviously, $H'_{2, h}$ is planar if and only if $H_{2, h}$ is planar.
  Then, for $h \geq 2$, the graph $H'_{2, h}$ can be created from two copies of the graph
    $H'_{2, h-1}$ in the following way:
  \begin{enumerate}
    \item In the $i$-th copy of the graph ($i \in \{1, 2\}$), replace each vertex $a_s$
      ($s \in \{1, 2\}^*$) by $a_{i \cdot s}$, and analogously each vertex $b_s$ by $b_{i \cdot s}$.
    \item Create new roots of $\Ac$ and $\Bc$ --- $a_\varepsilon$ and $b_\varepsilon$,
      respectively.
    \item Add the following edges: $a_\varepsilon a_1$, $a_\varepsilon a_2$,
      $b_\varepsilon b_1$, $b_\varepsilon b_2$, $b_1 a_2$.
  \end{enumerate}
  
  \begin{figure}[h]
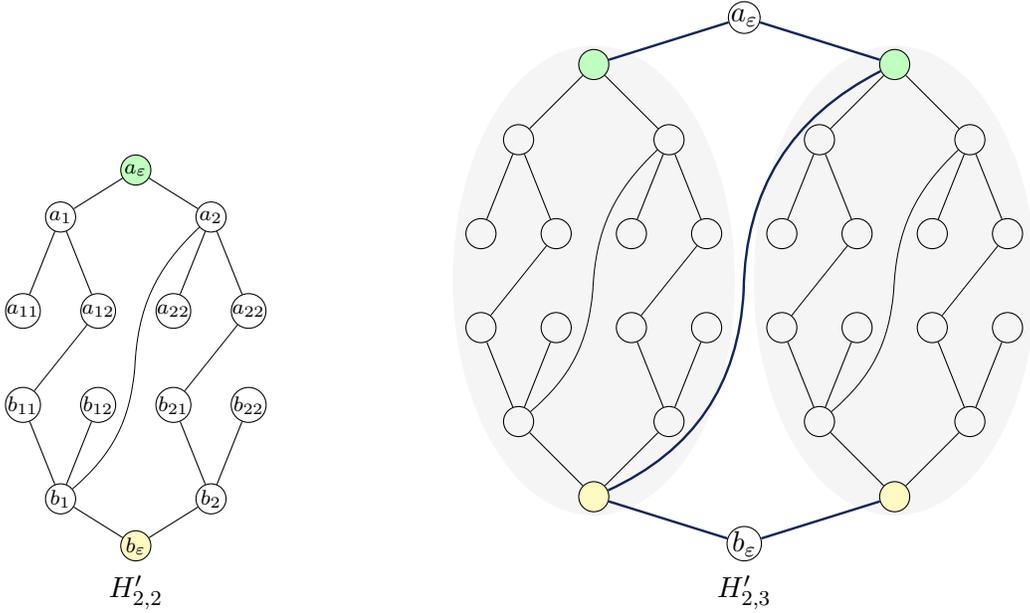

  \centering
  \begin{minipage}[b]{0.35\textwidth}
      \centering
      \input{figures/planar-conn-a.tex}
      $H'_{2,2}$
    \end{minipage}\begin{minipage}[b]{0.63\textwidth}
      \centering
      \input{figures/planar-conn-b.tex}
      $H'_{2,3}$
    \end{minipage}
  \caption{A step in the inductive process of creating $H'_{2,h}$ for each $h \geq 2$.
    Two copies of $H'_{2,2}$ (left; its roots are marked green and yellow) are put side by side,
    two new vertices are added, and five new edges (blue) are drawn, resulting in the graph
    $H'_{2,3}$ (right).}
  \label{planar-conn-fig}
  \end{figure}
  
  Note that $H'_{2, 1} = H_{2, 1}$ is planar (Figure \ref{degree-examples-fig}(a)), and moreover,
    in the presented embedding, its vertices $a_\varepsilon$ and $b_\varepsilon$ lie on the
    outer face of the graph.
  These two properties can be maintained inductively:
    for $h \geq 2$, we take the embeddings of two copies of $H'_{2, h-1}$ satisfying these
    properties, put them side by side,
    and add the required edges, making sure that the graph remains planar and
    that the new roots of the trees are on the outer face of the graph (Figure \ref{planar-conn-fig}).
  \end{proof}
\end{lemma}

We immediately infer the following:
\begin{corollary}
  \label{planar-lower-bound}
  For every $d \geq 1$, the maximum order of a distance-$d$ half graph ($d \geq 1$) 
  in the class of planar graphs is bounded from below by $2^{\left\lceil \frac{d}{2} \right\rceil}$.
  \begin{proof}
  Assume $d$ is odd.
  Then, by Lemma \ref{degree-four-is-planar-lemma}, the graph $H_{2, \frac{d+1}{2}}$ is planar.
  Also, it contains a distance-$d$ half graph of order
    $2^{\frac{d+1}{2}} = 2^{\left\lceil \frac{d}{2} \right\rceil}$.

  For even $d$, we construct a planar graph containing a distance-$(d-1)$ half graph
    of order $2^{\frac{d}{2}} = 2^{\left\lceil \frac{d}{2} \right\rceil}$, and then for
    each vertex $a_i$ in this graph ($1 \leq i \leq 2^{\frac{d}{2}}$), we add an edge whose
    one endpoint is $a_i$, and another is a new vertex $a'_i$.
  We now can see that the vertices $a'_1, b_1, a'_2, b_2, \dots, a'_{2^{\frac{d}{2}}},
    b_{2^{\frac{d}{2}}}$ form a distance-$d$ half graph of required order.
  \end{proof}
\end{corollary}

\subsection{Graphs with bounded pathwidth}
\label{pw-section}

In this section, we will find large half graphs in the class of graphs with bounded pathwidth.
Specifically, we will prove the following theorem:

\begin{theorem}
  \label{pw-lower-bound}
  For each $d \geq 1$ and $p \geq 0$, there exists a graph $P_{p,d}$ with pathwidth
    at most $p + 2$ which contains a distance-$(4d-1)$ half graph of order
    $(2d+1)^p$.
\begin{proof}

Let $\ell_{p,d} := (2d+1)^p$.
We define a slightly broader class of graphs $\{P_{p,d,k}\,\mid\,d \geq 1,\,
  p \geq 0,\,k\geq0\}$
  where $P_{p,d,k}$ has the following properties:

\begin{itemize}
  \item The graph has pathwidth at most $p + 2$.
  \item The graph contains a distance-$(4d-1)$ half graph
    $a_1, a_2, \dots, a_{\ell_{p,d}}$, $b_1, b_2,
    \dots, b_{\ell_{p, d}}$ in which every two vertices are at distance
    $2d$ or more from each other.
  \item Each of the vertices of the half graph
    has $k$ vertex disjoint paths of length $3d$ attached to it (called
    \emph{appendices});
    in other words, after removing any of these vertices, there are at least $k$
    connected components which are paths with $3d$ vertices.
\end{itemize}

Thanks to the first two conditions, the graph $P_{p, d} := P_{p, d, 0}$ will satisfy
  Theorem~\ref{pw-lower-bound}.
The final condition is purely technical and will help us argue that pathwidths
  of the graphs constructed in the process described below do not grow too large.
  
\medskip

The construction will be inductive over $p$.
First, for $d \geq 1$, $k \geq 0$, we choose the graph $P_{0, d, k}$
  to be the $3d$-subdivision of two disjoint stars with $k$ edges, with the
  distance-$(4d-1)$ half graph of order $\ell_{0,d} = 1$ given by the centers of both stars.
It can be easily verified that this graph satisfies the requirements; in particular,
  $\pw{P_{0,d,k}} \leq 2$.

\medskip

We shall now describe how to construct $P_{p, d, k}$ for $p, d \geq 1$ and $k \geq 0$:
\begin{enumerate}
  \item Take $2d+1$ disjoint copies of $P_{p-1, d, k+1}$, and renumber the vertices
    of the half graph in each of the copies, ordering all $2d+1$ half graphs one after another.
    Formally, the labels of the vertices
    $a_i$ and $b_i$ of the half graph in the $j$-th disjoint copy of the graph are changed to
    $a_{(j-1)\ell_{p-1,d} + i}$ and $b_{(j-1)\ell_{p-1,d} + i}$, respectively.
    Since $\ell_{p,d} = (2d+1)\ell_{p-1,d}$,
      the graph contains vertices $a_i$, $b_i$ for each $i \in [1, \ell_{p, d}]$.
    The following steps of the construction will ensure that these vertices form
      a distance-$(4d-1)$ half graph, provided that $P_{p-1,d,k+1}$ satisfies the required properties.
    
  \item Add a new vertex $r$ to the graph.
  
  \item For each vertex $b_i$ in the $j$-th disjoint copy of the graph, $j \in [1, 2d+1]$,
      we connect $b_i$ with $r$ with a~path of length $d + (j - 1)$.
    In order to do that, we use up one of $k+1$ appendices attached to $b_i$
      by connecting $r$ with the $(d+j-2)$-nd vertex of the appendix with an edge
      (Figure \ref{pw-appendix-join});
      the part of the appendix after this vertex can then be safely discarded.
  
  \item We apply a similar procedure to the vertices $a_i$ of the half graph, only that
    the vertex of the half graph in the $j$-th disjoint copy of the graph is connected
    with $r$ with a~path of length $3d - (j - 1)$.
    This finishes the construction of $P_{p, d, k}$.
\end{enumerate}

\begin{figure}[h]
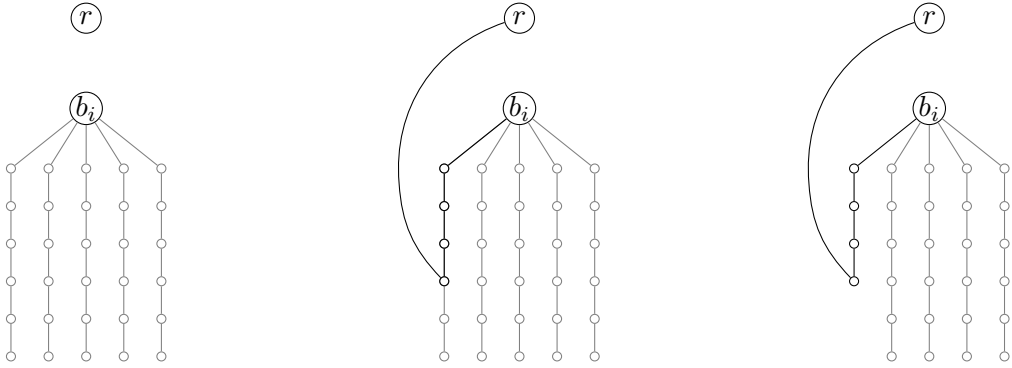

  \centering
  \begin{minipage}[b]{0.33\textwidth}
      \centering
      \input{figures/pw-appendix-join-a.tex}
  \end{minipage}\begin{minipage}[b]{0.33\textwidth}
      \centering
      \input{figures/pw-appendix-join-b.tex}
  \end{minipage}\begin{minipage}[b]{0.33\textwidth}
      \centering
      \input{figures/pw-appendix-join-c.tex}
  \end{minipage}
  \caption{The process of turning an appendix attached to $b_i$
    into a path of length $5$ connecting $b_i$ with $r$.}
  \label{pw-appendix-join}
\end{figure}

Obviously, each vertex $a_i$ and $b_j$ in $P_{p,d,k}$
  has $k$ remaining appendices.
We also note that $P_{p, d, k} \setminus \{r\}$ is a subgraph of a graph formed
  from $2d+1$ disjoint copies of $P_{p-1, d, k+1}$. Hence,
$$ \pw{P_{p, d, k}} \leq \pw{P_{p,d,k} \setminus \{r\}} + 1 \leq
  \pw{P_{p-1, d, k+1}} + 1 \leq p + 2. $$

We also argue that $a_1, a_2, \dots, a_{\ell_{p,d}}$, $b_1, b_2, \dots,
  b_{\ell_{p,d}}$ is a distance-$(4d-1)$ half graph satisfying the distance requirements
  of $P_{p,d,k}$.
In the construction, multiple disjoint copies of $P_{p-1,d,k+1}$ were connected
  by vertex disjoint paths to a new vertex $r$, where each path has length at least
  $d$ and at most $3d$.
Hence, the shortest path between any two vertices of the new half graph is either fully
  contained in a copy of $P_{p-1,d,k+1}$ (and it must not be shorter than $2d$
  by the inductive assumption), or it passes through $r$ and thus it follows two
  paths of length at least $d$, so its length is at least $2d$ as well.
From this property,
  we infer that if the shortest path between any two vertices of the half graph has length
  at most $4d - 1$, then this path must avoid all other vertices of the half graph.

Let us now consider whether there exists a path between two vertices $a_i$ and $b_j$
  of the half graph of length $4d - 1$.
Again, this path must either be fully contained in a single copy of
  $P_{p-1,d,k+1}$ (and this condition is only satisfied if $i > j$), or it passes through $r$.
In the latter case, this path must consist of a direct connection from $a_i$ to $r$,
  and then a direct connection from $r$ to $b_j$.
Assuming that $a_i$ is in the $x$-th disjoint copy of $P_{p-1,d,k+1}$, and $b_j$ is
  in the $y$-th copy of the graph ($x, y \in [1, 2d+1]$), we infer that the length of that path
  is exactly
  $$[d + (y - 1)] + [3d - (x - 1)] = 4d + (y - x).$$
Hence, such a path exists if and only if $x > y$; that is, only if $i > j$ and both
  $a_i$ and $b_j$ come from the different copies of $P_{p-1,d,k+1}$.

We conclude that $\dist(a_i, b_j) \leq 4d-1$ if and only if $i > j$.
Therefore, all the required conditions are satisfied, so $P_{p,d,k}$ is defined correctly.
This finishes the proof of the lower bound.
\end{proof}
\end{theorem}

On an end note, we remark the following asymptotic version of Theorem \ref{pw-lower-bound}:

\begin{corollary}
  \label{pw-lower-bound-asym}
  For each pair of integers $p \geq 2$, $d \geq 4$, the maximum order of a distance-$d$
    half graph in the class of graphs with pathwidth bounded by $p$ is bounded from below
    by $\Theta(d)^{p - O(1)}$.
\end{corollary}

\subsection{Graphs with bounded treewidth and minor-free graphs}
\label{tw-section}

In this section, we will construct a family of graphs with small treewidth containing large
  half graphs.
Since the graphs with bounded treewidth also avoid large cliques as minors, this
  family will also serve as examples of graphs avoiding relatively small cliques as minors,
  but still containing large half graphs.
  
To start with, we shall introduce a variant of a tree decomposition of a graph $G$
  containing a distance-$2d$ half graph.

\begin{definition}
  \label{pairing-td-def}
  A rooted tree decomposition $\Tc$ of a graph $G$ is
    a~\textbf{distance-$d$ pairing decomposition} if all the following conditions hold:
  \begin{itemize}
    \item The root of $\Tc$ is a bag containing exactly two elements, called $A$ and
      $B$.
    \item Each leaf of $\Tc$ is a bag containing exactly two elements, called $a_i$ and
      $b_i$, where $i \in \{1, 2, \dots, \ell\}$ and $\ell$ is the number of leaves in $\Tc$.
    \item All vertices $A, B, a_1, b_1, a_2, b_2, \dots, a_\ell, b_\ell$ are pairwise different.
    \item For each pair of integers $i, j$ such that $1 \leq i < j \leq \ell$, we have that
      $\dist(b_i, a_j) = 2d$.
    \item For each pair of integers $i, j$ such that $1 \leq j \leq i \leq \ell$, we have that
      $\dist(b_i, a_j) = 2d + 1$.  
  \end{itemize}
  For such a tree decomposition, we define two objects: its root
    $\Root(\Tc) = (A, B)$ and the sequence of vertices in its leaves
    $\Leaves(\Tc) = (a_1, b_1, a_2, b_2, \dots, a_\ell, b_\ell)$.
  We will refer to the value $\ell$ as the \textbf{order} of $\Tc$.
  We also define the \textbf{width} of the tree decomposition in the standard way.
\end{definition}

We remark that Definition \ref{pairing-td-def} implies that
  the set of vertices $a_1, a_2, \dots, a_\ell, b_1, b_2, \dots, b_\ell$ forms
  a distance-$2d$ half graph of order $\ell$ in $G$.
However, the conditions required here are even stronger
  as the definition of a half graph only requires that $\dist(b_i, a_j) \leq 2d$ for each pair of integers
  $i < j$, and $\dist(b_i, a_j) > 2d$ for $i > j$.
This strengthening will be essential in the proof of the correctness
  of the construction.

We also note that $\Tc$ --- in some way ---
  ``exposes'' the vertices of the underlying half graph
  in the form of the leaf bags of the tree decomposition ($\Leaves(\Tc)$),
  and two other vertices in the form of the root bag of the decomposition ($\Root(\Tc)$).
These outer bags of the tree decomposition are ``sockets'' which
  can connect to the ``sockets'' of pairing decompositions of other graphs, allowing
  a~recursive construction of large graphs.
  
\begin{example}
  \label{pairing-td-example}
  We consider the graph $G$ presented in Figure \ref{pairing-td-fig}(a)
    with the tree decomposition $\Tc$ depicted in Figure \ref{pairing-td-fig}(b).
  This decomposition is distance-$1$ pairing and has order~$2$ and width~$3$.
  The table in Figure \ref{pairing-td-fig}(c) shows the distances between
    $a_i$ and $b_j$ for each pair of indices $i, j$.
  We can see that $\dist(b_1, a_2) = 2$ and all the remaining distances are equal to~$3$.
  Moreover, $\Root(\Tc) = (A, B)$ and $\Leaves(\Tc) = (a_1, b_1, a_2, b_2)$.  
  \begin{figure}[h]
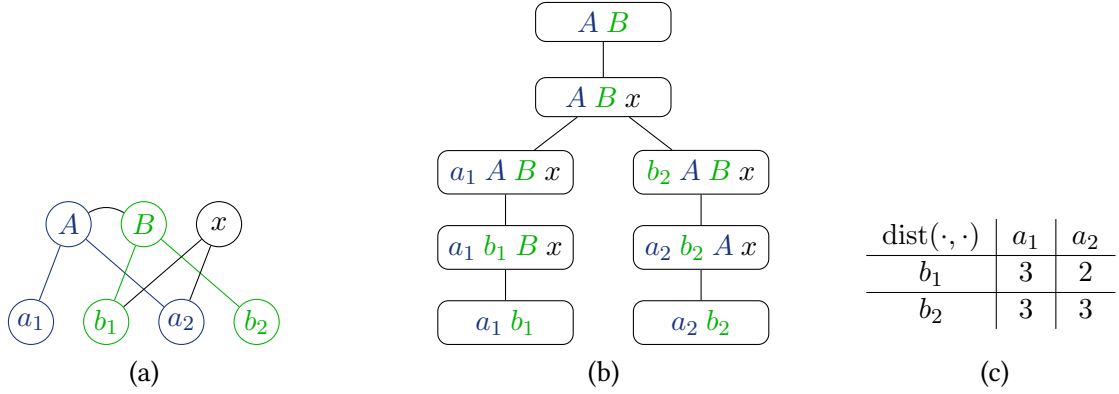

    \centering
    \begin{minipage}[b]{0.37\textwidth}
      \centering
      \input{figures/tw-basic-a.tex}
      (a)
    \end{minipage}\begin{minipage}[b]{0.37\textwidth}
      \centering
      \input{figures/tw-basic-b.tex}
      (b)
    \end{minipage}\begin{minipage}[b]{0.26\textwidth}
      \centering
      \begin{tabular}{c|c|c}
        $\dist(\cdot, \cdot)$ & $a_1$ & $a_2$ \\ \hline
        $b_1$ & $3$ & $2$ \\ \hline
        $b_2$ & $3$ & $3$
      \end{tabular}
      \vspace{1em}
      
      (c)
    \end{minipage}
    \caption{(a) --- the graph $G$ in Example \ref{pairing-td-example};
      (b) --- its tree decomposition; (c) --- its distance matrix.}
      \label{pairing-td-fig}
  \end{figure}
\end{example}

In the construction of our family of graphs, we will need to use a refinement
  of Definition~\ref{pairing-td-def}.
Intuitively, we want to expose the root vertices of the decomposition in an efficient way.
We will do it by creating a stronger variant of the pairing decomposition in which
  one of the root vertices is adjacent to all the vertices on one side of the half graph, while
  the other root vertex is adjacent to the other side.
  
\begin{definition}
  Consider an arbitrary undirected graph $G$ with a distance-$d$ pairing decomposition
    $\Tc$ where $\Root(\Tc) = (A, B)$ and $\Leaves(\Tc) = (a_1, b_1, \dots, a_\ell, b_\ell)$.
  We will say that $\Tc$ is \textbf{neighboring} if for each $i \in \{1, 2, \dots, \ell\}$, we have that:
      $$ \dist(a_i, A) = \dist(b_i, B) = 1, \quad\ \ \dist(a_i, B) \geq 2d + 1, \quad\ \
        \dist(b_i, A) \geq 2d + 1. $$
  Moreover, $A$ must be adjacent only
    to vertices $a_1, \dots, a_\ell$, and $B$ must be adjacent only to vertices
    $b_1, \dots, b_\ell$.
\end{definition}

\vspace{0.5em}

We can verify that the decomposition from Example \ref{pairing-td-example} is not neighboring
  since the vertices $A$ and $B$ are connected by an edge.
  
\vspace{0.5em}

We will now present a few lemmas.
Lemma \ref{make-neighboring-lemma} will allow us to turn an arbitrary decomposition into
  a neighboring decomposition by incurring a small cost on the width of the decomposition.
Lemma \ref{tw-combine-lemma} will be a method of combining two pairing decompositions ---
  a neighboring decomposition of width $t$ and order $\ell_1$, and a pairing decomposition
  of width $t-2$ and order $\ell_2$, into a pairing decomposition of order $\ell_1\ell_2$.
Then, the combination of these two lemmas (Lemma \ref{combine-any-tw}) will provide us
  with a procedure, which we will use to create graphs with enormous pairing decompositions.

\begin{lemma}
  \label{make-neighboring-lemma}
  There exists a procedure which takes an arbitrary undirected graph $G$ with
    a distance-$d$ pairing decomposition $\Tc$ of width $t$ and order $\ell$,
    and produces a graph $G'$ with a neighboring distance-$d$
    pairing decomposition $\Tc'$ of width $t+2$ and order $\ell$.
  \begin{proof}
    Let $\Leaves(\Tc) = (a_1, b_1, \dots, a_\ell, b_\ell)$.

    We construct $G'$ from $G$ by creating two new vertices $A_{\mathrm{new}}$ and
      $B_{\mathrm{new}}$, connecting $A_{\mathrm{new}}$ directly by an edge
      to each of the vertices $a_i$ for $i \in [1, \ell]$, and analogously connecting
      $B_{\mathrm{new}}$ to each of the vertices $b_i$.

    We will now construct a rooted tree decomposition $\Tc'$ from $\Tc$ in the following way
    \begin{enumerate}
    \item For each leaf bag $L_i \in V(\Tc)$ containing vertices $a_i$ and $b_i$,
      we clone it (without cloning the vertices in the bag), and connect it with $L_i$ in $\Tc$.
      This way, the clone becomes a leaf in the new decomposition instead of $L_i$.
    \item We create a new bag $(A_{\mathrm{new}}, B_{\mathrm{new}})$, which we directly
      connect to the root of $\Tc$.
    \item To each bag of $\Tc$ which is neither the root nor a leaf of the new tree decomposition,
      we add both $A_{\mathrm{new}}$ and $B_{\mathrm{new}}$.
    \end{enumerate}
    
    We can easily see that $\Tc'$ is a tree decomposition of $G'$ of width $t + 2$.
    
    We will now show that $(G', \Lc')$ is a neighboring decomposition
      by deducing all the required equalities and inequalities between the distances.
    Obviously, $\dist(a_i, A_{\mathrm{new}}) = \dist(b_i, B_{\mathrm{new}}) = 1$ for each
      $i \in [1, \ell]$.
    Furthermore, by considering all neighbors of $B_{\mathrm{new}}$, we can see
      that for each $i \in [1, \ell]$, we have that
    $$ \dist(a_i, B_{\mathrm{new}}) =  \min_{j \in \{1, \dots, \ell\}} \left(
      \dist(a_i, b_j) + 1\right) \geq 2d + 1. $$
    Analogously, $\dist(b_i, A_{\mathrm{new}}) \geq 2d+1$ for every $i \in [1, \ell]$.

    We can now prove that no path connecting $b_i$ and $a_j$ which has length $2d$ or less,
      for any choice of indices $i, j \in [1, \ell]$, can pass through $A_{\mathrm{new}}$:
    $$ \dist(b_i, A_{\mathrm{new}}) \geq 2d + 1 \qquad\Rightarrow\qquad
      \dist(b_i, A_{\mathrm{new}}) + \dist(A_{\mathrm{new}}, a_j) > 2d + 1. $$
    Analogously, no such path can pass through $B_{\mathrm{new}}$.

    If $i < j$, then $\dist_G(b_i, a_j) = 2d$ holds in $G$.
    Since $A_{\mathrm{new}}$ and $B_{\mathrm{new}}$ were the only vertices added to
      $G$ in the construction of $G'$, and no path with length $2d$ or less can pass through
      these vertices, it means that we also have $\dist_{G'}(b_i, a_j) = 2d$.
    Analogously, if $i \geq j$, then $\dist_G(b_i, a_j) = 2d+1$ in also implies
    $\dist_{G'}(b_i, a_j) = 2d+1$.
  \end{proof}
\end{lemma}

Let the procedure introduced in the statement of Lemma \ref{make-neighboring-lemma}
  be named $\mathsf{MakeNeighboring}$.
The procedure takes a graph $G$ together with its pairing distance-$d$ tree
  decomposition $\Tc$, and produces a modified graph $G'$ together with its
  neighboring pairing distance-$d$ decomposition $\Tc'$ of width larger by $2$ than that
  of $\Tc$:
$$ (G', \Tc') = \mathsf{MakeNeighboring}(G, \Tc).$$

\begin{lemma}
  \label{tw-combine-lemma}
  For every $d \geq 2$, $t \geq 1$ and $\ell_1, \ell_2 \geq 2$, there exists a procedure which takes:
  \begin{itemize}
    \item a graph $G$ with a distance-$(d-1)$ pairing decomposition $\Tc_G$
      of order $\ell_1$ and width $t$,
    \item a graph $H$ with a neighboring distance-$d$ pairing decomposition $\Tc_H$
      of order $\ell_2$ and width $t$,
  \end{itemize}
  and produces a graph with a distance-$d$ pairing decomposition of order $\ell_1\ell_2$
  and width $t$.
  
  \begin{proof}
  We fix the graph $G$ with its pairing decomposition $\Tc_G$, and let
    $\Root(\Tc_G) = (A, B)$, $\Leaves(\Tc_G) = (A_1, B_1, A_2, B_2, \dots,
    A_{\ell_1}, B_{\ell_1})$.
  We also create $\ell_1$ disjoint copies of the graph $H$ together with its tree decomposition
  $\Tc_H$.
  For each $i \in [1, \ell_1]$, we define the graph $H_i$ as the $i$-th copy of $H$.
  In its tree decomposition $\Tc_{H_i}$, we name the vertices of the root bag and the leaf bags:
  \begin{equation*}
  \begin{split}
    \Root(\Tc_{H_i}) &= (\alpha_i, \beta_i), \\
    \quad \Leaves(\Tc_{H_i}) &= (a_{(i-1)\ell_2 + 1},
    b_{(i-1)\ell_2 + 1}, a_{(i-1)\ell_2 + 2}, b_{(i-1)\ell_2 + 2}, \dots, a_{i\ell_2}, b_{i\ell_2}). 
  \end{split}
  \end{equation*}

  We remark that the vertices $A_1, B_1, A_2, B_2, \dots, A_{\ell_1}, B_{\ell_1}$ form a half graph
    in $G$.
  Moreover, for each $i \in [1, \ell_1]$, the set of vertices $a_j, b_j$ for all $j \in [(i-1)\ell_2+1,
    i\ell_2]$ forms a half graph in $H_i$.
  We also note that in the sequence
    $(a_1, b_1, a_2, b_2, \dots, a_{\ell_1\ell_2}, b_{\ell_1\ell_2})$, the vertices of $H_1$ form a prefix,
    followed by the vertices of $H_2$, then $H_3$, and so on,
    and a suffix of this sequence is the vertices of $H_{\ell_1}$.
    
  We create a new undirected graph $U$, together with its tree decomposition $\Tc_U$,
    in the following way (Figure \ref{combine-tw-algo-fig}):
  \begin{enumerate}
    \item We take the disjoint union of the graphs $G, H_1, \dots, H_{\ell_1}$, together with
      the forest created by the disjoint union of the corresponding tree decompositions
      $\Tc_G, \Tc_{H_1}, \Tc_{H_2}, \dots, \Tc_{H_{\ell_1}}$.
    \item For each $i \in \{1, 2, \dots, \ell_1\}$, we fuse the vertices $A_i$ and $\alpha_i$, and
      the vertices $B_i$ and $\beta_i$.
    \item In the forest of tree decompositions,
      for each $i \in \{1, 2, \dots, \ell_1\}$, we fuse the $i$-th leaf bag of $\Tc_G$
      (containing vertices $A_i$ and $B_i$) with the root bag of $\Tc_{H_i}$ (containing
      $\alpha_i$ and $\beta_i$).
      This step is possible since $A_i = \alpha_i$ and $B_i = \beta_i$.
      We call the resulting tree $\Tc_U$.
  \end{enumerate}
  
  \begin{figure}[h]
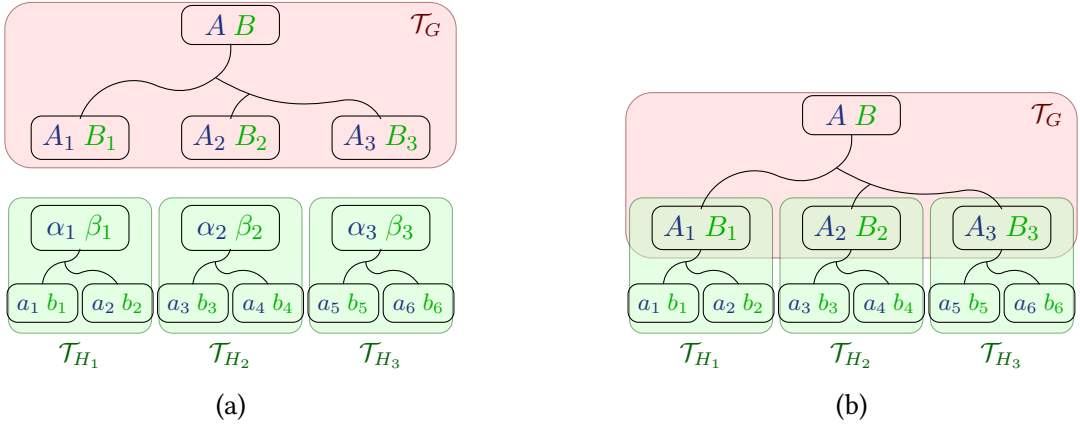

    \centering
    \begin{minipage}[b]{0.5\textwidth}
      \centering
      \input{figures/tw-combine-a.tex}
      (a)
    \end{minipage}\begin{minipage}[b]{0.5\textwidth}
      \centering
      \input{figures/tw-combine-b.tex}
      (b)
    \end{minipage}
    \caption{(a) --- tree decompositions after Step 1 of the algorithm;
      all the bags apart from the root and leaves have been omitted. \\      
      (b) --- the resulting tree decomposition $\Tc_U$.}
    \label{combine-tw-algo-fig}
  \end{figure}
  
  It is straightforward to see that $\Tc_U$ is a tree decomposition of $U$,
    and that its width remains equal to $t$.
  
  In $U$, we can identify the subgraphs $G, H_1, \dots, H_{\ell_1}$, originating from the disjoint
    union in Step 1.
  We stress that for each $i \in [1, \ell_1]$, subgraphs $G$ and $H_i$ share two vertices ---
    $A_i$ and $B_i$.
  These two vertices form a bag separating $\Tc_G$ from $\Tc_{H_i}$.
  This has a couple of consequences, directly implied by the properties of tree decompositions:
  
  \begin{itemize}
    \item If a simple path $P$ connects a vertex from $H_i$ with a vertex outside of $H_i$,
      then $P$ must pass through one of the vertices $A_i$, $B_i$ (possibly both);
    \item If a simple path $P$ has both endpoints in $H_i$, but contains a vertex outside of $H_i$,
      then $P$ must pass through both $A_i$ and $B_i$.
  \end{itemize}
  
  We now proceed to proving all the distance equalities and inequalities required by the
    definition of a distance-$d$ pairing decomposition.
    
  \begin{claim}
    \label{tw-claim-bottom-long-path}
    For each $k \in [1, \ell_1]$, the vertices $A_k$ and $B_k$ are at distance $2d+2$
      in $H_k$.
    \begin{claimproof}
      \cqed
      The set of neighbors of $A_k$ in $H_k$ is exactly
        $\{a_j\,\mid\, j \in [(k-1)\ell_2+1,\,k\ell_2]\}$, while
        the set of neighbors of $B_k$ in $H_k$ is
        $\{b_i\,\mid\, i \in [(k-1)\ell_2+1,\,k\ell_2]\}$.
      Hence, the second vertex on the shortest path from $A_k$ to $B_k$ in $H_k$
        is $a_j$ for some $j \in [(k-1)\ell_2+1,\,k\ell_2]$,
        while the penultimate vertex on this path is $b_i$
        for some $i \in [(k-1)\ell_2+1,\,k\ell_2]$.

      Since the vertices $a_i, b_i$ for $i \in [(k-1)\ell_2+1,\,k\ell_2]$ are in the distance-$d$ pairing
        decomposition $\Tc_{H_k}$, we get that $\dist_{H_k}(b_i, a_j) = 2d$
        if $i < j$; otherwise, this distance is larger than $2d$.
      As $\ell_2 \geq 2$, the vertices $b_{k\ell_2-1}$ and $a_{k\ell_2}$ both belong to
        $H_k$, and the distance between them is equal to $2d$.
      Hence,
      \[ \dist_{H_k}(A_k, B_k) = 2 + \min_{i, j \in [(k-1)\ell_2+1,\,k\ell_2]} \dist_{H_k} (b_i, a_j) = 2 +
        \dist_{H_k}(b_{k\ell_2 - 1}, a_{k\ell_2}) = 2d + 2. \]

    \end{claimproof}
  \end{claim}
  
  \begin{claim}
    \label{tw-claim-top-same-dists}
    For each pair of vertices $u, v \in V(G)$, we have that $\dist_U(u, v) = \dist_G(u, v)$.
    \begin{claimproof}
      \cqed
      We first note that for each $k \in [1, \ell_1]$, the vertices $A_k$ and $B_k$ are at distance
        $2d - 1$ in $G$; this is because $\Tc_G$ is a distance-$(d-1)$ pairing decomposition
        of $G$.
      
      Let us fix two vertices $u, v \in V(G)$ and assume that the shortest path in $U$
        between them passes through a vertex $x \not\in V(G)$.
      Then, $x \in V(H_k)$ for some $k \in [1, \ell_1]$.
      Hence, the path enters the subgraph $H_k$ through one of the vertices $A_k, B_k$,
        passes through $x$, and then leaves this subgraph through the other of these vertices.
      It means that this shortest path contains a subpath connecting $A_k$ and $B_k$ which
        is fully contained within $H_k$.
      By Claim \ref{tw-claim-bottom-long-path}, this subpath has length at least $2d + 2$,
        and therefore it can be replaced by the shortest path in $G$ connecting $A_k$ and $B_k$,
        which has length $2d - 1$.
      Therefore, the considered path was not a shortest path in the first place --- a contradiction.
      
      Hence, the shortest path between $u$ and $v$ lies fully within $G$ and thus
        $\dist_U(u, v) = \dist_G(u, v)$.
    \end{claimproof}
  \end{claim}
  
  Claim \ref{tw-claim-top-same-dists} significantly simplifies the structure of the distances
    between the vertices of $U$;
    no shortest path between a pair of vertices in $U$ can enter a subgraph $H_k$ for some
    $k \in [1, \ell_1]$ and then leave it, since there exists a strictly shorter path between
    these two vertices omitting this subgraph.
    
  Also, as a corollary of Claim \ref{tw-claim-top-same-dists}, we note the following:
    for each $k \in [1, \ell_1]$, we have that $\dist_U(A_k, B_k) = \dist_G(A_k, B_k) = 2d - 1$.
    
  We will now prove a claim that will directly imply that $\Tc_U$ is a distance-$d$
    pairing decomposition of $U$.
    
  \begin{claim}
    \label{tw-combined-distances}
    For each pair of indices $i, j \in [1, \ell_1\ell_2]$, we have that
    $$ \dist_U(b_i, a_j) = \begin{cases}
      2d & \text{if }i < j, \\
      2d+1 & \text{if }i \geq j.
    \end{cases} $$
    
    \begin{claimproof}
    \cqed
    Firstly, let us assume that $b_i$ and $a_j$
      are located in the same subgraph $H_k$ of $U$ for some $k \in [1, \ell_1]$;
      that is, $i, j \in [(k-1)\ell_2 + 1, k\ell_2]$.
    The shortest path in $U$ between these two vertices will have one of the following shapes:
    \begin{itemize}
      \item The path is fully contained within $H_k$. In this case,
        $$ \dist_U(b_i, a_j) = \dist_{H_k}(b_i, a_j) = \begin{cases}
        2d & \text{if }i < j, \\
        2d+1 & \text{if }i \geq j,
        \end{cases} $$
        as required.
      \item The path leaves $H_k$ and then reenters this subgraph.
        This would, however, require the path to contain a path between $A_k$ and $B_k$
          as a subpath.
        Since $\dist_U(A_k, B_k) = 2d - 1$, and neither $A_k$ nor $B_k$ can be an endpoint
          of the considered path, we get that
        $$ \dist_U(b_i, a_j) \geq 2 + \dist_U(A_k, B_k) \geq 2d + 1. $$
        However, there exists a path of equal or shorter length contained within $H_k$.
    \end{itemize}
    Therefore, there is a shortest path between $b_i$ and $a_j$ that is fully contained
      within $H_k$, which completes the proof in this case.
      
    \vspace{0.5em}
    Let us now assume that $a_i, b_i \in V(H_s)$ and $a_j, b_j \in V(H_t)$ for some
      $s, t \in [1, \ell_1]$,
      $s < t$, and let us compute the distances $\dist_U(b_i, a_j)$ and $\dist_U(b_j, a_i)$.
    By the construction of $U$, we infer that $i < j$.
    
    By our considerations above, we deduce that the shortest path from $b_i$ to $a_j$
      originates in $H_s$, then enters $G$ through either of the vertices $A_s$, $B_s$,
      then enters $H_t$ through either of the vertices $A_t$, $B_t$, and eventually
      terminates in $H_t$ (Figure \ref{tw-shortest-paths-fig}).
    The same applies to the shortest path from $b_j$ to $a_i$.
    
    We will split this shortest path into three parts: the first part, terminating at either of
      the vertices $A_s, B_s$; the second part, taking off from where the first path finished,
      and continuing to either of the vertices $A_t, B_t$; and the final part, continuing until
      the reaching $a_j$ or $b_j$.
    If there are multiple possible partitions, we pick the split minimizing the length of the first
      and the third part.
    This ensures us that the first path is fully contained within $H_s$, and the final part is
      fully contained within $H_t$.
    Since $A_s, B_s, A_t, B_t \in V(G)$, Claim \ref{tw-claim-top-same-dists} asserts that
      the second part of the path is also fully contained within $G$.
    
    \begin{figure}[h]
      \centering
      \input{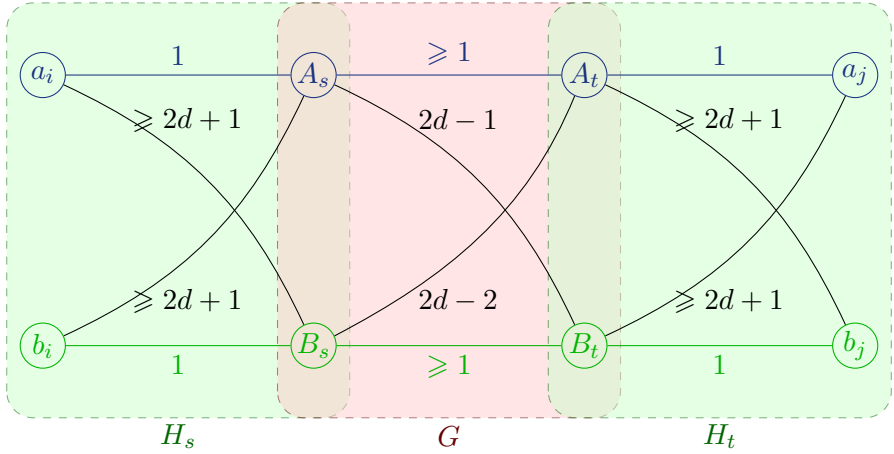}
      \caption{The important vertices on the shortest path from $b_i$ to $a_j$, and
        on the shortest path from $b_j$ to $a_j$.
        The edges are labeled by the length of the shortest path between the connected vertices
          or the lower bound on this length.}
      \label{tw-shortest-paths-fig}
    \end{figure}
    
    Since $\Tc_{H_s}$ is a neighboring tree decomposition of $H_s$, we figure that
    $$ \dist_{H_s}(a_i, A_s) = \dist_{H_s}(b_i, B_s) = 1, \quad \dist_{H_s}(a_i, B_s) \geq 2d+1, \quad
      \dist_{H_s}(b_i, A_s) \geq 2d + 1. $$
    Analogously, in $H_t$:
    $$ \dist_{H_t}(a_j, A_t) = \dist_{H_t}(b_j, B_t) = 1, \quad \dist_{H_t}(a_j, B_t) \geq 2d+1, \quad
      \dist_{H_t}(b_j, A_t) \geq 2d + 1. $$
    Finally, in $G$, since $A_s \neq A_t$, we get that $\dist_G(A_s, A_t) \geq 1$; similarly,
      $\dist_G(B_s, B_t) \geq 1$.
    By the definition of $G$ (which has a distance-$(d-1)$ pairing decomposition)
      and the fact that $s < t$, we infer that $\dist_G(B_s, A_t) = 2d - 2$
      and $\dist_G(B_t, A_s) = 2d - 1$.
      
    By verifying all possibilities, we can now conclude that (Figure \ref{tw-shortest-paths-fig}):
    \begin{equation*}
    \begin{split}
    \dist_U(b_i, a_j) &= \min_{X_s \in \{A_s, B_s\}} \min_{X_t \in \{A_t, B_t\}}
      \left[ \dist_{H_s}(b_i, X_s) + \dist_G(X_s, X_t) + \dist_{H_t}(X_t, a_j) \right] = \\
      & = \dist_{H_s}(b_i, B_s) + \dist_G(B_s, A_t) + \dist_{H_t}(A_t, a_j) = 2d
    \end{split}
    \end{equation*}
    and
    \begin{equation*}
    \begin{split}
    \dist_U(a_i, b_j) &= \min_{X_s \in \{A_s, B_s\}} \min_{X_t \in \{A_t, B_t\}}
      \left[ \dist_{H_s}(a_i, X_s) + \dist_G(X_s, X_t) + \dist_{H_t}(X_t, b_j) \right] = \\
      & = \dist_{H_s}(a_i, A_s) + \dist_G(A_s, B_t) + \dist_{H_t}(B_t, b_j) = 2d + 1.
    \end{split}
    \end{equation*}
    
    Hence, in this case, $i < j$ implies $\dist(b_i, a_j) = 2d$ and $\dist(b_j, a_i) = 2d+1$.
    Therefore, the statement of the claim holds in this case as well.
    
    \vspace{0.5em}
    Since the statement of the claim is true for each choice of a pair of indices $i, j$ ---
      no matter whether $b_i$ and $a_j$ are in the same subgraph $H_k$ for some $k$
      or not --- the proof is complete.
    \end{claimproof}
  \end{claim}
    
  Claim \ref{tw-combined-distances} directly proves the lemma.
  Therefore, there exists a procedure constructing a distance-$d$ pairing decomposition
    of $U$ of width $t$ and order $\ell_1\ell_2$.
  \end{proof}
\end{lemma}

We will name the procedure described in Lemma \ref{tw-combine-lemma} as
  $\mathsf{CombineNeighboring}$.
This procedure:
\begin{itemize}
  \item takes a graph $G$, together with its distance-$(d-1)$ pairing decomposition $\Tc_G$
    of order $\ell_1$ and width $t$;
  \item and a graph $H$, together with its neighboring distance-$d$
    pairing decomposition $\Tc_H$ of order $\ell_2$ and width $t$;
  \item and produces a graph $U$ and its distance-$d$ decomposition $\Tc_U$
    of order $\ell_1\ell_2$ and width $t$.
\end{itemize}
We apply this procedure in the following way:
$$ (U, \Tc_U) = \textsf{CombineNeighboring}((G, \Tc_G), (H, \Tc_H)). $$

We can now combine Lemmas \ref{make-neighboring-lemma} and \ref{tw-combine-lemma}
  in order to produce an even stronger procedure combining two pairing decompositions:
\begin{corollary}
  \label{combine-any-tw}
  There exists a procedure $\sf{Combine}$, which for $d \geq 2$, $t \geq 3$,
    $\ell_1, \ell_2 \geq 2$ takes:
  \begin{itemize}
  \item a graph $G$ with a distance-$(d-1)$ pairing decomposition $\Tc_G$
    of order $\ell_1$ and width $t$, and
  \item a graph $H$ with a distance-$d$ pairing decomposition $\Tc_H$
    of order $\ell_2$ and width $t - 2$,
  \end{itemize}
  and constructs a graph with a distance-$d$ pairing decomposition
  of order $\ell_1\ell_2$ and width $t$.
  
  \begin{proof}
  We set
  $$ \mathsf{Combine}((G, \Tc_G), (H, \Tc_H)) = \mathsf{CombineNeighboring}(
    (G, \Tc_G), \mathsf{MakeNeighboring}(H, \Tc_H)). $$
  By Lemma \ref{make-neighboring-lemma}, $\mathsf{MakeNeighboring}(H, \Tc_H)$
    is a graph with a distance-$d$ neighboring pairing decomposition of order
    $\ell_2$ and width $t$.
  By applying Lemma \ref{tw-combine-lemma} on $G$ and the newly produced graph,
    we create a graph with a distance-$d$ pairing decomposition of
    order $\ell_1\ell_2$ and width~$t$.
  \end{proof}
\end{corollary}

Using this corollary, we will now produce huge pairing decompositions.
Precisely, for each pair of integers $d \geq 1$, $k \geq 1$, we will define
  an undirected graph $G_{d,k}$, together with its distance-$d$ pairing decomposition
  $\Tc_{d,k}$ of width $2k+1$.

Since the tree decomposition of the graph in Example \ref{pairing-td-example}
  is distance-$1$ pairing and has width $3$, we shall take it as $\Tc_{1, 1}$, and the corresponding
  graph as $G_{1,1}$.
The decomposition $\Tc_{1, 1}$ has order $2$.

\begin{figure}[h]
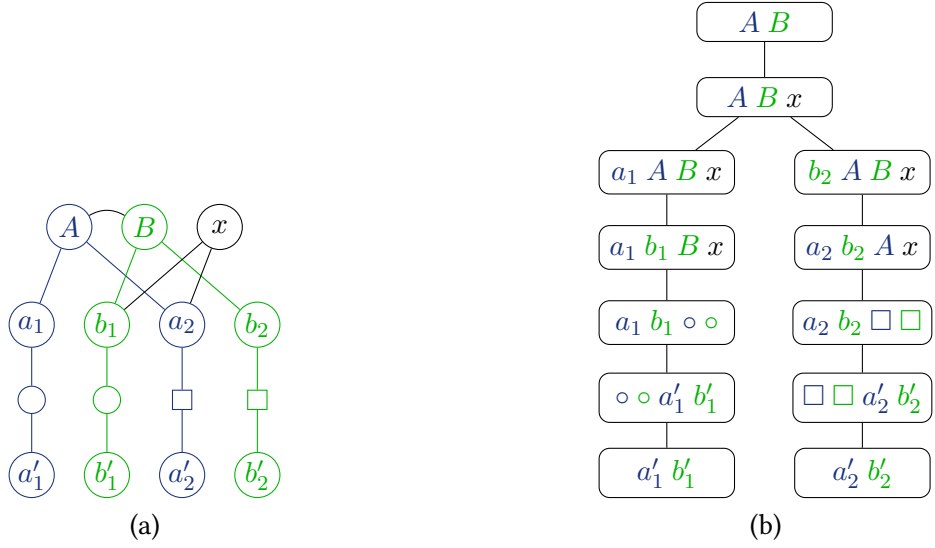

    \centering
    \begin{minipage}[b]{0.5\textwidth}
      \centering
      \input{figures/tw-longer-a.tex}
      (a)
    \end{minipage}\begin{minipage}[b]{0.5\textwidth}
      \centering
      \input{figures/tw-longer-b.tex}
      (b)
    \end{minipage}
    \caption{(a) --- the graph $G_{3, 1}$; (b) --- an example tree decomposition $T_{3, 1}$.}
    \label{tw-longer-fig}
  \end{figure}

For $d \geq 2$, we create the graph $G_{d, 1}$ and its distance-$d$ decomposition
  $\Tc_{d, 1}$ from $G_{1,1}$ and $\Tc_{1,1}$ in the
  following way: for each vertex $v \in \{a_1, b_1, a_2, b_2\}$, we append to the graph
  a path of length $d - 1$ whose one endpoint is $v$, and the other endpoint is a new vertex
  $v' \in \{a'_1, b'_1, a'_2, b'_2\}$.
We can easily create a distance-$d$ pairing
  decomposition $\Tc_{d, 1}$ of $G_{d, 1}$ of width $3$ (Figure \ref{tw-longer-fig}).
We remark that $\Tc_{d, 1}$ also has order $2$.

Also, for $k \geq 2$, we create the graph $G_{1, k}$ and its distance-$1$ pairing
  decomposition $\Tc_{1, k}$ from $G_{1,1}$ and $\Tc_{1,1}$ by adding $2(k-1)$
  isolated vertices to $G_{1, 1}$, and appending each new vertex
  to each bag of the decomposition $\Tc_{1, 1}$.
The resulting decomposition is obviously a distance-$1$ pairing decomposition of width
  $2k+1$ and order $2$.
  
Eventually, for $d \geq 2$, $k \geq 2$, we create the graph $G_{d,k}$ and the corresponding
  decomposition $\Tc_{d, k}$ by applying the procedure $\mathsf{Combine}$:
$$ (G_{d,k}, \Tc_{d, k}) := \mathsf{Combine}((G_{d-1,k}, \Tc_{d-1,k}), (G_{d,k-1}, \Tc_{d,k-1})). $$
We can easily verify that the input to this procedure complies with Corollary \ref{combine-any-tw},
  and leads to the construction of a graph with a distance-$d$ pairing decomposition
  of width $2k+1$,
  whose order is equal to the product of the orders of $\Tc_{d-1, k}$ and $\Tc_{d, k-1}$.
  
\vspace{0.5em}
For $d, k \geq 1$ we set $\ell_{d, k}$ as the order of $\Tc_{d, k}$.
\begin{lemma}
  \label{large-tw-order-lemma}
  For each $d, k \geq 1$, the following holds:
  $$ \ell_{d, k} = 2^{\displaystyle \binom{d + k - 2}{k - 1}}. $$
  \begin{proof}
  The formula can be directly verified for $d = 1$ and for $k = 1$.
  For $d, k \geq 2$, we know that the order of $\Tc_{d, k}$ is the product of the orders of
    $\Tc_{d-1, k}$ and $\Tc_{d, k-1}$.
  Therefore, by applying a simple induction on $d$ and $k$, we compute that
  \[
  \log_2 \ell_{d, k} = \log_2 \ell_{d-1,k} + \log_2 \ell_{d, k-1} =
    \binom{d + k - 3}{k - 1} + \binom{d + k - 3}{k - 2} = \binom{d + k - 2}{k - 1}.
  \]
  \end{proof}
\end{lemma}

Lemma \ref{large-tw-order-lemma} immediately leads to the following conclusion:

\begin{theorem}
  \label{tw-lower-bound}
  For every even $d \geq 2$ and odd $t \geq 3$, there exists a graph with
    treewidth not exceeding $t$ that contains a distance-$d$ half graph of order
  $$ 2^{\displaystyle \binom{\frac12(d + t - 5)}{\frac12(t-3)}}. $$
  \begin{proof}
  	We consider the graph $G_{\frac{d}{2}, \frac{t-1}{2}}$ together with its distance-$\frac{d}{2}$
  	  pairing decomposition $\Tc_{\frac{d}{2}, \frac{t-1}{2}}$.
  	The decomposition has order $\ell_{\frac{d}{2}, \frac{t-1}{2}}$
  	  (Lemma \ref{large-tw-order-lemma}), which matches the formula in the statement
  	  of the theorem; and width $t$, which ensures that the graph has treewidth not exceeding $t$.
  	Eventually, by the properties of pairing decompositions,
  	  $\Leaves(\Tc_{\frac{d}{2}, \frac{t-1}{2}})$
  	  is a distance-$d$ half graph of order $\ell_{\frac{d}{2}, \frac{t-1}{2}}$.
  \end{proof}
\end{theorem}

We also remark the asymptotic version of this theorem.

\begin{corollary}
  \label{tw-lower-bound-asym}
  Fix $t \geq 3$.
  Then for $d \geq 2$, the maximum order of a distance-$d$ half graph in the class
    of graphs with treewidth at most $t$ is bounded from below by
  $ 2^{\displaystyle d^{\Omega(t)}}. $
\end{corollary}

This discovery immediately leads to a lower bound in the class of $K_t$-minor-free graphs:

\begin{theorem}
  \label{kt-lower-bound}
  For every even $d \geq 2$ and odd $t \geq 5$, there exists a $K_t$-minor-free graph
    containing a distance-$d$ half graph of order 
  $$ 2^{\displaystyle \binom{\frac12(d + t - 7)}{\frac12(t-5)}}. $$
  \begin{proof}
    Theorem \ref{tw-lower-bound} asserts the existence of a graph with treewidth not exceeding
      $t-2$ containing a distance-$d$ half graph of order 
      $2^{\binom{\frac12(d + t - 7)}{\frac12(t-5)}}. $
    By Theorem \ref{tw-minor-free}, this graph is $K_t$-minor-free.
  \end{proof}
\end{theorem}

We immediately infer the following asymptotic result:

\begin{corollary}
  \label{kt-lower-bound-asym}
  Fix $t \geq 5$.
  Then for $d \geq 2$, the maximum order of a distance-$d$ half graph in the class
    of $K_t$-minor-free graphs is bounded from below by
  $ 2^{\displaystyle d^{\Omega(t)}}. $
\end{corollary}

\section{Polynomial bound on neighborhood complexity in planar graphs}\label{planar-neighborhood-chapter}
\label{noose-profile-lemma-section}
In this section, we prove that the neighborhood complexity of any vertex set $S$
  in planar graphs is bounded by a polynomial of $|S|$ and the distance $d$.
Namely, we prove the following theorem:

\begin{restatable}{theorem}{neighborhoodcomplthm}
  \label{neighborhood-complexity-planar}
  Consider any planar graph $G$, fix a non-empty set $A$ containing $c \geq 1$
  vertices of $G$.
  Then, the set $\{\pi_d[v, A]\,\mid\, v \in V(G)\}$ of all distance-$d$ profiles on $A$ has
    at most $64c^3(d+2)^7$ distinct elements.
\end{restatable}

We remark that it is already known that in each class $\Cc$ of graphs with bounded expansion
  the set $\{\pi_d[v, A]\,\mid\,v\in V(G)\}$ contains at most $f(d) \cdot |A|$ distinct
  elements for every non-empty set $A$, for some function $f$ depending only on
  $\Cc$~\cite{DBLP:journals/ejc/ReidlVS19}.
However, the proof of Reidl et al. requires $f$ to grow exponentially, even for planar graphs.
The polynomial dependence on~$d$ asserted by the Theorem \ref{neighborhood-complexity-planar}
  is crucial in the following parts of the proof.

The proof of Theorem \ref{neighborhood-complexity-planar} consists of two major parts.
Firstly, we prove the restricted version of the theorem, called the Noose Profile Lemma
  (Theorem \ref{noose-profile-lemma}).
Namely, we say that a non-empty subset $A$ of vertices of a graph $G$ is
  \emph{nice} if there exists an embedding of $G$ in the plane in which we can draw
  a~simple closed curve (called \emph{noose})
  passing through each vertex of $A$ exactly once, and not touching the interiors of
  edges at any point.
The Noose Profile Lemma states that the statement of Theorem~\ref{neighborhood-complexity-planar}
  holds for \emph{nice} sets of vertices, with a slightly better polynomial dependence
  on $|A|$ and $d$.
  
Then, in Subsection \ref{beyond-noose-profile-lemma-section}, we lift the lemma to
  Theorem \ref{neighborhood-complexity-planar}.
The proof of this part has been inspired by the concept of ,,cutting open'' the plane along
  the Steiner tree of a graph, present in the work of Pilipczuk
  et al.~\cite{DBLP:journals/talg/PilipczukPSL18}. 

The theorem is then used as a black-box in the upper bound proof in
  Section~\ref{planar-upper-bound-section}, so Sections~\ref{noose-profile-lemma-section}
  and \ref{planar-upper-bound-section} can be read independently.

\subsection{Statement of the Noose Profile Lemma}

In order to state Noose Profile Lemma, we first need to define nooses in plane graphs.

\begin{definition}
For any graph $G$ embedded in the plane,
  a {\bf noose} is a closed curve in the plane that passes through vertices and faces of this graph
  that does not intersect itself or the interiors of any edges.
Note that this definition disallows the curve to pass through any vertex multiple times.
\end{definition}

\begin{figure}[h]
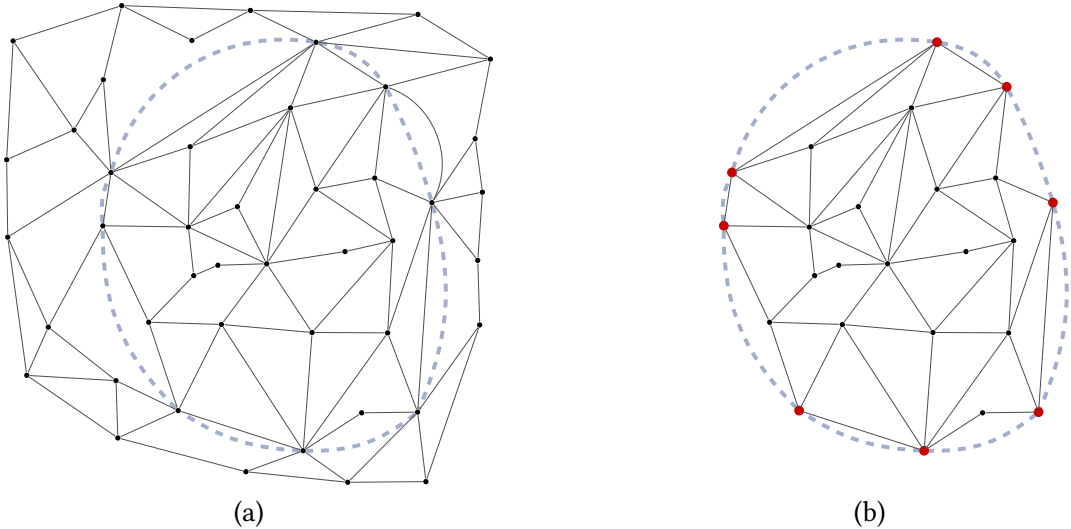

  \centering
  \begin{minipage}[b]{0.5\textwidth}
  \centering
  \input{figures/gl-subgraph-def-a.tex}
  (a)
  \end{minipage}\begin{minipage}[b]{0.5\textwidth}
  \centering
  \input{figures/gl-subgraph-def-b.tex}
  (b)
  \end{minipage}
  
  \caption{(a) --- an example graph $G$ with a noose $\Lc$ (blue); (b) --- the graph $G_\Lc$,
    with $V(\Lc)$ marked red.}
  \label{gl-subgraph-fig}
\end{figure}

We let the \textbf{length} of a noose to be the number of vertices on it.
For any noose $\Lc$ in a plane graph $G$, we define the set of vertices on the noose $V(\Lc)$.
Moreover, for each $v \in V(G)$, we define the distance-$d$ profile of $v$ on $\Lc$ in $G$:
  $\pi_d[v, \Lc] := \pi_d[v, V(\Lc)]$.
Also, we define the subgraph $G_\Lc$ as the graph containing all vertices and edges 
  within the bounded closed region whose boundary is $\Lc$ (Figure \ref{gl-subgraph-fig}).
The vertices and edges of $G_\Lc$ are said to be \textbf{enclosed} by~$\Lc$.

We are now ready to present the lemma in question:

\begin{theorem}[Noose Profile Lemma]
\label{noose-profile-lemma}
Consider any graph $G$ together with its embedding in the plane, and fix any noose
  $\Lc$ containing $c \geq 1$ vertices of $G$.
The set $\{\pi_d[v, \Lc]\, \mid\, v \in V(G_\Lc)\}$
  of different distance-$d$ profiles on $\Lc$,
  measured from the vertices enclosed by $\Lc$,
  has at most $c^3 (d+2)^4$ elements.
\end{theorem}

We remark that if $c$ is bounded by a polynomial of $d$, then the theorem above implies
  that the number of distance-$d$ profiles on $\Lc$ is also
  polynomially bounded on $d$.

Subsections \ref{closeness-def-section} through \ref{noose-profile-lemma-conclusion-subsection} are devoted to the proof of the lemma.



\subsection{Definition and properties of closeness}
\label{closeness-def-section}
Let us first describe some helpful notation.
For a closed region $\Omega \subseteq \mathbb{R}^2$ of the plane,
  we denote $\partial \Omega$ as the topological boundary of $\Omega$,
  and $\Int\Omega$ as its interior.
We also introduce the following majorization relation on distance-$d$ profiles in $G$:
  
\begin{definition}
  For a graph $G$, a subset of vertices $S \subseteq V(G)$, and two vertices
    $s, t \in V(G)$, we say that the distance-$d$ profile $\pi_d[s, S]$ {\bf majorizes} $\pi_d[t, S]$ if
    $\pi_d[s, S](v) \geq \pi_d[t, S](v)$ for all $v \in S$.
  We denote this partial order as $\pi_d[s, S] \succcurlyeq \pi_d[t, S]$.
\end{definition}

As in the statement of the lemma, fix any graph $G$ and a noose $\Lc$.
Let $v_1, v_2, \dots, v_c$ be the vertices of $\Lc$, listed
  in the clockwise order.
For simplicity, we set $v_0 = v_c$ and $v_{c+1} = v_1$.
  Let also $\Lc[i, j] \subseteq \Lc$
  be the curve from $v_i$ to $v_j$ following the noose clockwise.
That is, $\Lc[i, j]$ contains vertices $v_i, v_{i+1}, \dots, v_j$ in this order if $i \leq j$,
  and $v_i, v_{i+1}, \dots, v_c, v_1, v_2, \dots, v_j$ otherwise.
Moreover, let $\Lc(i, j) := \Lc[i, j] \setminus
  \{v_i, v_j\}$.
  
Let us say that a vertex $s \in V(G_\Lc)$ is \textbf{close} to a vertex $v_i \in \Lc$ if the distance
  between these two vertices in $G$ does not exceed $d$, and no shortest path connecting
  $v_i$ and $s$ contains any vertex from $\Lc$ as an internal vertex.
Otherwise, we say that $s$ is \textbf{far} from $v_i$.
For such vertex $s$ we also set $\mathrm{Close}(s) = \{v \in \Lc\,\mid\,
  s\text{ is close to }v\}$.
  
\begin{lemma}
  \label{lemma-majorization-close}
  If for any two vertices $s, t \in V(G_\Lc)$ we have $\pi_d[s, \mathrm{Close}(s)] \succcurlyeq
    \pi_d[t, \mathrm{Close}(s)]$, then $\pi_d[s, \Lc] \succcurlyeq \pi_d[t, \Lc]$.
  \begin{proof}
  Take any $v \in \Lc$. If also $v \in \mathrm{Close}(s)$, then
    $$ \pi_d[s, \Lc](v) = \pi_d[s, \mathrm{Close}(s)](v) \geq
      \pi_d[t, \mathrm{Close}(s)](v) = \pi_d[s, \Lc](v). $$
  If however $v \not\in \mathrm{Close}(s)$, then either of the following cases holds:
  \begin{itemize}
    \item The distance between $s$ and $v$ is longer than $d$.
      In this case, $\pi_d[s, \Lc](v) = +\infty \geq \pi_d[t, \Lc](v)$.
    \item The shortest path between $s$ and $v$ has length at most $d$, and
      some such shortest path passes through another vertex of $\Lc$.
      Pick a vertex of the noose $v' \in \Lc$ on any such shortest path.
      In case there are multiple such vertices, we pick any vertex minimizing its distance from $s$.
      As no shortest path between $s$ and $v'$ contains any other vertex of $\Lc$
        (otherwise, such vertex would be even closer to $s$),
        we have that $v' \in \mathrm{Close}(s)$.
      Hence,
      \begin{equation*}
        \begin{split}
          \pi_d[s, \Lc](v) &= \dist(s, v) = \dist(s, v') + \dist(v', v) = \pi_d[s, \mathrm{Close}(s)](v') +
          \dist(v', v) \geq \\ &\geq \pi_d[t, \mathrm{Close}(s)](v') + \dist(v', v) =
          \dist(t, v') + \dist(v', v) \geq \\
          &\geq \dist(t, v) = \pi_d[t, \Lc](v).
        \end{split}
      \end{equation*}
  \end{itemize}
  
  In all cases, we proved that $\pi_d[s, \Lc](v) \geq \pi_d[t, \Lc](v)$.
  Thus, $\pi_d[s, \Lc] \succcurlyeq \pi_d[t, \Lc]$.
  \end{proof}
\end{lemma}

This immediately leads to the following corollary:

\begin{corollary}
  \label{only-close-profile}
  If for any two vertices $s, t \in V(G_\Lc)$ we have $\mathrm{Close}(s) = \mathrm{Close}(t)$
    and $\pi_d[s, \mathrm{Close}(s)] = \pi_d[t, \mathrm{Close}(t)]$, then
    $\pi_d[s, \Lc] = \pi_d[t, \Lc]$.
\end{corollary}

Hence, $\pi_d[s, \Lc]$ can be uniquely deduced from the distances between $s$ and
  each vertex from $\mathrm{Close}(s)$.
This allows us to safely replace the graph $G$ with $G_\Lc$ (that is, to remove
  all vertices and edges not enclosed by $\Lc$); this replacement does not change
  $\mathrm{Close}(s)$ or any of the distances to the close vertices.

\vspace{1em}

In $G_\Lc$, let
  $\mathcal{C}_i := \{s \in V(G_\Lc)\,\mid\, |\mathrm{Close}(s)| = i\}$.
First, we examine the behavior of the vertices from $\mathcal{C}_0$ and $\mathcal{C}_1$.
Note that $\Lc \subseteq \mathcal{C}_1$.

\begin{lemma}
  There are at most $c(d + 1) + 1$ distinct distance-$d$ profiles on $\Lc$ measured
    from the vertices in $\mathcal{C}_0 \cup \mathcal{C}_1$.
  
\begin{proof}
Obviously, no vertex from $\mathcal{C}_0$ can reach $\Lc$ by a path consisting
  of at most $d$ edges,
  so each such vertex generates the same distance-$d$ profile on $\Lc$:
  the constant function $\lambda v.(+\infty)$.
  
Now consider any vertex $s \in \mathcal{C}_1$. Its distance-$d$ profile on $\Lc$
  can be unambiguously deduced from the only vertex $v_i \in \mathrm{Close}(s)$
  and the distance $\dist(s, v_i) \in [0, d]$ (Corollary \ref{only-close-profile}).
This leads to a maximum of $c(d + 1)$ distinct profiles created by $\mathcal{C}_1$.

Hence, at most $c(d+1) + 1$ distinct distance-$d$ profiles can be generated by the
  vertices in $\mathcal{C}_0 \cup \mathcal{C}_1$.
\end{proof}
\end{lemma}

From now on, we only consider remaining vertices.

\subsection{Rightmost shortest paths}

We will now extend slightly the definitions by Klein \cite{10.5555/1070432.1070454}.

\begin{figure}[h]
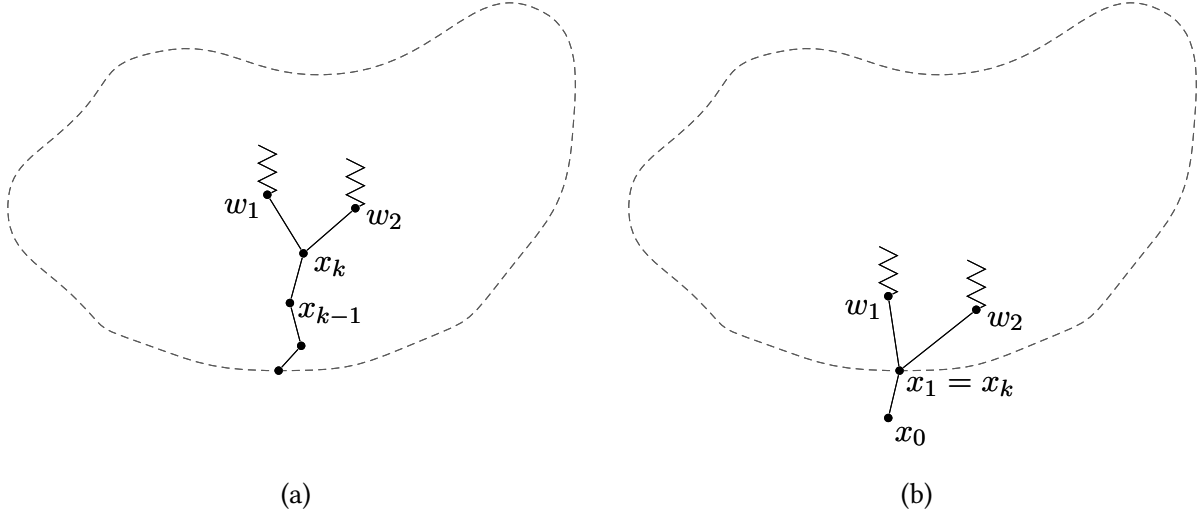

  \centering
  \begin{minipage}[b]{0.5\textwidth}
  \centering
  \input{figures/to-the-left-a.tex}
  (a)
  \end{minipage}\begin{minipage}[b]{0.5\textwidth}
  \centering
  \input{figures/to-the-left-b.tex}
  (b)
  \end{minipage}
  
  \caption{(a) --- the path containing $w_1$ is to the left of the path containing $w_2$. \\
  (b) --- if the first edges of two compared paths differ, we imagine an auxiliary vertex
    $x_0$ and an auxiliary edge $x_0x_1$ which we prepend to both paths.
    We then compare the paths as in (a).}
  \label{figure-to-the-left}
\end{figure}

\begin{definition}
\label{left-of-def}
  Fix any vertex $v \in \Lc$.
  We define a strict partial order $\ll$ on oriented paths originating from $v$, but not necessarily
    terminating at the same vertex.
  We say that $P_1 \ll P_2$ (or: \textbf{''$P_1$ is left of $P_2$''}) if:
  \begin{itemize}
    \item neither of the paths is a prefix of another,
    \item $P_1$ and $P_2$ have a common prefix $x_1, x_2, \dots, x_k$, where $v = x_1$,
    \item $P_i$ contains an oriented edge $x_k \to w_i$ for $i = 1, 2$,
    \item edge $x_k \to w_1$ is ``to the left'' of edge $x_k \to w_2$; formally,
        edges $x_kx_{k-1}$, $x_kw_1$ and $x_kw_2$ are in the clockwise order around $x_k$
        (Figure \ref{figure-to-the-left}(a)).
      However, if $k = 1$, the paths do not contain $x_0$.
      In order to cope with this problem, we imagine an auxiliary vertex $x_0$
        outside of the region enclosed by
        $\Lc$ and connect it by an edge to $x_1$, placed
        strictly outside of this region (Figure \ref{figure-to-the-left}(b)).
  \end{itemize}
  
  We can also interchangeably say $P_2 \gg P_1$ (or: \textbf{``$P_2$ is right of $P_1$''}.)
\end{definition}

The relation $\ll$ was defined in \cite{10.5555/1070432.1070454}
  only between the simple paths originating at $v$ and terminating at a fixed vertex $s$.
Then $\ll$ is a linear order.
In that setup, if $v$ and $s$ are in the same connected component of the graph,
  this relation limited only to the shortest paths between $v$ and $s$ contains the
  maximal element $P_{v, s}$, which we call the \textbf{rightmost shortest path}.
However, Definition \ref{left-of-def} naturally extends $\ll$ to any two paths $P_1$, $P_2$ which both originate
  at $v$ and neither is a prefix of the other, but which do not have to terminate at the
  same vertex.

We also define the rightmost shortest path tree rooted at $v$ as in
  \cite{10.5555/1070432.1070454}.
Notably, for any vertex $s$ of the tree, the path between $v$ and $s$ in this tree is $P_{v, s}$.
Obviously, if $v \in \Lc$ and for some vertices $a, b \in V(G_\Lc)$,
  $a$ belongs to $P_{v, b}$, then $P_{v, a} = P_{v, b}[v, a]$.
Moreover, for $v \in \Lc$ and $s \in V(G_\Lc)$, we obviously have that
  $\len(P_{v,s}) = \dist(v, s)$.

\begin{lemma}
  \label{prefix-paths}
  For any two distinct vertices $s, t \in \Lc$, paths
    $P_{a,s}$ and $P_{a,t}$ share a common prefix and are vertex-disjoint after removing
    this prefix.
  \begin{proof}
      Consider the rightmost shortest path tree rooted at $a$.
      The tree contains paths $P_{a,s}$ and $P_{a,t}$.
      The root $a$ belongs to both paths, hence they both
        must intersect at a common prefix and nowhere else.
  \end{proof}
\end{lemma}

\begin{lemma}
  \label{suffix-paths}
  Fix any $a, b \in \Lc$ and $s \in V(G_\Lc) \setminus V(\Lc)$
    such that $a, b \in \mathrm{Close}(s)$.
  Paths $P_{a,s}$ and $P_{b,s}$ share a common suffix as a subpath, and
    are vertex-disjoint after removing this suffix.
    
  \begin{proof}
    Take the first vertex $x \in P_{a,s} \cap P_{b,s}$ on the path between $a$ and $s$,
      that is, the one minimizing $\dist(x, a)$.
    If $x = s$, we are done.
    In the opposite case, we set $a^*$ and $b^*$ as the next vertices after $x$ on the paths
      $P_{a,s}$ and $P_{b,s}$, respectively.

    Moreover, we set $a_p$ and $b_p$ as the last vertices before $x$ on these paths.
    Vertices $a_p$ and $b_p$ exist since $P_{a,s}$ cannot pass through $b$ (otherwise,
      $s$ would not be close to $a$), and similarly $P_{b,s}$ cannot pass through $a$.
    Hence, $x \not \in \{a, b\}$ --- that is, $x$ cannot be the initial vertex of either
      of these paths.

    Obviously, the length of each of the paths $P_{a,s}[x,s]$ and $P_{b,s}[x,s]$
      must be the same, otherwise we would be able to swap one of the suffixes with another
      in one of the paths $P_{a,s}, P_{b,s}$, reducing its length.
    This is, however, impossible as we require $P_{a,s}$ and $P_{b,s}$ to be the
      shortest paths connecting the corresponding pairs of vertices.
    Also, note that paths $P_{a,x}$ and $P_{b,x}$ are vertex-disjoint apart from their common
      endpoint $x$.
      
    Consider the closed region $\Omega$ bounded by $P_{a,x}$, $P_{b,x}$ and the part of
      the noose between $a$ and $b$ --- either $\Lc[a, b]$ or $\Lc[b, a]$ ---
      chosen so that $s \in \Int\Omega$.
    (Obviously, $s$ cannot belong to $\Lc$, $P_{a,x}$ or $P_{b,x}$.)
    Assume now that $a^* \not\in \Omega$ (Figure \ref{figure-suffix-paths}(a)).
    Then, the path $P_{a, s}[a^*, s]$ must intersect $\partial \Omega$.
    Hence, one of the following cases occurs:
    \begin{itemize}
      \item The path intersects $P_{a,x}$ (Figure \ref{figure-suffix-paths}(b)).
        In this case, the path $P_{a, s}$ would intersect itself.
      \item The path intersects $P_{b,x}$ at some point $y \neq x$
        (Figure \ref{figure-suffix-paths}(c)).
        In this case, $P_{a,s}$ contains vertices
          $x, y$ in this order, while $P_{b, s}$ contains these vertices in the opposite order.
        Hence,
          $$ \len(P_{a,s}[y, s]) < \len(P_{a,s}[x, s]) = \len(P_{b,s}[x, s]) < \len(P_{b,s}[y, s]). $$
        Therefore, $P_{b,s}$ could be shortened by replacing its suffix $P_{b,s}[y, s]$ by
          $P_{a, s}[y, s]$.
      \item The path enters $\Omega$ through an edge whose either endpoint is in
        $\Lc \setminus \{a, b\}$.
        This is impossible as this edge would not be enclosed
        by the noose and therefore would not belong to $G_\Lc$.
    \end{itemize}
    \begin{figure}[h]
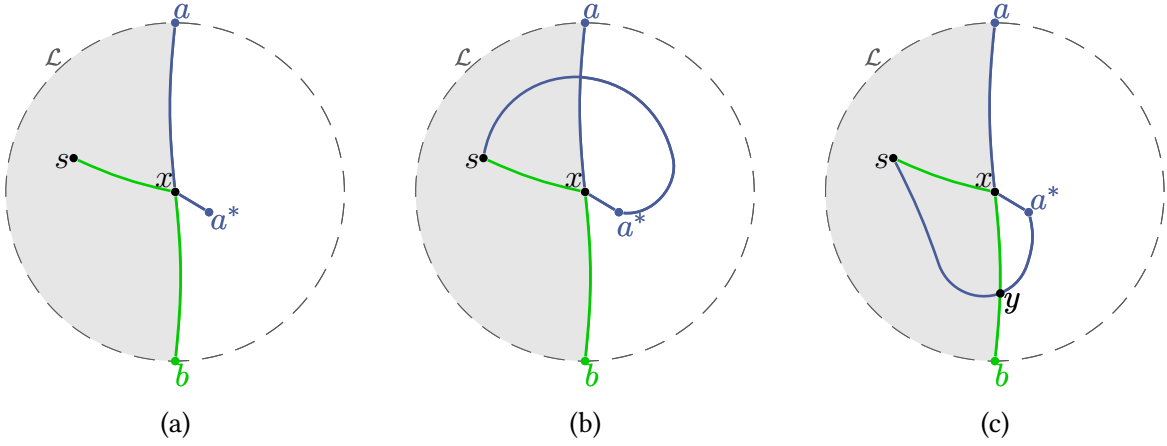

      \centering
      \begin{minipage}[b]{0.33\textwidth}
      \centering
      \input{figures/suffix-paths-fig-a.tex}
      (a)
      \end{minipage}\begin{minipage}[b]{0.33\textwidth}
      \centering
      \input{figures/suffix-paths-fig-b.tex}
      (b)
      \end{minipage}\begin{minipage}[b]{0.33\textwidth}
      \centering
      \input{figures/suffix-paths-fig-c.tex}
      (c)
      \end{minipage}
      \caption{(a) --- an example setup. The noose $\Lc$ is a gray dashed circle.
        The region $\Omega$ is shaded.
        Green path is $P_{b,s}$, and the blue path is the prefix of $P_{a,s}$. \\
        (b) --- a case where $P_{a,s}[a^*, s]$ intersects $P_{a,x}$. \\
        (c) --- a case where $P_{a,s}[a^*, s]$ intersects $P_{b,x}$.}
      \label{figure-suffix-paths}
    \end{figure}
    
    We arrived at a contradiction in each of the cases above.    
    Therefore $a^* \in \Omega$, and analogously $b^* \in \Omega$.

    Assume now that $a^* \neq b^*$.
    Then, we remind that $P_{a,s}$ is rightmost among all shortest paths from $a$ to $s$.
    Consider the path $Q_a = P_{a,x} \cdot P_{b,s}[x, s]$.
    As $\mathrm{len}(P_{a, s}[x, s]) = \mathrm{len}(P_{b, s}[x, s])$, we also have that
      $\mathrm{len}(Q_a) = \mathrm{len}(P_{a, s})$.
    Since $Q_a \neq P_{a, s}$, we get that $Q_a \ll P_{a,s}$.
    Notice that $a_p$ and $x$ are two last common vertices of $Q_a$ and $P_{a,s}$, while
      $b^*$ is the first vertex of $Q_a$ not belonging to $P_{a,s}$, and
      $a^*$ is the first vertex of $P_{a,s}$ not belonging to $Q_a$.
    By the definition of $\ll$, edges $xa_p$, $xb^*$ and $xa^*$
      form a clockwise order around $x$.
    We similarly deduce that edges $xb_p$, $xa^*$ and $xb^*$ form
      a clockwise order around $x$.
    This implies the clockwise order of four edges around $x$: $xa_p, xb^*, xb_p, xa^*$.
    As $a_px$ and $xb_p$ are two consecutive edges on the boundary of $\Omega$,
      we conclude that this boundary separates $a^*$ and $b^*$.
    That is however impossible as $a^*, b^* \in \Omega$.
    
    From the contradiction above, $a^* = b^*$.    
    This means that $P_{a,s}$ and $P_{b,s}$ share an oriented edge $xa^*$.
    We now prove that they must share the remaining suffix of the path from $x$ to $s$.
    Otherwise, there would be an oriented edge $fg$ shared by both $P_{a,s}$ and $P_{b,s}$,
      followed by an oriented edge $gh_a$ in $P_{a,s}$, and by an oriented edge $gh_b$
      in $P_{a,s}$ ($h_a \neq h_b$).
    Assume without loss of generality that the oriented edges $gf, gh_a, gh_b$ are oriented
      clockwise around $g$.
    By the Definition \ref{left-of-def}, this means that
    $$ P_{a, s} = P_{a,s}[a, g] \cdot P_{a,s}[g, s] \ll P_{a,s}[a,g] \cdot P_{b,s}[g,s] $$
    --- a contradiction.
    Hence, these paths must share the remaining suffix of the path
      from $x$ to~$s$.

    As $P_{a,s}[a,x] \cap P_{b,s}[b,x] = \{x\}$ and $P_{a,s}[x,s] = P_{b,s}[x,s]$,
     this finishes the proof.
  \end{proof}
\end{lemma}

Fix any vertex $s \in \mathcal{C}_k$, $k \geq 2$ and assume that it is close to vertices
  $v_{i_1}, v_{i_2}, \dots, v_{i_k}$ of the noose, listed in the clockwise order
  ($1 \leq i_1 < i_2 < \dots < i_k \leq c$).
For simplicity, we also denote $i_{k+1} = i_1$.
Now, consider $k$ paths $P_{v_{i_1}, s}, P_{v_{i_2}, s}, \dots, P_{v_{i_k}, s}$.

\begin{lemma}
  \label{disjoint-paths}
  Exactly one of the following conditions holds:
  \begin{enumerate}[(1)]
  \item There exists $j \in [1, k]$ such that paths
    $P_{v_{i_j}, s}$ and $P_{v_{i_{j+1}}, s}$ are vertex disjoint
    apart from~$s$,
  \item There exists a vertex $s' \neq s$ such that all the paths
    $P_{v_{i_1}, s}, P_{v_{i_2}, s}, \dots, P_{v_{i_k}, s}$ contain $s'$.
    This vertex is adjacent to $s$ and is the penultimate vertex in each of these paths.
  \end{enumerate}
  \begin{proof}
    Obviously, condition (2) implies that (1) is false.
    We only need to prove the opposite implication.
    
    Assume that property (1) holds for no $j \in [1, k]$.
    Therefore, by Lemma \ref{suffix-paths}, we see that for each $j \in [1,k]$, paths
      $P_{v_{i_j},s}$ and $P_{v_{i_{j+1}}, s}$ share the common suffix, containing at least
      two vertices.
    Hence, each path $P_{v_{i_1},s}, P_{v_{i_2},s}, \dots, P_{v_{i_k},s}$ shares the same
      common suffix, having at least two vertices.
    The penultimate vertex in this suffix is $s'$ as in property (2).
  \end{proof}
\end{lemma}

Let $V_1$ be the set of vertices from $\mathcal{C}_2 \cup \mathcal{C}_3 \cup \dots \cup
  \mathcal{C}_c$ for which the first condition in Lemma~\ref{disjoint-paths} holds,
  and $V_2$ --- the analogous set for the second condition.
  
\begin{lemma}
  \label{profile-shift}
  For every $s \in V_2$, there exists a vertex $s' \in V_1$ and an integer
    $m \in \{1, 2, \dots, d\}$ such that
  $$\forall_{v \in \Lc}\ \pi_d[s, \Lc](v) = \mathrm{Cap}_d\left(\pi_d[s',
     \Lc](v) +
    m\right), $$
  where
  $$\mathrm{Cap}_d(x) = \begin{cases}
    x & \text{ if } x \leq d, \\
    +\infty & \text{ otherwise.}
  \end{cases}$$
  
  \begin{proof}
  If $s \in V_2$, then there exists a sequence $(s_0, s_1, \dots, s_{m-1}, s_m)$ where
    $s = s_m$ and $m \geq 1$, which is the longest common suffix of all the paths
    $P_{v_{i_1}, s}, P_{v_{i_2}, s}, \dots, P_{v_{i_k}, s}$.
  The rightmost shortest paths from $v_{i_1}, v_{i_2}, \dots, v_{i_k}$ to $s_0$ can be
    obtained by dropping $m$ last vertices in each of the paths.
  In this setup, we have $s_0 \in V_1$ --- otherwise, we would be able to prepend the sequence
    defined above by the common neighbor of $s_0$ on all the paths terminating at $s_0$.
  We will now prove that $s' = s_0$ and $m$ satisfy the condition stated in the lemma.

  Pick any vertex $v \in \Lc$ and consider the following cases:
\begin{itemize}
\item If $s$ is close to $v$, then $\dist(s, v) = \dist(s_0, v) + m$ (because of the common
  prefix above).
\item If $s$ is far from $v$ and $\pi_d[s, \Lc](v) = \infty$, then
  $\dist(s, v) > d$.
  Hence, $\dist(s_0, v) > d-m$ as $\dist(s, v) \leq \dist(s, s_0) + \dist(s_0, v) = m + \dist(s_0, v)$.
  To put it in the other way, $\dist(s_0, v) + m > d$.
  Therefore $\pi_d[s_0, \Lc](v) + m > d$, which is consistent with the condition
    $\pi_d[s, \Lc](v) = \infty$ required by the statement of the lemma.
\item If $s$ is far from $v$ and $\pi_d[s, \Lc](v) \leq d$, then by definition of closeness
  some shortest path from $s$ to $v$ passes through another vertex $v' \in V(\Lc)
  \setminus \{v\}$.
  Assume that $v'$ is one of the vertices closest to $s$ with this property.

  We see that $s$ is close to $v'$.
  Thus, $v' \in \{v_{i_1}, \dots, v_{i_k}\}$, so $P_{v',s}$ is one of the paths in the set
    $\{P_{v_{i_1}, s}, \dots, P_{v_{i_k}, s}\}$, and as such, it contains $s_0$.
  Therefore, some shortest path from $v$ to $s$ contains $v, v', s_0, s$ in this order.
  Hence, $\pi_d[s, \Lc](v) = \dist(s, v) = \dist(s, s_0) + \dist(s_0, v) = m + \dist(s_0, v) =
  	\pi_d[s_0, \Lc](v) + m$.
\end{itemize}
Therefore, the distance-$d$ profile $\pi_d[s, \Lc]$ is constructed by considering
  the profile $\pi_d[s_0, \Lc]$ where $s_0 \in V_1$,
  and increasing each entry in the profile by $m$,
  while remembering to replace each integer exceeding $d$ by~$\infty$.
\end{proof}
\end{lemma}  

\begin{corollary}
  \label{profile-shift-cor}
  If
    \[ |\{\pi_d[v, \Lc]\, \mid\, v \in V_1\}| = B \quad\text{for an integer }B,\]
  then
    \[ |\{\pi_d[v, \Lc]\, \mid\, v \in V_1 \cup V_2\}| \leq B \cdot (d + 1). \]
  \begin{proof}
    By Lemma \ref{profile-shift}, the profile of a vertex $v \in V_2$ can be deduced
      from a distance-$d$ profile of a vertex $v' \in V_1$ and an integer $m \in \{1, 2, \dots, d\}$.
    Hence,
    \[ |\{\pi_d[v, \Lc]\, \mid\, v \in V_2\}| \leq
      |\{\pi_d[v, \Lc]\, \mid\, v \in V_1\}| \cdot |\{1, 2, \dots, d\}| = B \cdot d. \]
    We also include the vertices in $V_1$ into this inequality and get:
    \[ |\{\pi_d[v, \Lc]\, \mid\, v \in V_1 \cup V_2\}| \leq
      |\{\pi_d[v, \Lc]\, \mid\, v \in V_1\}| + |\{\pi_d[v, \Lc]\, \mid\, v \in V_2\}| \leq B \cdot (d+1). \]
  \end{proof}
\end{corollary}

Corollary \ref{profile-shift-cor} allows us to focus on 
  bounding the number of different distance-$d$ profiles
  measured only from the vertices in $V_1$.
  
\subsection{Buckets}

We will now group all the vertices in $V_1$ into a reasonably small number of buckets.
The buckets will be chosen in such a way that all the vertices within one group expose
  similar properties in $G$.
This will eventually allow us to bound the number of different distance-$d$ profiles measured
  by each individual bucket.

\begin{figure}
  \centering
  \input{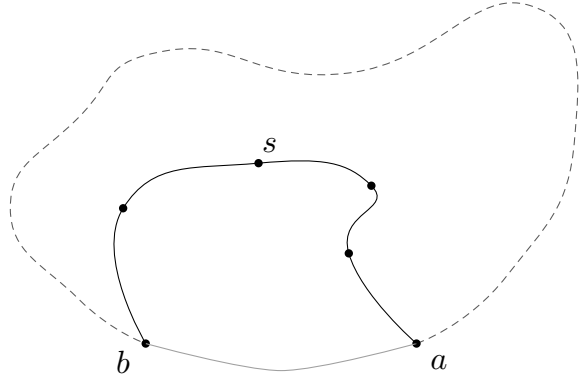}
  \caption{Definition of a bucket assigned to a vertex $s$. Here, $d_a = 3$ and $d_b = 2$.
    Vertex $s$ is far from any vertex on the dashed part of the noose.}
  \label{bucket-def-figure}
\end{figure}

\begin{definition}
  For a vertex $s \in V_1$, we define the \textbf{bucket} of $s$ as any quadruple
  $(a, b, d_a, d_b)$ satisfying the following conditions (Figure \ref{bucket-def-figure}):
  \begin{itemize}
    \item $a, b \in V(\Lc)$ and $a, b \in \mathrm{Close}(s)$,
    \item $\Lc(b, a) \cap \mathrm{Close}(s) = \varnothing$; in other words,
      each vertex $v \in \Lc(b, a)$ is far from $s$,
    \item paths $P_{a, s}$ and $P_{b, s}$ intersect only at $s$,
    \item $d_a = \dist(a, s)$ is the distance between $a$ and $s$ in $G_{\Lc}$,
      and $1 \leq d_a \leq d$,
    \item $d_b = \dist(b, s)$ is the distance between $b$ and $s$ in $G_{\Lc}$,
      and $1 \leq d_b \leq d$.
  \end{itemize}
  
  \noindent Lemma \ref{disjoint-paths} proves that such an assignment exists since $s \in V_1$.
\end{definition}

Notice that there at most $c(c - 1) d^2$ possible quadruples
  $(a, b, d_a, d_b)$ --- $a$ and $b$ are two distinct vertices of the noose in some order,
  and $d_a, d_b$ are two positive distances not exceeding $d$.
  
Fix a bucket $(a, b, d_a, d_b)$ and consider the set $\mathcal{S} =
  \mathcal{S}(a, b, d_a, d_b)$ of
  vertices $s$ that are assigned to the bucket.
Remember that each $s \in \mathcal{S}$ is far from each vertex of
  $\Lc(b, a)$.
  
\vspace{1em}
  
 \begin{lemma}
    \label{hard-cross-paths}
    For any two different vertices $s, t \in \mathcal{S}$, if paths
      $P_{a,t}$ and $P_{b,s}$ intersect at any vertex $x$, then $P_{a,s} \gg P_{a,t}$,
      $P_{b,t} \gg P_{b,s}$, and $\dist(x, t) = \dist(x, s) > 0$ (Figure \ref{hard-cross-paths-config}).
    
    \begin{figure}[h]
      \centering
      \input{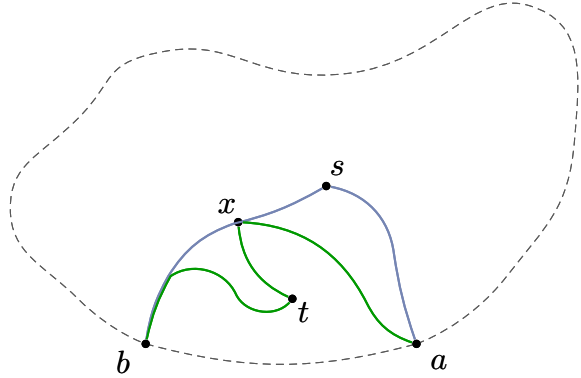}
      \caption{A possible configuration in Lemma \ref{hard-cross-paths}. Paths
        $P_{a,s}$ and $P_{b,s}$ are marked blue, while $P_{a,t}$ and $P_{b,t}$ are green.}
      \label{hard-cross-paths-config}
    \end{figure}
      
    \begin{proof}
    Assume that $P_{a, t}$ and $P_{b, s}$ intersect at $x$.
    Then we have a couple of simple equalities:
    $$
    \dist(a, t) = \dist(a, x) + \dist(x, t) \qquad\text{and}\qquad
    \dist(b, s) = \dist(b, x) + \dist(x, s).
    $$
    From the triangle inequality and the fact that $s, t \in \mathcal{S}(a, b, d_a, d_b)$, we infer:
    \begin{equation}
    \label{xs-geq-xt-eq}
    \dist(a, x) + \dist(x, s) \geq \dist(a, s) = d_a = \dist(a, t) = \dist(a, x) + \dist(x, t),
    \end{equation}
    \begin{equation}
    \label{xt-geq-xs-eq}
    \dist(b, x) + \dist(x, s) = \dist(b, s) = d_b = \dist(b, t) \leq \dist(b, x) + \dist(x, t).
    \end{equation}
    We infer that $\dist(x, s) \geq \dist(x, t)$ from (\ref{xs-geq-xt-eq}) and
      $\dist(x, s) \leq \dist(x, t)$ from (\ref{xt-geq-xs-eq}).
    Therefore, $\dist(x, s) = \dist(x, t)$ and all the inequalities above are satisfied
      with equality.
    As $x \neq s$ (otherwise we would have $x = s = t$), obviously $\dist(x, s) > 0$.

    Hence, $\dist(a, s) = \dist(a, x) + \dist(x, s)$ and $\dist(b,t) = \dist(b,x) +
      \dist(x,t)$.
      
    We now see two shortest paths from $a$ to $s$: the rightmost path $P_{a,s}$, and
      another path $Q := P_{a,x} \cdot P_{b, s}[x, s]$.
      
    Notice that $P_{a,s}$ cannot be the prefix of $P_{a,x}$ as it is strictly longer than $P_{a,x}$.
    Also $P_{a,x}$ cannot be the prefix of $P_{a,s}$
      because $x$ lies on $P_{b,s}$ (which is vertex disjoint with $P_{a,s}$ apart from $s$,
      and $x \neq s$).
    Hence, we have either $P_{a,s} \ll P_{a,x}$ or $P_{a,s} \gg P_{a,x}$.
    We however cannot have $P_{a,s} \ll P_{a,x}$,
      or otherwise we would have $P_{a,s} \ll P_{a,x} \cdot P_{b,s}[x,s] = Q$,
      which is impossible by the definition of $P_{a,s}$.
    Therefore $P_{a,s} \gg P_{a,x}$.
    As $P_{a,x} = P_{a,t}[a, x]$, we also get $P_{a,s} \gg P_{a,t}$.

    An analogous argument shows that $P_{b, t} \gg P_{b, s}$.
    \end{proof}
  \end{lemma}
  
  \begin{corollary}
    For any two distinct vertices $s, t \in \mathcal{S}$, if paths
      $P_{a,s}$ and $P_{b,t}$ intersect at any vertex $x$, then $P_{a,t} \gg P_{a,s}$,
      $P_{b,s} \gg P_{b,t}$, and $\dist(x, t) = \dist(x, s) > 0$.
    \label{hard-cross-paths-sym}
    \begin{proof}
      In the statement of Lemma \ref{hard-cross-paths}, swap $s$ and $t$.
    \end{proof}
  \end{corollary}
  
  \begin{corollary}
    \label{hard-cross-paths-both}
    It is impossible for $P_{a,t}$ and $P_{b,s}$ to intersect, and for $P_{a,s}$ and
    $P_{b,t}$ to intersect at the same time.
    \begin{proof}
      The former condition would imply $P_{a,s} \gg P_{a,t}$ (Lemma \ref{hard-cross-paths}),
      while the latter --- $P_{a,t} \gg P_{a,s}$ (Corollary \ref{hard-cross-paths-sym}).
    \end{proof}
  \end{corollary}

\subsection{Gamma-regions}

We fix any bucket $\mathcal{S} = \mathcal{S}(a, b, d_a, d_b)$.
We will now assign a region of the plane to each vertex $w \in S$.

\begin{lemma}
For any vertex $w \in \mathcal{S}$, the curve $\Lc[a, b] \cup P_{a,w} \cup P_{b,w}$
  is a closed Jordan curve.
\begin{proof}
Notice that $w \not\in V(\Lc)$ and
  $w$ is close to $a$, so $P_{a,w}$ cannot contain vertices from $\Lc$
  as internal vertices. Therefore,
  $P_{a,w} \cap \Lc[a, b] = \{a\}$. Analogously,
  $P_{b,w} \cap \Lc[a, b] = \{b\}$. Finally,
  $P_{a,w} \cap P_{b,w} = \{w\}$ due to the disjointness of these paths in the sense of
    Claim \ref{disjoint-paths}.
Thus the union of these curves is a closed Jordan curve.
\end{proof}
\end{lemma}

\begin{figure}[h]
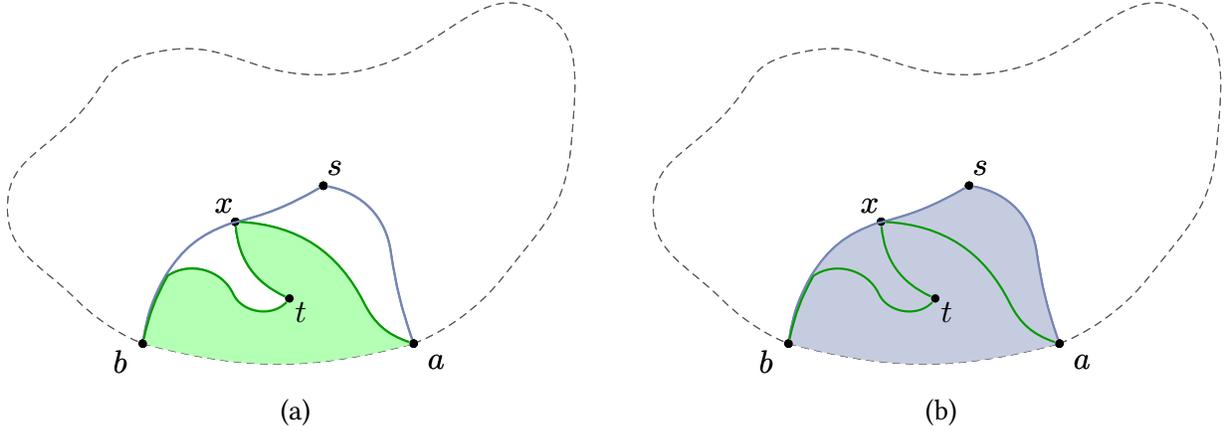

\centering
    \begin{minipage}[t]{0.48\textwidth}
	\centering
	\input{figures/cross-lemma-area-a.tex}
    (a)
    \end{minipage}\begin{minipage}[t]{0.04\textwidth}\ 
    \end{minipage}\begin{minipage}[t]{0.48\textwidth}
	\centering  
	\input{figures/cross-lemma-area-b.tex}
    (b)
    \end{minipage}
  \caption{Gamma regions in the configuration from Figure \ref{hard-cross-paths-config}:
    (a) --- $\Gamma_t$, (b) --- $\Gamma_s$.}
  \label{gamma-def-fig}
\end{figure}

\begin{definition}
For any vertex $w \in \mathcal{S}$, we define the \textbf{gamma-region}
  $\Gamma_w$ as the
  closed, bounded region of the plane whose boundary is $\Lc[a, b] \cup P_{a,w} \cup P_{b,w}$
  (Figure \ref{gamma-def-fig}).
\end{definition}

We will now prove the potential of this definition.
Intuitively, for any two vertices $s, t \in \mathcal{S}$ in the bucket,
  we will either have $t \in \Gamma_s$, which will imply a serious structural relation
    between the gamma-regions $\Gamma_s$ and $\Gamma_t$,
  or we will have $t \not\in \Gamma_s$, which will in turn have major implications
    for the distance-$d$ profiles $\pi_d[s, \Lc]$ and $\pi_d[t, \Lc]$.
    
  \begin{lemma}
    \label{gamma-contain}
    For any two different vertices $s, t \in \mathcal{S}$, if $t \in \Gamma_s$, then
      $\Gamma_t \subsetneq \Gamma_s$.
    \begin{proof}
      Assume that $t \in \Gamma_s$.
      We will first prove a couple of helper claims.
      
      \begin{claim}
      \label{gamma-contain-claim-b}
      \cqed
      If $P_{b,t}$ is disjoint with $P_{a,s}$,
        then $P_{b,t} \gg P_{b,s}$ and $P_{b,t}$ is fully contained within $\Gamma_s$.
        \begin{claimproof}
        Assume that $P_{b,t}$ is disjoint with $P_{a,s}$.
        As $\dist(b, t) = \dist(b, s)$ and $s \neq t$,
          neither of paths $P_{b,s}$ and $P_{b,t}$ is a prefix of another and thus
          they are ordered by $\ll$.
        We notice that $P_{b,s}$ follows $\partial \Gamma_s$ clockwise.
        This is because we defined $\Lc[a, b]$ to run clockwise from $a$ to $b$ on
          $\Lc$, and $P_{b,s}$ is a directed path originating from $b$, where
          $\Lc[a, b]$ terminated.
        Hence, both $P_{b,s}$ and $\Lc[a, b]$ have to be oriented the same way
          on $\partial \Gamma_s$.
        
        Let $(y_1, y_2, \dots, y_{k-1}, y_k)$ ($y_k = t$) be the maximal suffix of $P_{b,t}$
          which is disjoint with $P_{b,s}$.
        For any choice of $i \in \{1, 2, \dots, k\}$, we have that
          $y_i \not\in P_{a,s}$ since $P_{b,t}$ is disjoint with $P_{a,s}$,
          and $y_i \not\in \Lc[a,b]$ as $\Lc[a,b] \cap P_{b,t} = \{b\}$.
        Hence $\{y_1, y_2, \dots, y_k\} \cap \partial \Gamma_s = \varnothing$, and thus
          these vertices are either all strictly inside of $\Gamma_s$, or all strictly outside of it.

        Now, as $P_{b,s}$ is oriented clockwise, it can be easily seen that $y_1 \in \Gamma_s\ 
          \Leftrightarrow\ P_{b,y_1} \gg P_{b,s}$;
          in other words, the interior of $\Gamma_s$ is ``to the right'' of $P_{b,s}$.
        
        If $P_{b,y_1} \gg P_{b,s}$, then $P_{b,t} \gg P_{b,s}$ and $P_{b,s}$ is fully within
          $\Gamma_s$.
        
        If $P_{b,y_1} \ll P_{b,s}$, then $P_{b,t} \ll P_{b,s}$. Because $y_1 \not\in
          \Gamma_s$, we have that $\{y_1, y_2 \dots, y_k\}$ is disjoint with $\Gamma_s$.
        In particular, $t \not\in \Gamma_s$.
        This however contradicts our assumption that $t \in \Gamma_s$.
        \end{claimproof}
      \end{claim}
      
      \begin{claim}
      \label{gamma-contain-claim-a}
      \cqed
      If $P_{a,t}$ is disjoint with $P_{b,s}$,
        then $P_{a,t} \ll P_{a,s}$ and $P_{a,t}$ is fully contained within $\Gamma_s$.
        \begin{claimproof}
        The proof of this claim follows mostly the proof of Claim \ref{gamma-contain-claim-b}.
        However, $P_{a,s}$ follows $\partial \Gamma_s$ anti-clockwise.
        This is because $\Lc[a, b]$ runs clockwise along the noose from $a$ to $b$,
          and both $\Lc[a, b]$ and $P_{a, s}$ originate from the same vertex $a$.
        Hence, these two paths are oriented in the opposite way along $\partial \Gamma_s$.
        This means that the interior of $\Gamma_s$ is ``to the left'' of $P_{a,s}$.
        Hence, $t \in \Gamma_s$ if and only if $P_{a,t} \ll P_{a,s}$.
        \end{claimproof}
      \end{claim}
      
      Note that if both $P_{b,t}$ is disjoint with $P_{a,s}$, and $P_{a,t}$ is disjoint with $P_{b,s}$,
        then by Claims \ref{gamma-contain-claim-b} and \ref{gamma-contain-claim-a}
        we get that $\partial\Gamma_t \subseteq \Gamma_s$ and thus $\Gamma_t \subseteq
        \Gamma_s$.
      
      Assume that $t \in \Gamma_s$, but $\Gamma_t \not\subseteq \Gamma_s$.
      As $\partial\Gamma_t$ is a Jordan curve, there exists a vertex $y$ on
        $\partial\Gamma_t$ not belonging to $\Gamma_s$.
      We will consider multiple cases, depending on which part of $\partial\Gamma_t$
        this vertex belongs to.
      
      \vspace{1em}
      \underline{Case 1}: $y \in \Lc[a, b]$. This is however impossible as $\Lc[a, b]
        \subseteq \Gamma_s$.
        
      \vspace{1em}
      \underline{Case 2a}: $y \in P_{b,t}$ and $P_{b,t} \gg P_{b,s}$.
        In this case, $P_{b,t}$ cannot intersect $P_{a,s}$ by Corollary~\ref{hard-cross-paths-sym}.
        Hence, by Claim~\ref{gamma-contain-claim-b}, we have that $P_{b,t}$ is fully contained
        within $\Gamma_s$ --- a contradiction because $y \in P_{b,t}$ must lie outside
        of $\Gamma_s$.
          
      \vspace{1em}
      \underline{Case 2b}: $y \in P_{b,t}$ and $P_{b,t} \ll P_{b,s}$.
        In this case, Claim \ref{gamma-contain-claim-b} requires
          that the intersection between $P_{b,t}$ and $P_{a,s}$ must exist.
        Hence, we also get $P_{a,t} \gg P_{a,s}$ by Corollary \ref{hard-cross-paths-sym}.

        Now, $P_{a,t}$ cannot intersect $P_{b,s}$ as this is forbidden by Corollary
          \ref{hard-cross-paths-both}.
        This leads to a contradiction due to Claim~\ref{gamma-contain-claim-a}.
          
      \vspace{1em}
      \underline{Case 3a}: $y \in P_{a,t}$ and $P_{a,t} \gg P_{a,s}$.
        Similarly to the first part of Case 2b, we deduce that $P_{a,t}$ must intersect $P_{b,s}$.
        However, by Lemma \ref{hard-cross-paths}, we get that $P_{a,t} \ll P_{a,s}$ ---
          a contradiction.
      
      \vspace{1em}
      \underline{Case 3b}: $y \in P_{a,t}$ and $P_{a,t} \ll P_{a,s}$.
        For our convenience, we take $y$ to be the first vertex of $P_{a,t}$ not in $\Gamma_s$.
        We again deduce that $P_{a,t}$ intersects $P_{b,s}$; otherwise, Claim
          \ref{gamma-contain-claim-a} would imply that $P_{a,t}$ is
          contained within $\Gamma_s$, which contradicts the choice of $y$.
          
        Let $x_1$ and $x_2$ be the first and the last intersection of $P_{a,t}$ with $P_{b,s}$,
          respectively.
        Notice that $y$ cannot precede $x_1$ on $P_{a,t}$; as $P_{a,t} \ll P_{a,s}$, then
          $P_{a,y} \ll P_{a,s}$.
        As $P_{a,y}$ would not intersect $P_{b,s}$, we could prove (analogously to the proofs
          of Claims \ref{gamma-contain-claim-b} and \ref{gamma-contain-claim-a})
          that $P_{a,y}$ is fully contained within $\Gamma_s$,
          which contradicts our assumption.
        Vertex $y$ cannot also appear after $x_2$ on $P_{a,t}$ --- in this case,
          $y \not \in \Gamma_s$ would imply $t \not\in \Gamma_s$
          (note that $P_{a,t}[y, t]$ would not be able to intersect $\partial\Gamma_s$),
          which again is impossible.

        Since vertices $a, x_1, y, x_2, t$ appear on $P_{a,t}$ in this order,
          we have that $x_1 \in P_{a, x_2}$.
        Also, by Lemma \ref{hard-cross-paths} and the fact that $\dist(x_2, t) < \dist(x_1, t)$,
          we get that
        \begin{equation*}
        \begin{split}
        \dist(b, x_1) &= \dist(b,s) - \dist(x_1, s) = \dist(b,s) - \dist(x_1, t) < \dist(b,s) - \dist(x_2, t) =\\ 
          &= \dist(b,s) - \dist(x_2, s) = \dist(b, x_2).
        \end{split}
        \end{equation*}
        We infer that $x_1 \in P_{b,s}[b, x_2] = P_{b, x_2}$.
        
        Now, consider two paths: $P_{a,x_2} = P_{a,t}[a, x_2]$ and $P_{b,x_2} = P_{b,s}[b, x_2]$.
        By Lemma \ref{suffix-paths}, these two paths must intersect on their common suffix
          and nowhere else.
        As $x_1 \in P_{a,x_2} \cap P_{b,x_2}$, we get that $P_{a,x_2}[x_1, x_2] = P_{b,x_2}[x_1, x_2]$.
        This means that $P_{b,s}$ shares its vertices with $P_{a,t}$ on the segment between
          $x_1$ and $x_2$.
        Hence,
        $$ P_{a,t}[x_1, x_2] = P_{a, x_2}[x_1, x_2] = P_{b, x_2}[x_1, x_2] =
          P_{b, s}[x_1, x_2] \subseteq P_{b,s} \subseteq \partial \Gamma_s. $$  
        However, we get a contradiction because of our assumptions that
          $y \in P_{a, x_2}[x_1, x_2]$ and $y \not\in \Gamma_s$.
      \vspace{1em}
      
      After considering all cases above, we have that $\Gamma_t \subseteq \Gamma_s$.
      Moreover, $t \not\in P_{a,s}$ (since $\dist(a,s) = \dist(a,t)$) and analogously $t \not\in P_{b,s}$, and
        $t \not\in V(\Lc)$ (otherwise $t$ would only be close to itself).
      Thus, $t \not\in \partial \Gamma_s$.
      Since obviously $t \in \partial\Gamma_t$, we immediately
        get that $\Gamma_s \neq \Gamma_t$.
      Therefore, $\Gamma_t \subsetneq \Gamma_s$.
    \end{proof}
  \end{lemma}
  
  \vspace{1em}
  
  The claim above described the case $t \in \Gamma_s$.
  Now we try to see what happens when $t \not\in \Gamma_s$.

  \begin{lemma}
    \label{majorization-lemma}
    If $t \not\in \Gamma_s$ for some $s, t \in \mathcal{S}$, then
      $\pi_d[t, \Lc] \succcurlyeq \pi_d[s, \Lc]$.
  \begin{proof}
    We will prove that $\pi_d[t, \mathrm{Close}(t)] \succcurlyeq \pi_d[s, \mathrm{Close}(t)]$.
    Lemma \ref{lemma-majorization-close} will then conclude our proof.
    
    Recall, by the choice of $a$ and $b$, that $\mathrm{Close}(t) \subseteq \Lc[a, b]$.
    We will now prove that $\dist(v, s) \geq \dist(v, t)$ for each $v \in \Lc[a, b]$ such that
      $v \in \mathrm{Close}(t)$.
    For $v \in \{a, b\}$ the inequality is satisfied with equality; thus, we only need to care
      about the case $v \in \Lc(a, b) \cap \mathrm{Close}(t)$.

    Notice that $v \in \Gamma_s$, but on the other hand $t \not\in \Gamma_s$.
    Therefore, we can take $x$ as the last intersection of $P_{v,t}$ with
      $\partial \Gamma_s$.
    We cannot have $x = v$, otherwise the edge connecting $x$
      with the next vertex on the path would lie outside of the noose
      since $v \in \Lc(a, b)$.
    Also, since $v \in \mathrm{Close}(t)$, we get that $P_{v, t} \cap V(\Lc) = \{v\}$;
      hence, $P_{v, t}$ cannot intersect $\Lc$ at any vertex other than $v$.
    Therefore, $x \in P_{a,s} \cup P_{b,s}$.
    Without loss of generality, assume that $x \in P_{a, s}$.
    
    As $x$ lies both on the shortest path from $a$ to $s$, and the shortest path from
      $v$ to $t$, we have the following equalities:
    $$\dist(v, t) = \dist(v, x) + \dist(x, t) \qquad\text{and}\qquad
      \dist(a, s) = \dist(a, x) + \dist(x, s).$$
    Moreover, $\dist(a, t) = \dist(a, s) = d_a$.
      
    From the triangle inequality and the equations above, we get
    $$ \dist(a, x) + \dist(x, t) \geq \dist(a, t) = \dist(a, s) = \dist(a, x) + \dist(x, s), $$
    or equivalently, $\dist(x, t) \geq \dist(x, s)$.
    
    Then, from the conditions above and another use of triangle inequality,
      we have
    \[ \dist(v, t) = \dist(v, x) + \dist(x, t) \geq \dist(v, x) + \dist(x, s) \geq \dist(v, s). \]
  \end{proof}
  \end{lemma}

  \subsection{Conclusion of the proof of Noose Profile Lemma}
  \label{noose-profile-lemma-conclusion-subsection}
  
  Assume there are $M$ vertices
    $s_1, s_2, \dots, s_M \in \mathcal{S}$ with pairwise distinct distance-$d$ profiles
    on $\Lc$.
  Pick any vertices $s_i, s_j$, $i \neq j$, and let us describe their relation with respect to
    the regions $\Gamma_{s_i}, \Gamma_{s_j}$ defined by them:
  \begin{itemize}
  \item If $s_i \not\in \Gamma_{s_j}$ and $s_j \not\in \Gamma_{s_i}$, then from
    Lemma \ref{majorization-lemma} we get $\pi_d[s_i, \Lc] = \pi_d[s_j, \Lc]$ ---
    a contradiction.
  \item If $s_i \in \Gamma_{s_j}$ and $s_j \in \Gamma_{s_i}$, the from Lemma
    \ref{gamma-contain}
    we get $\Gamma_{s_i} \subsetneq \Gamma_{s_j} \subsetneq \Gamma_{s_i}$ ---
    another contradiction.
  \item Hence, either $\pi_d[s_i, \Lc] \prec \pi_d[s_j, \Lc]$
    (if $s_j \not\in \Gamma_{s_i}$) or $\pi_d[s_j, \Lc] \prec \pi_d[s_i, \Lc]$
    (if $s_i \not\in \Gamma_{s_j}$).
  \end{itemize}
  
  Therefore, the relation $\preccurlyeq$ defined on the distance-$d$ profiles for $s_1, s_2,
    \dots, s_M$
  is antisymmetric, transitive (by definition of $\preccurlyeq$) and connex
    (by our considerations above).
  Hence, it is a linear order and thus contains a chain of length $M$.
  However, as $\pi_d$ is a function whose domain contains $c$ elements,
    and its codomain contains $d+2$ possible values,
  the maximum length of such a chain is $c(d+1) + 1$.
  We thus get $M \leq c(d+1) + 1$.
    
  Summing everything up, we proved that each set $\mathcal{S}(a, b, d_a, d_b)$
    generates at most $c(d+1) + 1$ different distance-$d$ profiles.
  There are at most $c(c-1)d^2$ such sets, as $a \neq b$ and $1 \leq d_a, d_b \leq d$.
  Hence, the number of different profiles induced in $G$ by $V_1$ is bounded by

  $$c(c - 1)d^2 \cdot [c(d+1) + 1] \leq c^3 (d+1)^3.$$
  
  The number of profiles in $\mathcal{C}_2 \cup \mathcal{C}_3 \cup \dots \cup
    \mathcal{C}_c = V_1 \cup V_2$ is then bounded by
    $c^3 (d+1)^4$ (Corollary \ref{profile-shift-cor}).
  The vertices from $\mathcal{C}_0 \cup \mathcal{C}_1$ induce at most
    $c(d + 1) + 1$ additional profiles.
  The total number of profiles can be thus bounded by
  
  $$ c^3 (d+1)^4 + c(d + 1) + 1 \leq c^3(d+2)^4. $$
    
  This concludes the proof of the Noose Profile Lemma.
  We note a simple corollary which is a direct consequence of it:
  
  \begin{corollary}
  \label{noose-profille-lemma-rev}
    Consider any graph $G$ together with its embedding in the plane, and fix any noose
  $\Lc$ containing $c \geq 1$ vertices of $G$.
The set $\{\pi_d[v, \Lc]\, \mid\, v \not\in V(G_\Lc)\}$
  of different distance-$d$ profiles on $\Lc$,
  measured from the vertices \textbf{not enclosed} by $\Lc$,
  has at most $c^3 (d+2)^4$ elements.
    
    \begin{proof}
	  Let us take the embedding of $G$ in the plane, and apply to it an inversion
	    with any positive radius
	    and a center in any point enclosed by $\Lc$ that does not lie on any vertex or
	    edge of $G$.
	  After this transformation,
         the images of all vertices of $G$ that have not been enclosed by $\Lc$
         in the original embedding will now be enclosed by the inversive image $\Lc'$ of the noose.
       Noose Profile Lemma applies to the inversed embedding as well,
         which allows us to reason that there are at most $c^3 (d+2)^4$ different distance-$d$
         profiles on $\Lc$, measured from the vertices whose images are enclosed by $\Lc'$
         --- or, in other words, measured from the vertices which were not enclosed by $\Lc$
         in the original embedding.
    \end{proof}
  \end{corollary}

\subsection{Proof of Theorem \ref{neighborhood-complexity-planar}}
\label{beyond-noose-profile-lemma-section}
We now lift the Noose Profile Lemma to the general setting, where the noose can
  be replaced with any non-empty set of vertices.
For convenience, we restate the theorem here.

\neighborhoodcomplthm*
  \begin{proof}
    We fix $G$ and embed it in the plane.
    We also fix the set $A$ as in the statement of the theorem.
    Let $V_d = \{v \in V(G)\,\mid\,\exists_{s\in A} \,\dist(v, s) \leq d\}$ be the set of vertices
      of $G$ which are at~distance at most $d$ from any vertex in $A$.
    We observe that for each vertex $v \not\in V_d$, the distance-$d$ profile
      $\pi_d[v, A]$ is the constant function equal to $+\infty$.
    Now, we define the graph $G_d = G[V_d]$ induced by $V_d$.
    We can easily see that for each vertex $v \in V_d$, the distance-$d$ profile of $v$ on $A$
      in $G_d$ is identical to the distance-$d$ profile of $v$ on $A$ in $G$.
    Hence, the number of different distance-$d$ profiles on $A$ in $G_d$ is at most $1$ less
      than the number of different distance-$d$ profiles on $A$ in $G$.    
      
    We first assume that $G_d$ is connected;
      we will resolve the case where $G_d$ is disconnected at the end of the proof.
    Let $\Tc$ be the minimum Steiner tree on $A$ in $G_d$; that is, the smallest possible tree
      which is a connected subgraph of $G_d$ and spans all vertices of $A$.
    \begin{claim}
    \label{claim-small-steiner-tree}
    $\Tc$ has at most $(c-1)(2d+1)$ edges.
    \begin{claimproof}
    \cqed
    Let $G_A$ be the complete weighted graph whose set of vertices is $A$, and for every
      pair of vertices $u, v \in A$, the edge between these vertices has weight $\dist_{G_d}(u, v)$.
    Let also $\Tc_A$ be the minimum-weight spanning tree of $G_A$.

    We now prove that no edge of $\Tc_A$ has weight greater than $2d + 1$.
    Assume for contradiction that $uv \in E(\Tc_A)$ is an edge with weight at least $2d + 2$.
    In this case, take a shortest path from $u$ to $v$ in $G_d$:
      $(x_0, x_1, x_2, \dots, x_\delta)$, where $x_0 = u$, $x_\delta = v$, $\delta \geq 2d+2$.
    We have that $\dist_{G_d}(x_{d+1}, u) = d + 1$ and $\dist_{G_d}(x_{d+1}, v) = \delta - (d+1)$.
    As $x_{d+1} \in V_d$, there exists a vertex $s \in A$ at distance at most $d$ from $x_{d+1}$.
    By the triangle inequality, we infer that
    \begin{equation*}
    \begin{split}
    \dist_{G_d}(s, u)& \leq \dist_{G_d}(s, x_{d+1}) + \dist_{G_d}(x_{d+1}, u) \leq
      d + (d + 1) < \delta, \\
    \dist_{G_d}(s, v)& \leq \dist_{G_d}(s, x_{d+1}) + \dist_{G_d}(x_{d+1}, v) \leq
      d + (\delta - (d+1)) < \delta.
    \end{split}
    \end{equation*}
    Hence, $G_A$ contains a three-vertex cycle $(u, v, s)$, in which the edge $uv$ is strictly
      the heaviest.
    Therefore, $uv$ cannot belong to $\Tc_A$ --- a contradiction.
    
    \smallskip

    Now, for each edge $uv \in \Tc_A$, we connect $u$ with $v$ in $G_d$ using the shortest
      path between these two vertices.
    We can easily see that the union of these paths is a connected subgraph of $G_d$ containing
      all the vertices in $A$.
    Hence, $\Tc$, the minimum Steiner tree on $A$, is at most as long as the union of these paths.
    Since we used $|A|-1 = c - 1$ paths, and each of them has length at most $2d + 1$,
      we infer that $\Tc$ has at most $(c-1)(2d+1)$ edges.
    \end{claimproof}
    \end{claim}    
    \newcommand{\Gcut}{G_{\mathrm{cut}}}
    We now create a new graph $\Gcut$, which is a modified version of $G_d$,
      by ``cutting the plane open'' along the Steiner tree $\Tc$ (Figure \ref{cutting-fig}).
    This process is inspired by a similar idea by Pilipczuk et al.
      \cite{DBLP:journals/talg/PilipczukPSL18}.
    Formally, $\Gcut$ is the result of the following process:
  
  \begin{figure}[h]
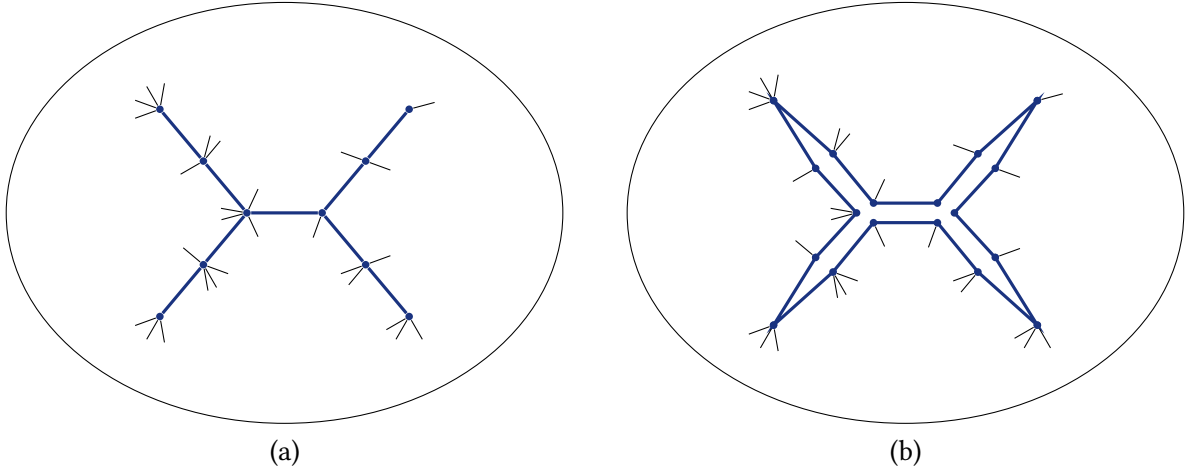

    \centering
    \begin{minipage}[b]{0.5\textwidth}
      \centering
      \input{figures/cut-tree-a.tex}
      (a)
    \end{minipage}\begin{minipage}[b]{0.5\textwidth}
      \centering
      \input{figures/cut-tree-b.tex}
      (b)
    \end{minipage}
    \caption{Creating $\Gcut$. (a) --- original graph $G_d$, with tree $\Tc$ (blue).
      (b) --- the resulting graph, with an Euler tour of $\Tc$.}
    \label{cutting-fig}
  \end{figure}
    
    \begin{enumerate}
      \item We find an Euler tour of $\Tc$ which does not intersect itself in the plane embedding
        of $G_d$. This tour traverses each edge of $\Tc$ twice, in the opposite directions.
      \item For each vertex $v \in V(\Tc)$, we let $s_1, s_2, \dots, s_\Delta$ to be the neighbors of
        $v$ in $\Tc$, ordered anti-clockwise around $v$.
        For simplicity, we assume $s_{\Delta+1} = s_1$.
        If $\Delta \geq 2$, we create $\Delta$~copies of $v$ (including the original vertex $v$):
          $v_1, v_2, \dots, v_\Delta$, and we put the $i$-th copy, $v_i$, very close to the
          original vertex $v$, and
          between the edges $vs_i$ and $vs_{i+1}$ in the embedding of $G_d$.
        Moreover, for each neighbor $x$ of $v$ such that $vx \not \in E(\Tc)$, we find the only
          index $i$ such that the edges $vs_i$, $vx$, $vs_{i+1}$ are ordered anti-clockwise around
          $v$;
          then, we replace the edge $vx$ with $v_ix$.
      \item We duplicate each edge $uv$ of the Euler tour of $\Tc$, so that we have one copy for
        each orientation of $uv$ in the Euler tour.
        We then use both copies to connect the corresponding copies of vertices $u, v$
          in $\Gcut$.
        We can easily see that all new edges form a cycle enclosing~$\Tc$, which we call $C$.
        Naturally, $C$ has twice as many edges as $\Tc$.
      \item We remove the edges of $\Tc$ from the graph.
    \end{enumerate}
    Let $V(C)$ be the set of vertices of $C$, and $E(C)$ --- the set of edges of $C$.
    Since $|E(\Tc)| \leq (c-1)(2d+1)$ (Claim \ref{claim-small-steiner-tree}),
      we have that $|E(C)| = 2|E(\Tc)| \leq 2(c-1)(2d+1)$.
    Also, for every vertex $v \in V(\Tc)$, we let $\deg_\Tc(v)$ to be the degree of $v$ in $\Tc$.
    Hence, for $v \in V(\Tc)$, the graph $\Gcut$ contains vertices $v_1, \dots, v_{\deg_\Tc(v)}$.

    \begin{claim}
    \label{claim-steiner-dist-profile}
    For any two different vertices $p, q \in V(\Gcut) \setminus V(C)$,
      if $\pi_d[p, V(C)] = \pi_d[q, V(C)]$ in $\Gcut$, then
      $\pi_d[p, V(\Tc)] = \pi_d[q, V(\Tc)]$ in $G_d$.
      \begin{claimproof}
      \cqed
      We will show that the distance-$d$ profile of any vertex $p \in V(\Gcut) \setminus V(C)$
        on $V(C)$ in $\Gcut$ uniquely determines the distance-$d$ profile of $p$ on $V(\Tc)$
        in $G_d$.
      The statement of the lemma will follow.
      
      We first remark that for each $x \in V(\Tc)$ and $i \in [1, \deg_\Tc(v)]$, we have that
        $\dist_{\Gcut}(p, x_i) \geq \dist_{G_d}(p, x)$.
      This is true since we can pick a shortest path from $p$ to $x_i$ in $\Gcut$, and replace
        each vertex of this path belonging to $V(C)$ with the corresponding vertex of $V(\Tc)$
        in $G_d$.
      The resulting path connects $p$ with $x$ in $G_d$, and has length $\dist_{\Gcut}(p, x_i)$.
      
      We fix a vertex $v \in V(\Tc)$, and we want to determine $\pi_d[p, V(\Tc)](v)$ in $G_d$ from
        $\pi_d[p, V(C)]$ in $\Gcut$.
      Obviously, for each $x \in V(\Tc)$ and $i \in [1, \deg_\Tc(v)]$ we have
      $$ \dist_{G_d}(p, v) \leq \dist_{G_d}(p, x) + \dist_{G_d}(x, v) \leq 
        \dist_{\Gcut}(p, x_i) + \dist_{G_d}(x, v). $$
      On the other hand, we pick a shortest path $P$ from $p$ to $v$ in $G_d$.
      We let $x$ to be the first intersection of $P$ with $V(\Tc)$, and $y$ --- the vertex
        immediately preceding $x$ on $P$.
      Since $y \not\in V(\Tc)$, then by the construction of $\Gcut$ we infer that $y$ is adjacent to
        a copy of $x$ in $\Gcut$ --- say, $x_i$ for some $i \in [1, \deg_\Tc(x)]$.
      For this choice of $x$ and $i$ we have that
      \begin{equation*}
      \begin{split}
      \dist_{G_d}(p, v) &= \dist_{G_d}(p, y) + 1 + \dist_{G_d}(x, v) =
        \dist_{\Gcut}(p, y) + \dist_{\Gcut}(y, x_i) + \dist_{G_d}(x, v) \geq \\
        &\geq \dist_{\Gcut}(p, x_i) + \dist_{G_d}(x, v).
      \end{split}
      \end{equation*}
      We conclude that
      $$ \dist_{G_d}(p, v) = \min \{ \dist_{\Gcut}(p, x_i) + \dist_{G_d}(x, v)\,\mid\,
        x \in V(\Tc), i \in [1, \deg_\Tc(x)]\}. $$
      It is now straightforward to determine $\pi_d[p, V(\Tc)](v)$ in $G_d$ from
        $\pi_d[p, V(C)]$ in $\Gcut$ and the distances in $G_d$ between each pair of vertices
        of $\Tc$ (which are independent on the choice of~$p$). 
      \end{claimproof}
    \end{claim}
    
    We now observe that $C$ is a cycle in $\Gcut$, which does not enclose any
      other vertices of $\Gcut$; hence, there exists a noose
      $\Lc$, closely following $C$, that passes through each vertex of $C$ exactly once
      and that does not enclose any other vertex of $\Gcut$.
    The variant of Noose Profile Lemma (Corollary \ref{noose-profille-lemma-rev}) applies
      to $\Lc$.
    Hence, there are at most $|C|^3 (d+2)^4$ different distance-$d$ profiles on $C$ in
      $\Gcut$, measured from the vertices not enclosed by $\Lc$.

    We fix two vertices $p, q \in V(\Gcut) \setminus V(C)$ such that
      $\pi_d[p, V(C)] = \pi_d[q, V(C)]$ in $\Gcut$.
    Firstly, we see that neither $p$ nor $q$ are enclosed by $\Lc$.
    We know from Claim \ref{claim-steiner-dist-profile} that $\pi_d[p, V(\Tc)] = \pi_d[q, V(\Tc)]$
      in $G_d$.
    Since $A \subseteq V(\Tc)$, we infer that $\pi_d[p, A] = \pi_d[q, A]$ in $G_d$.
    By the construction of $G_d$, we also get that $\pi_d[p, A] = \pi_d[q, A]$ in $G$.

    Since there are at most $|C|^3 (d+2)^4$ different distance-$d$ profiles on $C$
      in $\Gcut$, measured from the vertices not enclosed by $\Lc$,
      we get that there are at most $|C|^3 (d+2)^4$ different distance-$d$ profiles on $A$
      in $G_d$, measured from the vertices in $V(G_d) \setminus V(\Tc)$.
    Therefore, $G_d$ has at most $|C|^3 (d+2)^4 + |V(\Tc)|$ different distance-$d$ profiles
      on $A$ in total.
    Recall that $G$ has at most one more distance-$d$ profile on $A$ compared to $G_d$
      --- that is, the constant function equal to $+\infty$.
    Hence, the number of distance-$d$ profiles on $A$ in $G$ is bounded by
    \begin{equation*}
    \begin{split}
    |C|^3 (d+2)^4 + |V(\Tc)| + 1 &\leq
      [2(c-1)(2d+1)]^3(d+2)^4 + (c-1)(2d+1) + 1 \leq \\ &\leq
      8c^3[(2d+1)^3(d+2)^4 + (2d+1) + 1] \leq 64c^3(d+2)^7.
    \end{split}
    \end{equation*}
    Let us finally consider the case where $G_d$ is disconnected.
    Assume that the vertices of $A$ are spread among $k$ different connected components
      in $G_d$, and the $i$-th connected component contains $c_i$ vertices of $A$ for each
      $i \in [1, k]$ ($c_1, c_2, \dots, c_k \geq 1$, $\sum c_i = c$).
    Then, we apply this theorem to each connected component of $G_d$ separately.
    We get that the number of different distance-$d$ profiles on $A$ in $G_d$ is bounded by
    $$ \sum_{i=1}^k 64c_i^3(d+2)^7 < 64c^3(d+2)^7. $$
    Hence, the number of different distance-$d$ profiles on $A$ in $G$ is also bounded by
      $64c^3(d+2)^7$.
  \end{proof}

\section{Upper bound on the orders of semi-ladders in planar graphs}
\label{planar-upper-bound-section}

This section is fully dedicated to the proof of an upper bound on the maximum
  order of a distance-$d$ semi-ladder in planar graphs:

\begin{theorem}
\label{planar-upper-bound}
  The order of any distance-$d$ semi-ladder in a planar graph
    is bounded from above by $d^{O(d)}$.
\end{theorem}
  
We refer to the introduction (Section \ref{introduction-chapter}) for a rough sketch
  of the proof.
Here, we provide more details on the structure of the proof.

In Subsection \ref{quasi-cages-subsection}, we introduce two structural objects:
  quasi-cages and cages.
These objects occur in planar graphs with sufficiently large semi-ladders.
Namely, we prove in Lemma \ref{quasi-cage-exists} that each graph containing
  a~distance-$d$ semi-ladder of order $d(d\ell + 2)^d + 1$ also contains a quasi-cage
  of order $\ell$.
Then, we establish in Lemma \ref{cage-exists} that each graph containing
  a~quasi-cage of order $\Theta(d) \cdot \ell$ also contains a cage of order $\ell$.
  
In Subsection \ref{ordered-cages-section}, we prove multiple structural properties of cages.
These properties allow us to introduce two variants of a cage: ordered cage and identity ordered
  cage, which are structures utilizing the topology of the embedding of the graph in the plane.
We then prove that each graph with a cage of order $\ell$ also contains an ordered cage of
  order $\ell$ (Lemma~\ref{order-lemma}), and that each graph with an ordered cage
  of order $\Theta(\ell^2)$ contains an identity ordered cage of order $\ell$
  (Lemma \ref{identity-ordered-cage-exists}).

In Subsections \ref{neighbor-area-section} and \ref{neighbor-cages-section}, we find more
  structural properties of identity ordered cages, building on prior knowledge about cages.
Thanks to these properties, we establish another variant of a cage --- neighbor cage ---
  and we prove that each graph with an identity ordered cage of order
  $\Theta(d^7) \cdot \ell$ also contains a neighbor cage of order $\ell$
  (Lemma \ref{large-neighbor-lemma}).
The polynomial bound on neighborhood complexity in planar graphs
  becomes essential in this part of the proof.
We then proceed to prove a~set of properties of neighbor cages.

Finally, in Subsection \ref{separating-cages-lemma}, we introduce the final variant of a cage
  --- separating cage --- and prove that if a graph contains a neighbor cage of order
  $\Theta(d^2) \cdot \ell$, then it also contains a separating cage of order $\ell$
  (Lemma \ref{large-separating-cage-lemma}).
We finalize the proof by directly showing that each separating cage has order smaller
  than $2d + 5$ (Lemma \ref{no-large-separating-cage}).
Retracing all the steps, we can easily figure out the $d^{O(d)}$ upper bound on the maximum
  order of a distance-$d$ semi-ladder in the class of planar graphs (Corollary
  \ref{planar-upper-bound-cor}).

\vspace{1em}
Throughout this section, we consider a fixed planar graph $G$ together with its embedding
  in the plane; the embedding may be altered later in the proof.
We also fix the parameter $d \geq 1$.
  
\subsection{Quasi-cages and cages}
\label{quasi-cages-subsection}

We begin by introducing a prototype for cages --- quasi-cages.
In order to do this, we first need to define a special kind of a rooted tree.

\begin{definition}
Oriented paths $P_1, P_2, \dots, P_\ell$ in a planar graph $G$ create a \textbf{geodesic tree}
  $\Tc = \bigcup P_i = \mathsf{Tree}(P_1, P_2, \dots, P_\ell)$ of
  order $\ell$ rooted at a vertex $r$ if the following conditions hold:
\begin{itemize}
  \item the union of all the paths is an oriented subtree $\Tc$ of $G$,
  \item all paths $P_i$ are shortest paths in $G$ and have the same length,
    not exceeding $d$,
  \item each path $P_i$ originates at $r$ and terminates at a different leaf of the tree.
\end{itemize}

We say that $\Tc$ \textbf{avoids} $s$ if none of the paths contain $s$ as a vertex.
Also, if for all different $i, j$, the intersection of $P_i$ and $P_j$ is equal to $\{r\}$,
  we call $\Tc$ a \textbf{simple geodesic tree}.
\end{definition}

We also use the following notation:
  $\Root(\Tc)$ is the root $r$ of the geodesic tree, and
  $\Path(\Tc, i)$ is the oriented path $P_i$.
Moreover, for any two vertices $u, v \in \Tc$, we define
  $\lca_\Tc(u, v)$ as the lowest common ancestor of these two vertices in $\Tc$,
  where $\Tc$ is considered rooted at $r$.
We may write $\lca$ instead of $\lca_\Tc$ if the tree is clear from the context.
\medskip

\begin{figure}[h]
  \centering
  \input{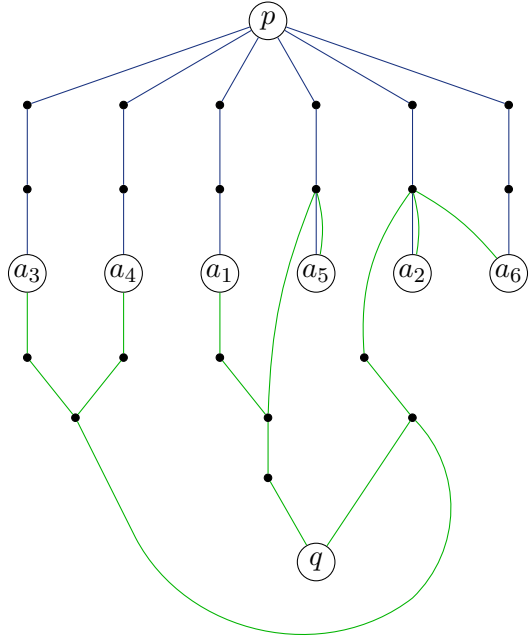}
  \caption{An example quasi-cage of order $6$. Tree $\Pc$ is marked blue and tree $\Qc$ is
    marked green. Vertices $b_1, b_2, \dots, b_6$ were omitted from the figure.
    Note that trees $\Pc$ and $\Qc$ can share vertices and edges, but $\Pc$ must not contain
      $q$, and $\Qc$ must not contain $p$.}
  \label{quasi-cage-sample-fig}
\end{figure}

We can now define a quasi-cage:

\begin{definition}
\label{def-quasi-cage}
A \textbf{quasi-cage} of order $\ell$ in a planar graph $G$ is a structure consisting of:
\begin{itemize}
\item a semi-ladder of order $\ell$ consisting of distinct vertices
  $a_1, a_2, \dots, a_\ell$, and $b_1, b_2, \dots, b_\ell$,
\item two distinct root vertices $p, q \not\in \{a_1, a_2, \dots, a_\ell, b_1, b_2, \dots, b_\ell\}$,
\item shortest paths $P_1, P_2, \dots, P_\ell$, where $P_i$ is an oriented path from $p$ to $a_i$ for
  each $i \in [1, \ell]$, which form a simple geodesic tree $\Pc = \Tree(P_1, P_2, \dots, P_\ell)$
  rooted at $p$ and avoiding $q$,
\item shortest paths $Q_1, Q_2, \dots, Q_\ell$, where $Q_i$ is a directed path from $q$ to $a_i$
  for each $i \in [1, \ell]$, which form a geodesic tree $\Qc = \Tree(Q_1, Q_2, \dots, Q_\ell)$
  rooted at $q$ and avoiding $p$.
\end{itemize}

We stress that while $\Pc$ must be a simple geodesic tree, $\Qc$ is not necessarily simple
  (Figure~\ref{quasi-cage-sample-fig}).

With the notation above, we denote a quasi-cage by
 $$\QuasiCage((a_1, a_2, \dots, a_\ell), (b_1, b_2, \dots, b_\ell), \Pc, \Qc). $$
\end{definition}

We remark that the definitions of semi-ladders, geodesic trees, simple
  geodesic trees, and quasi-cages are closed to taking subsets.
Formally, fix a subsequence of indices $i_1, i_2, \dots, i_m$ ($1 \leq i_1 < i_2 < \dots < i_m
  \leq \ell$).
Then:
\begin{itemize}
  \item For a semi-ladder $(a_1, a_2, \dots, a_\ell), (b_1, b_2, \dots, b_\ell)$,
    vertices
    $(a_{i_1}, a_{i_2}, \dots, a_{i_m})$, $(b_{i_1}, b_{i_2}, \dots, b_{i_m})$ form
    a semi-ladder of order $m$.
  \item For a geodesic tree $\mathsf{Tree}(P_1, P_2, \dots, P_\ell)$ of order $\ell$, also
    $\mathsf{Tree}(P_{i_1}, P_{i_2}, \dots, P_{i_m})$ is a geodesic tree of order $m$.
    If the original geodesic tree was simple, its subset is simple as well.
  \item For a quasi-cage of order $\ell$
    $$\QuasiCage((a_1, a_2, \dots, a_\ell), (b_1, b_2, \dots, b_\ell),
      \mathsf{Tree}(P_1, P_2, \dots, P_\ell), \mathsf{Tree}(Q_1, Q_2, \dots, Q_\ell)), $$
    the following describes a quasi-cage of order $m$:
    $$\QuasiCage((a_{i_1}, \dots, a_{i_m}), (b_{i_1}, \dots, b_{i_m}),
      \mathsf{Tree}(P_{i_1}, \dots, P_{i_m}), \mathsf{Tree}(Q_{i_1}, \dots, Q_{i_m})). $$
\end{itemize}

This fact will also be true for each further refinement of the definition of a quasi-cage.
Therefore, from now on, we will refer to the subsets of any object (semi-ladder, geodesic tree,
  quasi-cage and any its refinement) by taking the original object
  and specifying the indices $i_1, i_2, \dots, i_m$.

\vspace{1em}

As the first step of our proof, we will demonstrate that planar graphs with large distance-$d$
  semi-ladders also contain large quasi-cages.

\begin{lemma}
\label{quasi-cage-exists}
If $G$ contains a distance-$d$ semi-ladder of order $d(d\ell + 2)^d + 1$, then
  $G$ also contains a quasi-cage of order $\ell$.
  \begin{proof}
    We will state and prove a series of claims, each imposing more structure on the semi-ladder.
      
    We will call a distance-$d$ semi-ladder $(a_1, a_2, \dots, a_k), (b_1, b_2, \dots, b_k)$ of order $k$,
      together with another vertex $r$,
      \textbf{R-equidistant} if:
    \begin{itemize}
      \item $r \not\in \{a_1, a_2, \dots, a_k, b_1, b_2, \dots, b_k\}$,
      \item each vertex $a_i$ is at the same distance from $r$, not exceeding $d$.
    \end{itemize}
    We will denote such semi-ladder as $\RSemiLadder((a_1, \dots, a_k), (b_1, \dots, b_k),
      r)$.

    \begin{claim}
    \label{quasi-cage-exists-r-claim}
    Every distance-$d$ semi-ladder of order $d\ell + 1$ contains an R-equidistant semi-ladder of order $\ell$
      as a subset.
      \begin{claimproof}
      \cqed
      Take a distance-$d$ semi-ladder of order $d\ell + 1$ in $G$: $(a_1, a_2, \dots, a_{d\ell + 1}),
        (b_1, b_2, \dots, b_{d\ell + 1})$.
      Notice that $1 \leq \dist(a_i, b_1) \leq d$ for each $i \geq 2$.
      This means that the set $\{a_2, a_3, \dots, a_{d\ell + 1}\}$ of vertices contains
        a subset with at least $\frac{(d\ell + 1) - 1}{d} = \ell$ elements
        where each element is at the same distance from $b_1$, not exceeding $d$.
      Let this subset be $a_{i_1}, a_{i_2}, \dots, a_{i_\ell}$ ($2 \leq i_1 < i_2 < \dots < i_\ell \leq
        d\ell + 1$).
      Then, the following structure is an R-equidistant semi-ladder:
      $$ \RSemiLadder((a_{i_1}, \dots, a_{i_\ell}), (b_{i_1}, \dots, b_{i_\ell}), b_1). $$
      By our definition of a semi-ladder, the vertices of any semi-ladder must be different.
      Hence, $b_1 \not\in \{a_{i_1}, \dots, a_{i_\ell}, b_{i_1}, \dots, b_{i_\ell}\}$.
      \end{claimproof}
    \end{claim}
    
    A \textbf{simple geodesic semi-ladder} of order $k$ is a distance-$d$ semi-ladder
      $(a_1, a_2, \dots, a_k)$, $(b_1, b_2, \dots, b_k)$, together with a simple geodesic tree
      $\Pc = \mathsf{Tree}(P_1, P_2, \dots, P_k)$,
      where for each $i \in [1, k]$, the path $P_i$ terminates at $a_i$.
    We denote such a simple geodesic semi-ladder as
    $\SGSemiLadder((a_1, \dots, a_k), (b_1, \dots, b_k), \Pc)$.
    
    \begin{claim}
    \label{quasi-cage-exists-sg-claim}
    Every R-equidistant semi-ladder of order $\ell^d$ contains a
      simple geodesic semi-ladder of order $\ell$ as a subset.
      \begin{claimproof}
      \cqed
      Take an R-equidistant semi-ladder of order $\ell^d$:
        $$\RSemiLadder((a_1, \dots, a_{\ell^d}), (b_1, \dots, b_{\ell^d}), r).$$
      Let $\delta := \dist(r, a_1)$. We also take 
        $\Tc$ as a shortest paths tree rooted at $r$ whose leaves are exactly
        $a_1, a_2, \dots, a_{\ell^d}$; this tree exists since all hypothetical leaves
        are at the same distance $\delta \leq d$ from $r$.
      This tree has depth $\delta$ and $\ell^d$ leaves.
      Thus, $\Tc$ contains a vertex which has at least $\ell^{d / \delta} \geq \ell$ children.
      Let $p$ be any such vertex.
      
      Since $p$ has at least $\ell$ children in $\Tc$, and each child is a root of a subtree
        in $\Tc$ containing at least one leaf of $\Tc$,
        we can pick a subset $A = \{a_{i_1}, a_{i_2}, \dots, a_{i_\ell}\}$ ($1 \leq i_1 < i_2 < \dots <
          i_\ell \leq \ell^d$)
        of leaves of $\Tc$ where each leaf is in a different subtree of $\Tc$ rooted at a child of $p$.
      Also, for each $j \in [1, \ell]$, we define $P_j$ as the shortest oriented path in $\Tc$
        which originates at $p$ and terminates at $a_{i_j}$.
      By the choice of $A$, all paths $P_1, P_2, \dots, P_\ell$ are vertex-disjoint apart
        from their common origin~$p$.
      They also have the same length since all leaves of $\Tc$ are at the same depth,
        and all leaves in $A$ are in the subtree of $\Tc$ rooted at $p$.
      Finally, they are shortest paths in $G$ by the choice of $\Tc$ as a shortest paths tree.
      Therefore, $\mathsf{Tree}(P_1, P_2, \dots, P_\ell)$ is a simple geodesic tree.
      
      We can now see that the following structure is a simple geodesic semi-ladder:
      \[ \SGSemiLadder((a_{i_1}, \dots, a_{i_\ell}), (b_{i_1}, \dots, b_{i_\ell}),
        \mathsf{Tree}(P_1, P_2, \dots, P_\ell)). \]
      \end{claimproof}
    \end{claim}
    
    A simple geodesic semi-ladder
      $\SGSemiLadder((a_1, \dots, a_k), (b_1, \dots, b_k), \Pc)$ of order $k$,
      together with another vertex $q$, is \textbf{Q-equidistant} if:
    \begin{itemize}
      \item $\Pc$ avoids $q$,
      \item $q \not\in \{a_1, a_2, \dots, a_k, b_1, b_2, \dots, b_k\}$,
      \item $q$ is at the same distance from each vertex $a_i$, and
      \item $\dist(q, \Root(\Pc)) + \dist(\Root(\Pc), a_i) > d$ for each $i \in [1, k]$.
    \end{itemize}
    We use the following notation for a Q-equidistant simple geodesic semi-ladder:
      $\QSemiLadder(\allowbreak (a_1, \dots, a_k),\allowbreak (b_1, \dots, b_k),\allowbreak \Pc, q)$.

    \begin{claim}
      \label{quasi-cage-exists-q-claim}
      Every simple geodesic semi-ladder of order $d\ell + 2$ contains
        a Q-equidistant simple geodesic semi-ladder of order $\ell$
        as a subset.
      \begin{claimproof}
      \cqed
      Take a simple geodesic semi-ladder of order $d\ell + 2$:
      $$ \SGSemiLadder((a_1, \dots, a_{d\ell + 2}), (b_1, \dots, b_{d\ell + 2}), \Pc) \qquad
      \text{where }\Pc = \mathsf{Tree}(P_1, P_2, \dots, P_{d\ell + 2}).$$

      Notice that $b_1 \neq \Root(\Pc)$; otherwise, we would have
        $d < \dist(b_1, a_1) = \dist(\Root(\Pc), a_1) = \dist(\Root(\Pc), a_2) =
        \dist(b_1, a_2) \leq d$.
      Since paths $P_2, \dots, P_{d\ell + 2}$ are vertex-disjoint apart from $\Root(\Pc)$,
        this means that at most one of them passes through $b_1$.
      Hence, at least $d\ell$ such paths do not contain $b_1$.

      Moreover, for each $i \in [2, d\ell + 2]$, we have $\dist(a_i, b_1) \in [1, d]$.
      By the pigeonhole principle, we can fix $\ell$ indices $i_1, i_2, \dots, i_\ell$
        ($2 \leq i_1 < i_2 < \dots < i_\ell \leq d\ell + 2$) with the following properties:
      \begin{itemize}
        \item for each $j \in [1, \ell]$, the path $P_{i_j}$ does not contain $b_1$,
        \item $b_1$ is at the same distance from each vertex $a_{i_j}$ ($j \in [1, \ell]$).
      \end{itemize}
      
      We will prove that the following structure is a Q-equidistant simple geodesic semi-ladder:
      $$ \QSemiLadder((a_{i_1}, \dots, a_{i_\ell}), (b_{i_1}, \dots, b_{i_\ell}),
        \mathsf{Tree}(P_{i_1}, \dots, P_{i_\ell}), b_1). $$
      Since the choice of the indices $i_1, \dots, i_\ell$ already met the first and the third
        condition in the definition of a Q-equidistant simple geodesic semi-ladder,
        we only need to verify the remaining requirements.

      As any two vertices in the semi-ladder are different,
        we have $b_1 \neq a_{i_j}$ and $b_1 \neq b_{i_j}$ for each $j \in [1, \ell]$, which
        confirms the second condition.
      Moreover, for every $j \in [1, \ell]$, we have that
      \begin{equation*}
      \begin{split}
      \dist(b_1, \Root(\Pc)) + \dist(\Root(\Pc), a_{i_j}) &= 
            \dist(b_1, \Root(\Pc)) + \dist(\Root(\Pc), a_1) \geq \\ &\geq
            \dist(b_1, a_1) > d.
      \end{split}
      \end{equation*}
      All the required conditions are satisfied, so the structure defined above
        is a Q-equidistant simple geodesic semi-ladder.
      \end{claimproof}
    \end{claim}
    
    Now we will see that a large Q-equidistant simple geodesic semi-ladder in $G$
      directly witnesses the existence of a large quasi-cage in $G$.
    
    \begin{claim}
      \label{quasi-cage-exists-qc-claim}
      Every graph with a Q-equidistant simple geodesic semi-ladder of order $\ell$ also
        contains a quasi-cage of order $\ell$.
      \begin{claimproof}
        \cqed
        Fix a Q-equidistant simple geodesic semi-ladder as in the statement of the claim:
          $$\Lc = \QSemiLadder((a_1, \dots, a_\ell),\allowbreak (b_1, \dots, b_\ell), \Pc, q).$$
        Let also $\Sc$ denote the semi-ladder $a_1, \dots, a_\ell$, $b_1, \dots, b_\ell$ underlying
          $\Lc$, and set $p := \Root(\Pc)$.
        Let us also take $\Qc$ as any shortest-path tree rooted at $q$ whose leaves
          are $a_1, a_2, \dots, a_\ell$.
          
        We will now verify that according to Definition \ref{def-quasi-cage}, the following structure
          is a quasi-cage:
        $$ \QuasiCage((a_1, \dots, a_\ell), (b_1, \dots, b_\ell), \Pc, \Qc). $$
        
        We first note that $p \not\in \Sc$ since $\Lc$ is also an R-equidistant semi-ladder.
        Similarly, $q$ is disjoint with $\Sc$, which follows immediately from the definition of $\Lc$ as
          a Q-equidistant simple geodesic semi-ladder.
        The geodesic tree $\Pc$ is simple (since $\Lc$ is a simple geodesic semi-ladder)
          and avoids $q$ (since $\Lc$ is a Q-equidistant simple geodesic semi-ladder).

        We only need to prove that $\Qc$ avoids $p$.
        Assume for contradiction that for some $i \in [1, \ell]$, a path $Q_i = \Path(\Qc, i)$
          contains $p$ as a vertex.
        This would mean that the shortest path between $q$ and $a_i$ contains
          $p$ as a vertex.
        However, by the definition of a Q-equidistant simple geodesic semi-ladder:
        $$ \dist(q, a_i) = \dist(q, p) + \dist(p, a_i) = \dist(q, \Root(\Pc)) + \dist(\Root(\Pc), a_i) > d $$
        --- a contradiction.
      \end{claimproof}
    \end{claim}

    We can now combine the claims above to prove our lemma.
    Take any distance-$d$ semi-ladder $L_1$ in $G$ of order $d(d\ell + 2)^d + 1$.

    By Claim \ref{quasi-cage-exists-r-claim}, $L_1$ contains an R-equidistant semi-ladder
      $L_2$ of order $(d\ell + 2)^d$ as a subset.

    By Claim \ref{quasi-cage-exists-sg-claim}, $L_2$ contains a simple geodesic semi-ladder
      $L_3$ of order $d\ell + 2$ as a subset.
      
    By Claim \ref{quasi-cage-exists-q-claim}, $L_3$ contains a Q-equidistant
      simple geodesic semi-ladder $L_4$ of order $\ell$ as a subset.
      
    Hence, by Claim \ref{quasi-cage-exists-qc-claim}, $L_4$ asserts the existence of
      a quasi-cage of order $\ell$ in $G$.
  \end{proof}
\end{lemma}

Next, we refine the definition of a quasi-cage by defining cages --- variants of
  quasi-cages ensuring that both geodesic trees in a definition of a quasi-cage intersect
  in a controlled way.

\begin{definition}
A \textbf{cage} of order $\ell$ is a quasi-cage
  $\QuasiCage((a_1, \dots, a_\ell), (b_1, \dots, b_\ell), \Pc, \Qc)$
  in which for every two indices $i, j \in [1, \ell]$, the paths $\Path(\Pc, i)$ and $\Path(\Qc, j)$
  have a non-empty intersection if and only if $i = j$.

We denote this structure in the following way:
$$ \Cage((a_1, \dots, a_\ell), (b_1, \dots, b_\ell), \Pc, \Qc). $$
\end{definition}

\begin{figure}[h]
  \centering
  \input{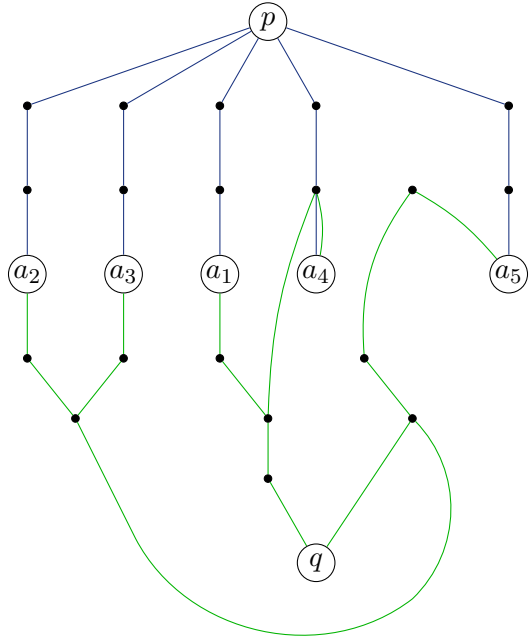}
  \caption{An example cage.
    This cage is a subset of the quasi-cage from Figure \ref{quasi-cage-sample-fig}.
    The pair of vertices $(a_2, b_2)$ has been removed from the original semi-ladder,
      and the remaining indices have been renumbered.
    Analogously to Figure \ref{quasi-cage-sample-fig}, vertices $b_1, b_2, \dots, b_5$ were
      removed from the picture.}
  \label{cage-sample-fig}
\end{figure}

To see that the definition of a cage is slightly stronger than the one of a quasi-cage,
  we can verify that a quasi-cage in Figure \ref{quasi-cage-sample-fig} is not a cage ---
  its path connecting $q$ with $a_6$ intersects the path connecting $p$ with $a_2$.
However, a large subset of it is a cage (Figure \ref{cage-sample-fig}).
This is not a coincidence, as proved in the following lemma.

\begin{lemma}
\label{cage-exists}
Every graph with a quasi-cage of order $(2d - 1)\ell$
  contains a cage of order $\ell$ as a subset.
  
  \begin{proof}
    Fix a quasi-cage:
    $$\QuasiCage((a_1, a_2, \dots, a_{(2d - 1)\ell}),
      (b_1, b_2, \dots, b_{(2d - 1)\ell}), \Pc, \Qc)$$
    where $\Pc = \mathsf{Tree}(P_1, \dots, P_{(2d-1)\ell})$,
      $\Qc = \mathsf{Tree}(Q_1, \dots, Q_{(2d-1)\ell})$.
  
	Firstly, note that for every $i \in [1, (2d-1)\ell]$, paths $P_i$ and $Q_i$ intersect, because
	  they share their final vertex $a_i$.
	Our aim is now to find a subset of indices $I \subseteq \{1, 2, \dots, (2d - 1)\ell\}$
	  such that $|I| = \ell$ and for any two different indices $i, j \in I$, paths
	  $P_i$ and $Q_j$ do not intersect.
	As soon as we achieve this, the subset of the quasi-cage given by the set of indices $I$
	  will form the sought cage.

    Let us create an auxiliary directed graph $H$ where $V(H) = \{1, 2, \dots, (2d - 1)\ell\}$,
      such that $H$ contains an oriented edge $i \to j$ for $i \neq j$ if and only if
      $P_i$ and $Q_j$ intersect.
    Let also $\widehat{H}$ be the undirected graph underlying $H$.
    We plan to find an independent set $I \subseteq V(\widehat{H})$ of size at least $\ell$.
    Then for any two indices $i, j \in I$, we will have that paths $P_i$ and $Q_j$ are disjoint,
      as well as paths $P_j$ and $Q_i$.
    Therefore, $I$ will be a subset of indices forming a cage.
              
    We note that each path $Q_j$ has length at most $d$.
    As $Q_j$ does not pass through $p$ (the root of the simple geodesic tree $\mathcal{P}$),
      each vertex other than its origin $q$ can lie on at most one of the paths $P_i$.
    Thus, $Q_j$ intersects at most $d$ different paths in $\mathcal{P}$.
    Since $Q_j$ intersects with $P_j$ for every $j \in [1, (2d-1)\ell]$, we get that
      $Q_j$ has a non-empty intersection with at most $d-1$ paths $P_i$ for $i \neq j$.
    This means that no vertex of $H$ has in-degree larger than $d-1$.

    In every non-empty subset $S$ of vertices of $\widehat{H}$,
      the total number of edges in the subgraph
      of $\widehat{H}$ induced by $S$ does not exceed $|S| \cdot (d-1)$;
      this is because in the subgraph of $H$ induced by $S$, each vertex has in-degree not
      exceeding $d-1$.
    Thus, $\Sc$ contains a vertex with degree not exceeding the average degree of the vertices
      in $S$, which can be bounded from above by
      $$\frac{2 \cdot [|S| \cdot (d - 1)]}{|S|} = 2d - 2.$$
    Hence, $\widehat{H}$ is $(2d-2)$-degenerate \cite{lick_white_1970}.
    Since $(2d-2)$-degenerate graphs are $(2d-1)$-colorable~\cite{SZEKERES19681},
      $\widehat{H}$ contains an independent set of size at least
      \[\frac{|V(\widehat{H})|}{2d-1} = \ell.\]
  \end{proof}
\end{lemma}

\subsection{Ordered cages}
\label{ordered-cages-section}

So far, the definition of a cage does not involve any topological properties of the plane.
In this subsection, we will introduce the topology to the cages by defining
  ordered cages.
We will, however, need some notation beforehand.

\begin{definition}
  For a cage $\Cc := \Cage((a_1, \dots, a_\ell), (b_1, \dots, b_\ell), \Pc, \Qc)$
    of order $\ell$ in a graph $G$, and an index $i \in \{1, 2, \dots, \ell\}$, we define:
    
    \begin{itemize}
      \item the \textbf{$i$-th splitting vertex} $\mathsf{Split}(\Cc, i)$
        as the vertex of the intersection of $\Path(\Pc, i)$ and $\Path(\Qc, i)$
        which is the closest to $\Root(\Pc)$ (this vertex is unique since $\Path(\Pc, i)$ is
        a shortest path in $G$);
      \item the \textbf{$i$-th splitting path} $\SPath(\Cc, i)$
        as the oriented simple path originating at $\Root(\Pc)$,
        following the prefix of $\Path(\Pc, i)$ until $\mathsf{Split}(\Cc, i)$,
        and then --- the prefix of $\Path(\Qc, i)$ in the opposite direction.
        Formally,
        $$ \SPath(\Cc, i) = \Path(\Pc, i)[\Root(\Pc),
          \mathsf{Split}(\Cc, i)]\ \cdot\ (\Path(\Qc, i))^{-1}[\mathsf{Split}(\Cc, i),
          \Root(\Qc)].$$
    \end{itemize}
\end{definition}

\begin{figure}[h]
  \centering
  \input{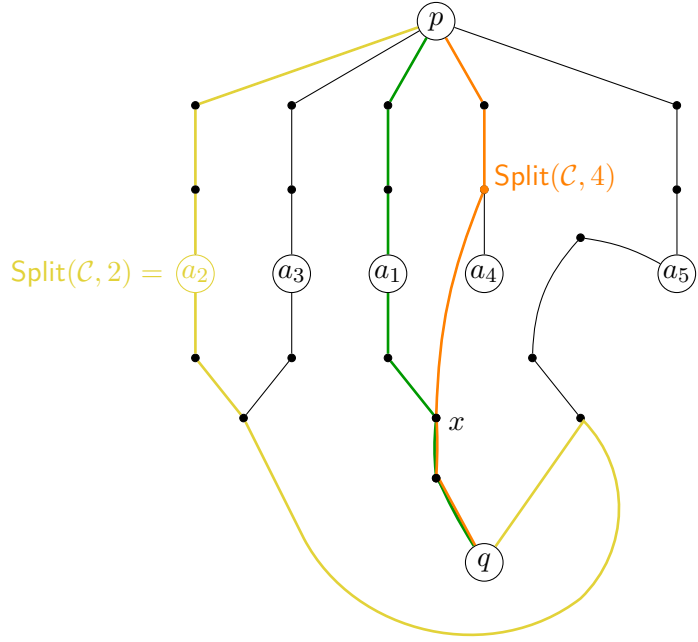}
  \caption{Example splitting vertices and splitting paths in the cage $\Cc$ from
    Figure~\ref{cage-sample-fig}.
    Objects $\mathsf{Split}(\Cc, 2)$ and $\SPath(\Cc, 2)$ are marked yellow,
    objects $\mathsf{Split}(\Cc, 4)$ and $\SPath(\Cc, 4)$ are marked orange,
    and $\SPath(\Cc, 1)$ is marked green.}
  \label{splitting-sample-fig}
\end{figure}

The definitions are also explained in Figure \ref{splitting-sample-fig}.

We derive the following facts about splitting paths:

\begin{lemma}
\label{two-splitting-paths-lemma}
  Fix a cage $\Cc = ((a_i)_{i=1}^\ell, (b_i)_{i=1}^\ell, \Pc, \Qc)$,
    two different indices $i, j \in [1, \ell]$ and a vertex $x \in V(\Qc)$.
  If $x \in Q_i$ and $x \in Q_j$, then $x \in \SPath(\Cc, i) \cap \SPath(\Cc, j)$.
  \begin{proof}
    Note that $\Split(\Cc, i) \in P_i \cap Q_i$.
    Hence, path $Q_j$ cannot contain $\Split(\Cc, i)$ or otherwise
      paths $P_i$ and $Q_j$ would intersect.
    Therefore, $x$ is an ancestor of $\Split(\Cc, i)$ in $\Qc$, so
      $$x \in (\Path(\Qc, i))^{-1}[\Split(\Cc, i), \Root(\Qc)].$$
    Thus, $x \in \SPath(\Cc, i)$.
      
    We prove that $x \in \SPath(\Cc, j)$ analogously.
  \end{proof}
\end{lemma}

\begin{lemma}
  \label{splitting-path-intersection}
  For a cage $\Cc = ((a_i)_{i=1}^\ell, (b_i)_{i=1}^\ell, \Pc, \Qc)$ in a graph $G$,
    and two different indices $i, j \in [1, \ell]$:
  $$ \SPath(\Cc, i) \cap \SPath(\Cc, j)\ =\ \{\Root(\Pc)\} \cup
    \SPath(\Cc, i)[\lca_{\Qc}(a_i, a_j), \Root(\Qc)]. $$
  In other words, any two distinct paths $\SPath(\Cc, i)$, $\SPath(\Cc, j)$
    intersect exactly in their first vertex $p$ and in their common suffix.
    
  \begin{proof}
    We first note that the following condition is simple to verify.
    By Lemma \ref{two-splitting-paths-lemma}, we have that
      $$\SPath(\Cc, i)[\lca_{\Qc}(a_i, a_j), \Root(\Qc)] =
        \SPath(\Cc, j)[\lca_{\Qc}(a_i, a_j), \Root(\Qc)].$$
    Hence:
      $$ \SPath(\Cc, i) \cap \SPath(\Cc, j)\ \supseteq\ \{\Root(\Pc)\} \cup
    \SPath(\Cc, i)[\lca_{\Qc}(a_i, a_j), \Root(\Qc)]. $$
    We only need to prove the $\subseteq$ part of the equality.
  
    Let $\Pc = \mathsf{Tree}(P_1, P_2, \dots, P_\ell)$ and $\Qc = \mathsf{Tree}(Q_1, Q_2, \dots,
      Q_\ell)$, and fix two different indices $i, j$.
    By the definition of a cage, paths $P_i$ and $Q_j$ are vertex-disjoint,
      and paths $P_j$ and $Q_i$ are vertex-disjoint.
    Since $\SPath(\Cc, i) \subseteq P_i \cup Q_i$ and $\SPath(\Cc, j) \subseteq P_j \cup Q_j$, we
      conclude that
    $$ \SPath(\Cc, i) \cap \SPath(\Cc, j) \subseteq (P_i \cap P_j) \cup (Q_i \cap Q_j). $$
    
    We conclude the proof by noting that $P_i \cap P_j = \{p\}$ ($\Pc$ is a simple geodesic tree)
      and $Q_i \cap Q_j = \SPath(\Cc, i)[\lca_\Qc(a_i, a_j), \Root(\Qc)]$.
  \end{proof}
\end{lemma}

In the example in Figure \ref{splitting-sample-fig}, paths $\mathsf{Split}(\Cc, 2)$
  and $\mathsf{Split}(\Cc, 4)$ intersect at $p$ and $q$ only,
whereas paths $\mathsf{Split}(\Cc, 1)$ and $\mathsf{Split}(\Cc, 4)$ intersect
  at $p$ and their common suffix, $\mathsf{Split}(\Cc, 4)[x, q]$
  where $x = \lca_\Qc(a_1, a_4)$.
  
\vspace{1em}
We will now introduce the concept of an order of a cage $\Cc$.
Roughly speaking, if we let $\Cc = \Cage((a_i)_{i=1}^\ell, (b_i)_{i=1}^\ell, \Pc, \Qc)$, then
  the paths $\Path(\Pc, 1), \Path(\Pc, 2), \dots, \Path(\Pc, \ell)$ originate at the common vertex
  $\Root(\Pc)$ and are vertex disjoint otherwise.
Therefore, we can order all these paths cyclically around $\Root(\Pc)$.
We stress that this cyclic order might not necessarily be consistent with
  the numbering of the paths.
Fortunately, it turns out that we can order these paths in such a way that
  the paths $\Path(\Qc, 1), \Path(\Qc, 2), \dots, \Path(\Qc, \ell)$ are also sorted
  in a meaningful way.
This is formalized by the following definition:

\begin{definition}
\label{order-def}
In a cage $\Cc = \Cage((a_i)_{i=1}^\ell, (b_i)_{i=1}^\ell, \Pc, \Qc)$ of order $\ell$ in a graph $G$
  embedded in the plane, where $\Pc = \Tree(P_1, \dots, P_\ell)$,
  $\Qc = \Tree(Q_1, \dots, Q_\ell)$,
  an \textbf{order} is a permutation $\sigma = (\sigma_1, \sigma_2, \dots, \sigma_\ell)$ of
  $\{1, 2, \dots, \ell\}$ such that:
\begin{itemize}
  \item the first edges on paths $P_{\sigma_1}, P_{\sigma_2}, \dots, P_{\sigma_\ell}$
    are ordered anti-clockwise around $\Root(\Pc)$ in the embedding of $G$ in the plane, and
  \item for each vertex $v \in V(\mathcal{Q})$, there exist two indices $L, R$
    ($1 \leq L \leq R \leq \ell$) such that
    $v$ lies on paths $Q_{\sigma_L}, Q_{\sigma_{L+1}}, \dots, Q_{\sigma_R}$ and on no other path
    $Q_i$.
\end{itemize}
\end{definition}

For instance, the cage in Figure \ref{splitting-sample-fig} has an order; permutations
  $\sigma = (1, 4, 5, 2, 3)$ and $\sigma = (5, 2, 3, 1, 4)$ are example orders.
However, not all cyclic shifts of $(1, 4, 5, 2, 3)$ are correct orders --- for instance,
  $\sigma = (4, 5, 2, 3, 1)$ is \emph{not} an order because the vertex marked $x$ in
  Figure \ref{splitting-sample-fig} would belong to paths $P_{\sigma_1}$ and $P_{\sigma_5}$
  and no other paths.

Using Definition \ref{order-def}, we can naturally refine the definition of a cage.

\begin{definition}
\label{ordered-cage-def}
An \textbf{ordered cage} of order $\ell$ in a graph $G$ embedded in the plane
  is a cage together with its order $\sigma$.  
We denote this structure as
$$ \OrderedCage((a_1, \dots, a_\ell), (b_1, \dots, b_\ell), \Pc, \Qc, \sigma). $$
\end{definition}

We are now ready to present a lemma which allows us to turn each cage into an ordered cage:

\begin{lemma}
  \label{order-lemma}
  Each cage in a graph $G$ embedded in the plane has an order.
  
  \begin{proof}  
  Fix a cage $\Cc = ((a_i)_{i=1}^\ell, (b_i)_{i=1}^\ell, \Pc, \Qc)$ in $G$.
  We let $\Pc = \mathsf{Tree}(P_1, \dots, P_\ell)$, $\Qc = \mathsf{Tree}(Q_1, \dots, Q_\ell)$,
    $p = \Root(\Pc)$, $q = \Root(\Qc)$.
    
  If $\ell = 1$, the lemma is trivial.
  From now on, we assume that $\ell \geq 2$.

  Consider a permutation $\mu^1 = (\mu^1_1, \mu^1_2, \dots, \mu^1_\ell)$ of
    $\{1, 2, \dots, \ell\}$ for which the first edges on paths $P_{\mu^1_1}, P_{\mu^1_2}, \dots,
    P_{\mu^1_\ell}$ are ordered anti-clockwise around $p$.
  We also consider all cyclic shifts of $\mu^1$ --- for each $k \in [1, \ell]$, we define
    $\mu^k = (\mu^k_1, \mu^k_2, \dots, \mu^k_\ell)$ as follows:
  
  $$ \mu^k = (\mu^1_k, \mu^1_{k+1}, \dots, \mu^1_\ell, \mu^1_1, \mu^1_2, \dots, \mu^1_{k-1}). $$
  
  Clearly, for each $k \in [1, \ell]$, the first edges on the paths
    $P_{\mu^k_1}, P_{\mu^k_2}, \dots, P_{\mu^k_\ell}$ are ordered anti-clockwise around $p$.

  We will prove that at least one of the permutations $\mu^1, \mu^2, \dots, \mu^\ell$ is an order.
  We begin by taking a vertex $v \in V(\Qc)$ and describing the set
    of paths $Q_i$ containing $v$.
  
  \begin{claim}
    \label{order-lemma-claim-no-interleave}
    For every $v \in V(\Qc)$,
    no cyclic shift $\mu^k$ ($k \in [1,\ell]$) of $\mu^1$ contains four distinct indices
    $\mu^k_x, \mu^k_y, \mu^k_z, \mu^k_t$ ($1 \leq x < y < z < t \leq \ell$) such that
    $$ v \in Q_{\mu^k_x}, \quad v\not\in Q_{\mu^k_y}, \quad v \in Q_{\mu^k_z}, \quad
      v \not\in Q_{\mu^k_t}. $$
    \begin{claimproof}
    \cqed
    For the sake of contradiction let us pick a vertex $v \in V(\Qc)$, an index $k \in [1, \ell]$,
      and indices $1 \leq x < y < z < t \leq \ell$ for which the conditions above hold.
    We simplify the notation by setting $\ddot{x} = \mu^k_x$, $\ddot{y} = \mu^k_y$,
      $\ddot{z} = \mu^k_z$, $\ddot{t} = \mu^k_t$ (Figure \ref{no-interleave-claim-fig}(a)).

    We first remark that since $q$ is the root of $\Qc$, it belongs to all paths $Q_i$.
    Therefore, we must have that $v \neq q$.
    
    \begin{figure}
      \centering
      \begin{minipage}[b]{0.33\textwidth}
        \centering
        \input{figures/ordering-a.tex}
        (a)
      \end{minipage}\begin{minipage}[b]{0.33\textwidth}
        \centering
        \input{figures/ordering-b.tex}
        (b)
      \end{minipage}\begin{minipage}[b]{0.33\textwidth}
        \centering
        \input{figures/ordering-c.tex}
        (c)
      \end{minipage}
      \caption{(a) --- the setup in the proof of Claim \ref{order-lemma-claim-no-interleave}. \\
        (b) --- the edges added in $G'$ (green). \\
        (c) --- sets {$\color{green!70!white}A_1$},
          {$\color{cyan!90!white}A_2$},
          {$\color{pink!90!black}A_3$},
          {$\color{brown!90!black}A_4$},
          {$\color{violet!80!white}A_5$}
          of vertices used in locating the $K_5$ clique as a minor of $G'$. }
      \label{no-interleave-claim-fig}
    \end{figure}

    We will create another planar graph $G'$ which is a version of $G$ with some edges added
      and some edges removed.
    Let $f_{\ddot{x}}$, $f_{\ddot{y}}$, $f_{\ddot{z}}$, $f_{\ddot{t}}$ be the second
      vertices on each of the paths $P_{\ddot{x}}$, $P_{\ddot{y}}$,
      $P_{\ddot{z}}$, and $P_{\ddot{t}}$, respectively --- that is, the vertices immediately
      following $p$ on the corresponding paths.
    We create a graph $G'$ from $G$ with the same set of vertices and the same
      embedding in the plane, but with the following
      modifications (Figure \ref{no-interleave-claim-fig}(b)):
    \begin{itemize}
      \item We disconnect $p$ (the root of $\mathcal{P}$) from all its neighbors apart
        from $f_{\ddot{x}}$, $f_{\ddot{y}}$, $f_{\ddot{z}}$, $f_{\ddot{t}}$.
      \item We add the following edges to $G'$: $f_{\ddot{x}}f_{\ddot{y}}$,
        $f_{\ddot{y}}f_{\ddot{z}}$, $f_{\ddot{z}}f_{\ddot{t}}$, and
        $f_{\ddot{t}}f_{\ddot{x}}$.
    \end{itemize}
    Removing edges from the graph does not spoil its planarity.
    Afterwards, four edges are inserted to $G'$.
    In this graph, the anti-clockwise cyclic order of all edges incident to $p$ is
      $pf_{\ddot{x}}, pf_{\ddot{y}}, pf_{\ddot{z}}, pf_{\ddot{t}}$.
    Each new edge $f_if_j$ connects a neighbor of $p$ with its next neighbor
      in this cyclic order.
    In order to preserve the planarity of $G'$, we draw an edge $f_if_j$ so that
      the triangle created by oriented edges $f_if_j$, $f_jp$, $pf_i$ is oriented anti-clockwise,
      and $f_if_j$ lies very close to edges $pf_i$ and $pf_j$ in the embedding
      (Figure \ref{no-interleave-claim-fig}(b)).
    Since $p$ has only four neighbors in $G'$, we can draw each edge without breaking
      the planarity of the embedding of $G'$.
    
    \vspace{0.5em}
    We will now find the clique $K_5$ as a minor of $G'$.
    We define five subgraphs of $G'$, induced by the following sets of vertices
      (Figure \ref{no-interleave-claim-fig}(c)):
    \begin{equation*}
      \begin{split}
      A_1 &= \{p\}, \\
      A_2 &= \SPath(\Cc, \ddot{x})[f_{\ddot{x}}, v], \\
      A_3 &= \SPath(\Cc, \ddot{y})[f_{\ddot{y}}, q], \\
      A_4 &= \SPath(\Cc, \ddot{z})[f_{\ddot{z}}, v]\ \setminus\ A_2, \\
      A_5 &= \SPath(\Cc, \ddot{t})[f_{\ddot{t}}, q]\ \setminus\ A_3. \\
      \end{split}
    \end{equation*}
    
    Since $v \in Q_{\ddot{x}}$ and $v \in Q_{\ddot{z}}$,
      Lemma \ref{two-splitting-paths-lemma} asserts that
      $v \in \SPath(\Cc, \ddot{x})$ and $v \in \SPath(\Cc, \ddot{z})$.
    It is also obvious that all other ends of paths in the definitions of $A_2, A_3, A_4, A_5$
      belong to the respective paths.
    Hence, the sets above are defined correctly.
    These sets are also non-empty since
    
    $$ p \in A_1,\ \ f_{\ddot{x}} \in A_2,\ \ f_{\ddot{y}} \in A_3,\ \ f_{\ddot{z}} \in A_4,\ \ 
      f_{\ddot{t}} \in A_5. $$
    
    We now prove that each of the induced subgraphs $G'[A_1], G'[A_2], G'[A_3], G'[A_4], G'[A_5]$
      is connected.
    This fact is obvious for $A_1$, $A_2$, and $A_3$.
    For $A_4$, we need to observe that because of Lemma \ref{splitting-path-intersection},
      paths $\SPath(\Cc, \ddot{x})$ and $\SPath(\Cc, \ddot{z})$
      intersect exactly at $p$ and their common suffix.
    Therefore, removing $A_2$ from $\SPath(\Cc, \ddot{z})[f_{\ddot{z}}, v]$
      erases a suffix from the path, and leaves a subpath as a result.
    An analogous argument applies to $A_5$.
    Therefore, each of the subgraphs is connected.
    
    \vspace{0.5em}
    We can also prove the pairwise disjointness of $A_1, \dots, A_5$.
    Obviously, $A_1$ is disjoint with all the remaining sets $A_2, A_3, A_4, A_5$.
    Also, $A_2$ is disjoint with $A_4$ by the definition of $A_4$,
    and $A_3$ is disjoint with $A_5$ by the definition of $A_5$.
    Now, take any indices $i \in \{\ddot{x}, \ddot{z}\}$, $j \in \{\ddot{y}, \ddot{t}\}$, and let us prove
      that paths $\SPath(\Cc, i)[f_i, v]$ and $\SPath(\Cc, j)[f_j, q]$
      are vertex-disjoint.
    By Lemma \ref{splitting-path-intersection} we have that
    \begin{equation}
    \label{k5-spath-intersection-eq}
    \SPath(\Cc, i) \cap \SPath(\Cc, j) = \{p\} \cup \SPath(\Cc, i)[\lca_\Qc(a_i, a_j), q]
    \end{equation}
    (where $\SPath(\Cc, i)[\lca_\Qc(a_i, a_j), q]$ is the common suffix of $\SPath(\Cc, i)$ and
      $\SPath(\Cc, j)$).

    Naturally, $p \not\in \SPath(\Cc, i)[f_i, v]$.
    Moreover, since $v \not\in Q_j$, we infer that
      $v \not\in \SPath(\Cc, i)[\lca_\Qc(a_i, a_j), q]$,
      so $\SPath(\Cc, i)[f_i, v]$ is also disjoint with
      $\SPath(\Cc, i)[\lca_\Qc(a_i, a_j), q]$.
    Hence, by \ref{k5-spath-intersection-eq}, $\SPath(\Cc, i)[f_i, v]$ is disjoint with
      $\SPath(\Cc, j)$ as well.
    This means that we proved the vertex-disjointness of each pair of sets $A_1, \dots, A_5$.
    
    \vspace{0.5em}
    We are left to see that there is an edge between each pair of sets:
    
    \begin{itemize}
      \item There is a direct edge $pf_{\ddot{x}}$ connecting sets $A_1$ and $A_2$;
        we analogously show edges connecting $A_1$ with each remaining set.
      \item There is a direct edge $f_{\ddot{x}}f_{\ddot{y}}$ connecting sets $A_2$ and $A_3$;
        we analogously show edges connecting $A_3$ with $A_4$, $A_4$ with $A_5$,
        and $A_5$ with $A_2$.
      \item Since the sets $A_2 = \SPath(\Cc, \ddot{x})[f_{\ddot{x}}, v]$ and
        $\SPath(\Cc, \ddot{z})[f_{\ddot{z}}, v]$ induce connected graphs and
        have a non-empty intersection, then
        naturally $A_4 = \SPath(\Cc, \ddot{z})[f_{\ddot{z}}, v] \setminus A_2$ is connected
        by an edge with $A_2$.
        We analogously show that $A_3$ and $A_5$ are connected by an edge.
    \end{itemize}
    \vspace{0.5em}
    Summing everything up, we created a modified planar graph $G'$ and indicated
      five non-empty connected subgraphs such that each pair of subgraphs is vertex-disjoint,
      yet there exist edges connecting each pair of subgraphs.
    Therefore, $G'$ contains $K_5$ as a minor and hence --- by Wagner's theorem
      \cite{agnarsson2007graph} --- it cannot be planar.
    The contradiction concludes the proof of the claim.
    \end{claimproof}
  \end{claim}
  
  As a direct consequence of Claim \ref{order-lemma-claim-no-interleave},
    we can see that no cyclic shift $\mu^k$ of $\mu^1$ can contain four direct indices
    $\mu^k_x, \mu^k_y, \mu^k_z, \mu^k_t$ (defined as in Claim
    \ref{order-lemma-claim-no-interleave}) such that
  $$ v \not\in Q_{\mu^k_x}, \quad v\in Q_{\mu^k_y}, \quad v \not\in Q_{\mu^k_z}, \quad
       v \in Q_{\mu^k_t}. $$
  If this set of conditions was true, we could consider the cyclic shift of $\mu^k$
    with $\mu^k_y$ at the first position.
  We can easily see that this cyclic shift would be contradictory with
    Claim \ref{order-lemma-claim-no-interleave}.
  
  \vspace{1em}
  For a vertex $v \in V(\mathcal{Q})$ and integer $k \in [1, \ell]$,
    we set $\Kc^k_v := \{i \in [1, \ell]\, \mid\, 
    v \in Q_{\mu^k_i}\}$; in other words, $\Kc^k_v$ is the set of positions $i$
    in the permutation $\mu^k$ for which the corresponding path $Q_{\mu_i^k}$
    contains $v$ as a vertex.

  Note that for every $k_1, k_2 \in [1, \ell]$ and every vertex $v \in V(\Qc)$,
    set $\Kc^{k_1}_v$ is a cyclic shift of $\Kc^{k_2}_v$; this is because the permutations
    $\mu^{k_1}$ and $\mu^{k_2}$ are also mutual cyclic shifts.
    
  \begin{claim}
    \label{order-lemma-cyclic-interval}
    Each set $\Kc^k_v$ is a cyclic interval in $[1, \ell]$; that is, it is either an interval or
      a set of the form $[x, \ell] \cup [1, y]$ for some $x > y$.
    \begin{claimproof}
    \cqed
    Fix a vertex $v \in V(\mathcal{Q})$, integer $k \in [1, \ell]$,
      and the corresponding set $\Kc^k_v$.
    Assume that $\Kc^k_v$ is not a (non-cyclic) interval in $[1, \ell]$.
    In this setup, there must exist
      three indices $x, y, z$ ($1 \leq x < y < z \leq \ell$) for which
      $x \in \Kc^k_v$, $y \not\in \Kc^k_v$, and $z \in \Kc^k_v$.
    Without loss of generality, we can assume that $[x+1, z-1]$ is disjoint with $\Kc^k_v$.

    Note that there cannot exist an index $t > z$ such that $t \not\in \Kc^k_v$ as this would
      directly contradict Claim \ref{order-lemma-claim-no-interleave}.
    Therefore, we have $[z, \ell] \subseteq \Kc^k_v$.
    
    Also, there cannot exist an index $t < x$ such that $t \not\in \Kc^k_v$ as this would
      pose a contradiction with the direct consequence of Claim
      \ref{order-lemma-claim-no-interleave}.
    Thus, we have $[1, x] \subseteq \Kc^k_v$.

    Since $[1,x] \subseteq \Kc^k_v$, $[x+1, z-1] \cap \Kc^k_v = \varnothing$, and
      $[z, \ell] \subseteq \Kc^k_v$, we have that $\Kc^k_v = [1, x] \cup [z, \ell]$.
    \end{claimproof}
  \end{claim}
  
  In order to complete our proof, we also need to know how the cyclic intervals
    $\Kc^k_v$ interact with each other.
  It turns out that these intervals form a laminar family:
  
  \begin{claim}
    \label{order-lemma-laminar-family}
    For any two vertices $u, v \in V(\mathcal{Q})$ and $k \in [1, \ell]$, the sets $\Kc^k_u$ and
      $\Kc^k_v$ are either disjoint, or one of them is a subset of the other.
    \begin{claimproof}
    \cqed
    We remark that for $v \in V(\mathcal{Q})$, we have the following equivalent definition
      of $\Kc^k_v$:
    $$ \Kc^k_v = \{i \in [1, \ell]\, \mid\, a_{\mu^k_i}\text{ is in the subtree of }
      \Qc\text{ rooted at }v\}. $$
    This is true since the path $Q_{\mu^k_i}$ in the original definition of $\Kc^k_v$
      originates at $\Root(\Qc)$ and terminates at $a_{\mu^k_i}$ (a leaf in $\Qc$).
    Hence, $v$ belongs to this path if and only if its endpoint --- $a_{\mu^k_i}$ ---
      is in the subtree rooted at $v$.

    It is however well-known that the family of sets of leaves in every rooted tree $T$, where each
      set corresponds to a subset of leaves in a rooted subtree of $T$, is laminar.
    The proof of the claim follows.
    \end{claimproof}
  \end{claim}
  
  Having proved all the necessary observations, we can find the cyclic shift $\mu^k$ of $\mu^1$
    which forms an order in the embedding of the graph.
  Fix any vertex $v^* \in V(\mathcal{Q})$ for which the cyclic interval $\Kc^1_{v^*}$
    is not equal to $[1, \ell]$, but otherwise contains the largest number of indices.
  Such a vertex exists since $\ell \geq 2$, so the leaves in $\Qc$ do not cover the whole interval
    $[1, \ell]$.

  We now take the cyclic shift $\mu^k$ for which $\Kc^k_{v^*}$ is a prefix of $[1, \ell]$;
    this is possible since $\Kc_{v^*}$ is a cyclic interval in $\mu^1$.
    
  \begin{claim}
    \label{order-lemma-found-claim}
    $\mu^k$ is an order in $\Cc$.
    \begin{claimproof}
    \cqed
    Let $M$ be an integer such that $\Kc^k_{v^*} = [1, M]$.
    We will now prove that $\mu^k$ is an order.
    Firstly, since $\mu^k$ is a cyclic shift of $\mu^1$, the first edges on the paths
      $P_{\mu^k_1}, P_{\mu^k_2}, \dots, P_{\mu^k_\ell}$ form an anti-clockwise order
      around $p$.

    In order to comply with Definition \ref{order-def},
      we still need to prove that for each vertex $v \in V(\Qc)$, there exist two indices $L, R$
      ($1 \leq L \leq R \leq \ell$) such that $v$ lies exactly on paths
      $Q_{\mu^k_L}, Q_{\mu^k_{L+1}}, \dots, Q_{\mu^k_R}$.

    Fix a vertex $v \in \Qc$.
    Note that $\Kc^k_v$ is a cyclic interval (Claim \ref{order-lemma-cyclic-interval}).
    If $\Kc^k_v = [1, \ell]$, then the proposition we are proving is trivial.
    On the other hand, assume that $\Kc_v \neq [1, \ell]$.
    From our assumptions, we have that $|\Kc^1_v| \leq |\Kc^1_{v^*}|$, so
      also $|\Kc^k_v| \leq |\Kc^k_{v^*}|$.
    From Claim \ref{order-lemma-laminar-family}, we get that $\Kc^k_v$ is either disjoint
      or fully contained within $\Kc^k_{v^*}$; it cannot contain $\Kc^k_{v^*}$ strictly due to
      $|\Kc^k_v| \leq |\Kc^k_{v^*}|$.

    Recall that $\Kc^k_{v^*} = [1, M]$.
    If $\Kc^k_v$ is disjoint with $\Kc^k_{v^*}$, then $\Kc^k_v \subseteq [M+1, \ell]$,
      and $\Kc^k_v$ is a non-cyclic interval.
    If $\Kc^k_v$ is fully contained within $\Kc^k_{v^*}$, then $\Kc^k_v \subseteq [1, M]$,
      and again $\Kc^k_v$ is a non-cyclic interval.

    Hence, every $v \in V(\mathcal{Q})$ can be assigned a subinterval $
      \Kc^k_v = [L, R] \subseteq [1, \ell]$
      such that $v$ lies exactly on paths $Q_{\mu^k_L}, Q_{\mu^k_{L+1}}, \dots, Q_{\mu^k_R}$.
    This finishes the proof that $\mu^k$ is an order.
    \end{claimproof}
  \end{claim}
  Claim \ref{order-lemma-found-claim} finishes the proof of the lemma.
  \end{proof}
\end{lemma}

While it is possible to proceed to the next parts of the proof
  using ordered cages, the fact that we need to use
  two different orderings --- one imposed by the semi-ladder, and one imposed by
  the order $\sigma$ in the ordered cage --- will be excessively unwieldy.
In order to remedy this inconvenience, we will introduce a~variant of ordered cages unifying
  both orderings.
This variant will lead us to a slightly weaker upper bound on the maximum semi-ladder size,
  though the bound will still remain asymptotically the same.
  
\begin{definition}
  An \textbf{identity ordered cage} of order $\ell$ is an ordered cage
  $$\textsf{OrderedCage}(\allowbreak(a_1, \dots, a_\ell), (b_1, \dots, b_\ell), \Pc, \Qc, \sigma)\ \ \text{where }\sigma = (1, 2, \dots, \ell)\text{ is the identity permutation.}$$
  
  \noindent We denote this structure as
    $$ \textsf{IdOrderedCage}((a_1, \dots, a_\ell), (b_1, \dots, b_\ell), \Pc, \Qc). $$
\end{definition}

We immediately follow with a proof that large identity ordered cages exist in planar graphs
  given that large ordered cages exist.
We, however, remark that the proof might alter the embedding of a graph in the plane.

\begin{lemma}
  \label{identity-ordered-cage-exists}
  If a graph $G$ embedded in the plane contains an ordered cage of order $(\ell - 1)^2 + 1$
    for $\ell \geq 1$,
    then there exists an embedding of $G$ in the plane in which $G$ contains an
    identity ordered cage of order $\ell$.
  \begin{proof}
    Take any graph $G$ with an ordered cage $\Cc$ of order $(\ell-1)^2 + 1$:
      $$\Cc = \OrderedCage((a_i)_{i=1}^{(\ell-1)^2+1},
      (b_i)_{i=1}^{(\ell-1)^2+1}, \Pc, \Qc, \sigma).$$
    By Erd\H{o}s-Szekeres Theorem \cite{10.1007/978-1-4612-0801-3_9},
      $\sigma$ contains either an increasing subsequence
      of length $\ell$, or a decreasing subsequence of length $\ell$.
      
    If $\sigma$ contains an increasing subsequence $\sigma_{\mathrm{inc}}$ of length $\ell$,
      then a subset of $\Cc$ with indices $\sigma_{\mathrm{inc}}$
      is an identity ordered cage of order $\ell$ in $G$ with the original embedding in the plane.
      
    Assume now that $\sigma$ contains a decreasing subsequence $\sigma_{\mathrm{dec}}$
      of length $\ell$.
    In this case, we modify the embedding of $G$ in the plane by reflecting it through
      any line in the Euclidean plane.
    Since in the original embedding, the first edges on paths $P_{\sigma_1}, P_{\sigma_2}, \dots,
      P_{\sigma_{(\ell-1)^2+1}}$ were ordered anti-clockwise around $\Root(\Pc)$,
      these edges are ordered clockwise in the reflected embedding.
    Hence, the reversed permutation $\sigma^{\rev} =
      (\sigma_{(\ell-1)^2 + 1}, \sigma_{(\ell-1)^2}, \dots, \sigma_1)$ is an order in the
      modified embedding, since
      the first edges on paths $P_{\sigma_{(\ell-1)^2+1}}, P_{\sigma_{(\ell-1)^2}}, \dots, P_{\sigma_1}$
      are ordered anti-clockwise in the altered embedding.
    The remaining condition in Definition \ref{ordered-cage-def} can be verified in a
        straightforward way.
    Now, we notice that $\sigma_{\mathrm{dec}}^{\rev}$ (the reversed subsequence
      $\sigma_{\mathrm{dec}}$) is an increasing subsequence in the newly formed order
      $\sigma^{\rev}$.
    Hence, the subset of $\Cc$ indicated by $\sigma_{\mathrm{dec}}^{\rev}$
      in the transformed embedding is an identity ordered cage.
  \end{proof}
\end{lemma}

An example identity ordered cage is depicted in Figure \ref{figure-neighbor-areas}(a).

In identity ordered cages, the following useful fact holds:

\begin{lemma}
\label{lca-composition-lemma}
  Fix an identity ordered cage $\IdOrderedCage((a_i)_{i=1}^\ell, (b_i)_{i=1}^\ell, \Pc, \Qc)$
    in a graph $G$ and three indices $i, j, k$ ($1 \leq i < j < k \leq \ell$).
  In the set $\{\lca_\Qc(a_i, a_j), \lca_\Qc(a_j, a_k)\}$,
    one of the elements is equal to $\lca_\Qc(a_i, a_k)$ and is an ancestor of
    the other element in the set (it can be, in particular, equal to it).
  \begin{proof}
    Since both elements of the set $\{\lca_\Qc(a_i, a_j), \lca_\Qc(a_j, a_k)\}$
      are the ancestors of $a_j$ in $\Qc$, one of these elements is an ancestor of the other.
    Without loss of generality, assume that $\lca_\Qc(a_i, a_j)$ is an ancestor of the other
      element of the set.
    This means that $a_k$ is in the subtree rooted at $\lca_\Qc(a_i, a_j)$.
    As both $a_i$ and $a_k$ are in this subtree, we infer that $\lca_\Qc(a_i, a_j)$ is an ancestor
      of $\lca_\Qc(a_i, a_k)$.

    Since $\lca_\Qc(a_i, a_k) \in Q_i$, $\lca_\Qc(a_i, a_k) \in Q_k$, and $j \in [i, k]$,
      by the definition of an identity ordered cage we also have $\lca_\Qc(a_i, a_k) \in Q_j$.
    Therefore, $\lca_\Qc(a_i, a_k)$ is also an ancestor of $a_j$ in $\Qc$.
    Hence, $\lca_\Qc(a_i, a_k)$ is an ancestor of $\lca_\Qc(a_i, a_j)$.
    
    We conclude that $\lca_\Qc(a_i, a_j) = \lca_\Qc(a_i, a_k)$.
  \end{proof}
\end{lemma}
  
\subsection{Neighbor areas}
\label{neighbor-area-section}

Until now, we were concerned about the structure imposed only by
  some vertices of the semi-ladder --- that is, the vertices $a_1, a_2, \dots, a_\ell$.
It is now time to include the corresponding vertices $b_1, b_2, \dots, b_\ell$ in the picture.
It turns out that each sufficiently large ordered cage contains a smaller ordered cage as a subset
  in which for every $i$, the vertex $b_i$ is, in some topological sense, \emph{near} $a_i$.
This is where the upper bounds on the neighborhood complexity in planar graphs proved in Section
  \ref{noose-profile-lemma-section}, will come into play.
  
First, we need to design technical tools allowing us to formalize this idea.

\begin{definition}
Let $\Cc = \IdOrderedCage((a_i)_{i=1}^\ell, (b_i)_{i=1}^\ell, \Pc, \Qc)$ be an identity ordered
  cage of order $\ell$ in $G$.
For every two integers $i, j$ ($ 1\leq i < j \leq \ell$), we set $v_{i,j} = \lca_\Qc(a_i, a_j)$.
Note that by Lemma \ref{splitting-path-intersection}, $v_{i,j}$ is the first vertex of the common
  suffix of $\SPath(\Cc, i)$ and $\SPath(\Cc, j)$.

For these integers $i, j$,
  we define the \textbf{neighbor area} $\Area(\Cc, i, j)$ as the closed
  part of the plane whose boundary is traversed anti-clockwise by the following directed cycle:
  $$ \SPath(\Cc, i)[\Root(\Pc), v_{i, j}] \cdot
    (\SPath(\Cc, j))^{-1}[v_{i, j}, \Root(\Pc)]. $$
\end{definition}

\begin{figure}[h]
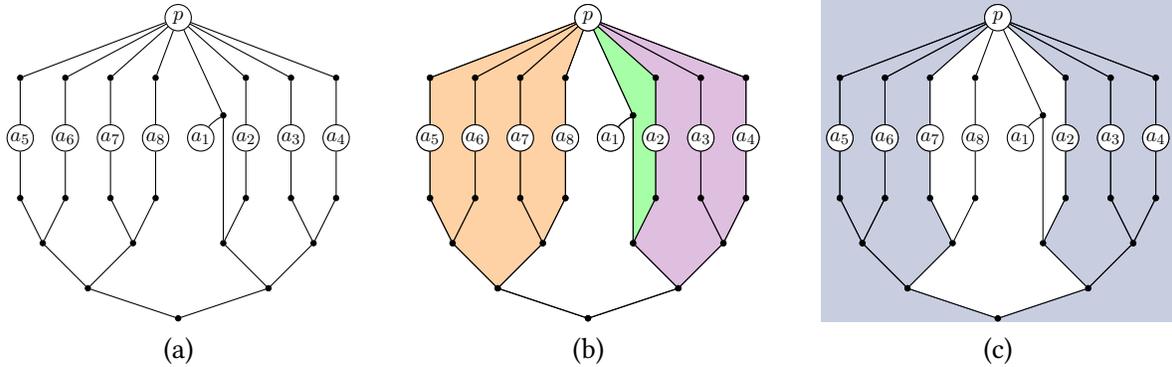

  \centering
  \begin{minipage}[b]{0.33\textwidth}
  \centering
  \input{figures/neighbor-area-a.tex}
  (a)
  \end{minipage}\begin{minipage}[b]{0.33\textwidth}
  \centering
  \input{figures/neighbor-area-b.tex}
  (b)
  \end{minipage}\begin{minipage}[b]{0.33\textwidth}
  \centering
  \input{figures/neighbor-area-c.tex}
  (c)
  \end{minipage}
  \caption{(a) --- an example identity ordered cage $\Cc$ of order $8$ embedded in the plane. \\
  (b) --- neighbor areas $\Area(\Cc, 1, 2)$ (green), $\Area(\Cc, 2, 4)$ (violet), and
    $\Area(\Cc, 5, 8)$ (orange). \\
  (c) --- an unbounded neighbor area $\Area(\Cc, 2, 7)$ (blue). The area is unbounded
    since its boundary is traversed anti-clockwise first by $\SPath(\Cc, 2)$, and later
    by $(\SPath(\Cc, 7))^{-1}$.}
  \label{figure-neighbor-areas}
\end{figure}

Note that neighbor areas are not necessarily bounded (Figure \ref{figure-neighbor-areas}).

We can easily see from Figure \ref{figure-neighbor-areas} that the union
  of two neighbor areas $\Area(\Cc, 1, 2)$ and $\Area(\Cc, 2, 4)$ in the presented example is
  $\Area(\Cc, 1, 4)$, and that the intersection of these is a prefix of $\SPath(\Cc, 2)$.
It turns out that this observation is satisfied in general:

\begin{lemma}
\label{area-sum-lemma}
In an identity ordered cage $\Cc$ of order $\ell$, for every indices $i, j, k$ such that
  $1 \leq i < j < k \leq \ell$, we have that:
  \begin{equation*}
  \begin{split}
  \Area(\Cc, i, j) \cup \Area(\Cc, j, k) &= \Area(\Cc, i, k), \\
  \Area(\Cc, i, j) \cap \Area(\Cc, j, k) &= \partial\Area(\Cc,i,j) \cap \partial\Area(\Cc,j,k) \subseteq
    \SPath(\Cc, j).
  \end{split}
  \end{equation*}
  \begin{proof}
  We first define vertices $v_{i, j} = \lca_\Qc(a_i, a_j)$, $v_{j, k}=\lca_\Qc(a_j, a_k)$
    as in the definition of neighbor areas.
  We also let $\Cc = \IdOrderedCage((a_i)_{i=1}^\ell, (b_i)_{i=1}^\ell, \Pc, \Qc)$,
    $p := \Root(\Pc)$, and $\Qc = \mathsf{Tree}(Q_1, \dots, Q_\ell)$.
  These definitions tell us that
  \begin{equation*}
  \begin{split}
  \partial\Area(\Cc, i, j) &= \SPath(\Cc, i)[p, v_{i, j}] \cdot
    (\SPath(\Cc, j))^{-1}[v_{i, j}, p], \\
  \partial\Area(\Cc, j, k) &= \SPath(\Cc, j)[p, v_{j, k}] \cdot
    (\SPath(\Cc, k))^{-1}[v_{j, k}, p].
  \end{split}
  \end{equation*}

  Since $v_{i, j}$ and $v_{j, k}$ are both ancestors of $a_j$ in $\Qc$, either of these is
    an ancestor of the other (in particular, we might have $v_{i, j} = v_{j, k}$).
  Without loss of generality, we assume that $v_{j, k}$ is an ancestor of $v_{i, j}$ in $\Qc$
    or equal to $v_{i, j}$ (Figure \ref{area-sum-layout-fig}).
  Then, $\SPath(\Cc, j)[p, v_{i, j}] \subseteq \SPath(\Cc, j)[p, v_{j, k}]$.
  Moreover, by Lemma \ref{lca-composition-lemma}, $v_{i,k} := \lca(a_i, a_k)$ is equal to $v_{j,k}$.

  \begin{figure}[h]
  \centering
    \input{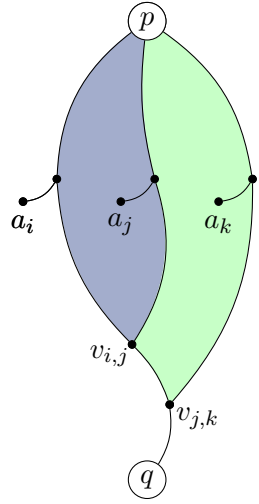}
    \caption{The configuration assumed in Lemma \ref{area-sum-lemma}. The vertices
      $v_{i,j}$ and $v_{j,k}$ may overlap. The blue area is $\Area(\Cc, i, j)$. The green area is
      $\Area(\Cc, j, k)$. The brown path is $\SPath(\Cc, j)$.}
    \label{area-sum-layout-fig}
  \end{figure}
    
  \vspace{0.5em}
  
  We first note that $\partial\Area(\Cc, i, j) \cap \partial\Area(\Cc, j, k) \subseteq \SPath(\Cc, j)$;
    this is because any intersection of both boundaries outside of $\SPath(\Cc, j)$ would
    have to be an intersection of $\SPath(\Cc, i)$ with $\SPath(\Cc, k)$.
  These paths intersect only at $p$ and $\SPath(\Cc, i)[v_{i, k}, q]$ (Lemma
    \ref{splitting-path-intersection}).
  However, $v_{i, k} = v_{j, k} \in \SPath(\Cc, j)$, so also $\SPath(\Cc, i)[v_{i, k}, q] \subseteq
    \SPath(\Cc, j)$.
  Obviously, $p \in \SPath(\Cc, j)$ as well.
  Therefore, no intersection of $\partial\Area(\Cc, i, j)$ and $\partial\Area(\Cc, j, k)$ can occur
    outside of $\SPath(\Cc, j)$.
    
  Thanks to this fact, we can now compute $\partial\Area(\Cc, i, j) \cap \partial\Area(\Cc, j, k)$:
  \begin{equation*}
  \begin{split}
    \partial\Area(\Cc, i, j) & \cap \partial\Area(\Cc, j, k) = \\
    & = (\partial\Area(\Cc, i, j) \cap \SPath(\Cc, j)) \cap (\partial\Area(\Cc, j, k) \cap \SPath(\Cc, j))
      = \\
    & = \SPath(\Cc, j)[p, v_{i, j}] \cap \SPath(\Cc, j)[p, v_{j, k}] = \SPath(\Cc, j)[p, v_{i, j}]
  \end{split}
  \end{equation*}
  --- that is, this intersection is a prefix of $\SPath(\Cc, j)$.
  
  Therefore, the two boundaries intersect at exactly one segment.
  We can use this information to prove that the interiors of $\Area(\Cc, i, j)$ and $\Area(\Cc, j, k)$
    are disjoint.
  Let us name $f_i, f_j, f_k$ as the second vertices on the paths $P_i, P_j, P_k$, respectively;
    that is, the vertices immediately following $p = \Root(\Pc)$ on the respective paths.
  Since the counter-clockwise orientation of $\partial\Area(\Cc, i, j)$ contains
    oriented edges $f_jp$ and $pf_i$ in this order,
    and the counter-clockwise orientation of $\partial\Area(\Cc, j, k)$ contains
    oriented edges $f_kp$ and $pf_j$ in this order,
    we only need to verify the cyclic order of the edges $pf_i, pf_j, pf_k$ around vertex $p$.
  These edges are ordered anti-clockwise around $p$ (by the definition of an ordered
    cage), so the interiors of the areas $\Area(\Cc, i, j)$, $\Area(\Cc, j, k)$
    are disjoint (Figure \ref{area-sum-layout-fig}).
    
  From the considerations above, we conclude that
  $$ \Area(\Cc, i, j) \cap \Area(\Cc, j, k) = \partial\Area(\Cc, i, j) \cap \partial \Area(\Cc, j, k)
    \subseteq \SPath(\Cc, j). $$
    
  Moreover, since the interiors of both neighbor areas are disjoint,
    given the anti-clockwise orientations of $\partial\Area(\Cc, i, j)$ and
    $\partial\Area(\Cc, j, k)$, both starting and finishing at their common vertex $p$,
    we can compute the anti-clockwise orientation of $\partial[\Area(\Cc, i, j) \cup
    \Area(\Cc, j, k)]$ by composing the anti-clockwise orientations of
    $\partial\Area(\Cc, i, j)$ and $\partial\Area(\Cc, j, k)$:
  \begin{equation*}
  \begin{split}
    & \partial[\Area(\Cc, i, j) \cup \Area(\Cc, j, k)] = \\
    & = \SPath(\Cc, i)[p, v_{i, j}] \cdot (\SPath(\Cc, j))^{-1}[v_{i, j}, p] \cdot
        \SPath(\Cc, j)[p, v_{j, k}] \cdot (\SPath(\Cc, k))^{-1}[v_{j, k}, p] = \\
    & = \SPath(\Cc, i)[p, v_{i, j}] \cdot \SPath(\Cc, j)[v_{i, j}, v_{j, k}] \cdot
      (\SPath(\Cc, k))^{-1}[v_{j, k}, p] \stackrel{(\star)}{=} \\
    & \stackrel{(\star)}{=} \SPath(\Cc, i)[p, v_{i, j}] \cdot \SPath(\Cc, i)[v_{i, j}, v_{j, k}] \cdot
      (\SPath(\Cc, k))^{-1}[v_{j, k}, p] = \\
    & = \SPath(\Cc, i)[p, v_{j, k}] \cdot (\SPath(\Cc, k))^{-1}[v_{j, k}, p] = \\
    & =
      \SPath(\Cc, i)[p, v_{i, k}] \cdot (\SPath(\Cc, k))^{-1}[v_{i, k}, p] =
      \partial\Area(\Cc, i, k).
  \end{split}
  \end{equation*}
  
  In the equality $(\star)$, we use the definition of $v_{i, j}$ --- it is the first vertex
    of the common suffix of $\SPath(\Cc, i)$ and $\SPath(\Cc, j)$.
  
  Therefore, $\Area(\Cc, i, j) \cup \Area(\Cc, j, k) = \Area(\Cc, i, k)$.
  \end{proof}
\end{lemma}

We can immediately infer the following fact from Lemma \ref{area-sum-lemma}:

\begin{corollary}
\label{area-int-lemma}
In any identity ordered cage $\Cc$ of order $\ell$, for every triple of indices $i, j, k$ such that
  $1 \leq i < j < k \leq \ell$, we have that:
$$ \Area(\Cc, i, j) \cap \mathrm{Int}\,\Area(\Cc, j, k) = \varnothing,
  \qquad \mathrm{Int}\,\Area(\Cc, i, j) \cap \Area(\Cc, j, k) = \varnothing. $$
\begin{proof}
  Since areas $\Area(\Cc, i, j)$ and $\Area(\Cc, j, k)$ intersect only at their boundaries,
    the interior of either neighbor area is disjoint with the other neighbor area.
\end{proof}
\end{corollary}

Now we will examine the importance of neighbor areas in identity ordered cages.
For an identity ordered cage $\Cc$ of order $\ell$ and two indices $i, j$ ($i < j$),
  we can observe which vertices of the semi-ladder underlying $\Cc$:
  $a_1, a_2, \dots, a_\ell, b_1, b_2, \dots, b_\ell$ belong to $\Area(\Cc, i, j)$,
  and which do not.
It turns out that there cannot be too many pairs of vertices $(a_k, b_k)$ for $k \in [1, \ell]$
  such that one vertex in the pair belongs to $\Area(\Cc, i, j)$ and the other does not.
For simplicity of the further part of the proof, we use the following definition:

\begin{definition}
\label{separated-from-def}
Fix a neighbor area $\Area(\Cc, i, j)$ in an identity ordered cage $\Cc$.
For a~pair of vertices $u, v$, if exactly one of them belongs to $\Area(\Cc, i, j)$,
  we will say that $u$ \textbf{is separated from} $v$ by $\Area(\Cc, i, j)$.
\end{definition}

For instance, in Figure \ref{figure-neighbor-areas}, $a_5$ is separated from $a_3$ by
  $\Area(\Cc, 2, 4)$, but $a_5$ is not separated from $a_7$ by the same neighbor area.

It turns out that the definition of ``being separated'' is connected in a non-trivial way
  to the distinctness of distance-$d$ profiles.
This will allow us to employ the Noose Profile Lemma (Theorem \ref{noose-profile-lemma}),
  which bounds the neighborhood complexity of nooses in planar graphs,
  in order to limit the number of pairs
  of the vertices of the semi-ladder separated by neighbor areas.

\begin{lemma}
\label{neighbor-separation-lemma}
In an identity ordered cage $\Cc = \IdOrderedCage((a_i)_{i=1}^\ell,
  (b_i)_{i=1}^\ell, \Pc, \Qc)$ of order~$\ell$, for every pair of indices $i, j$ ($1 \leq i < j \leq
  \ell$), the neighbor area $\Area(\Cc, i, j)$ separates at most $128d^3(d+2)^4$ distinct pairs
  of vertices $(a_k, b_k)$ for $k \in [1, \ell]$.
  \begin{proof}
  Fix two integers $i, j$ ($1 \leq i < j \leq \ell$) for the proof.
  Denote the set of vertices on the boundary of $\Area(\Cc, i, j)$ as $B$.
  \begin{claim}
    \label{large-neighbor-small-noose}
    $|B| \leq 4d$.
    \begin{claimproof}
    \cqed
      We first note that $|E(\SPath(\Cc, k))| \leq 2d$ for each $k \in [1, L]$ since
        $\SPath(\Cc, k)$ is the concatenation of subpaths of two paths,
        $\Path(\Pc, k)$ and $\Path(\Qc, k)$, whose lengths do not exceed $d$.
      The boundary of $\Area(\Cc, i, j)$ is a cycle and
        consists of two prefixes of paths $\SPath(\Cc, i)$ and $\SPath(\Cc, j)$.
      Therefore,
      \[ |B| = |V(\partial \Area(\Cc, i, j))| = |E(\partial \Area(\Cc, i, j))| \leq
        |E(\SPath(\Cc, i))| + |E(\SPath(\Cc, j))| \leq 4d. \]
    \end{claimproof}
  \end{claim}
  
  Let $I$ denote the set of indices $k$ for which $a_k$ is separated from $b_k$ by
    $\Area(\Cc, i, j)$; we need to prove that $|I| \leq 128d^3(d + 2)^4$.
    
  For each $k \in I$, two cases are possible --- either $a_k$ belongs to $\Area(\Cc, i, j)$,
    which implies that $b_k \not\in \Area(\Cc, i, j)$; or $a_k \not\in \Area(\Cc, i, j)$, which
    implies that $b_k \in \Area(\Cc, i, j)$.
  We will partition $I$ into two subsets --- $I_1$ and $I_2$ --- depending on which
    case holds for a given $k$:

  $$ I_1 = \{k \in I\, \mid\, a_k \in \Area(\Cc, i, j),\ b_k \not\in \Area(\Cc, i, j)\}, $$
  $$ I_2 = \{k \in I\, \mid\, a_k \not\in \Area(\Cc, i, j),\ b_k \in \Area(\Cc, i, j)\}. $$
  
  We will prove that both $I_1$ and $I_2$ contain no more than $64d^3(d + 2)^4$ indices;
    the statement of the lemma will follow immediately.
  
  \begin{claim}
    \label{large-neighbor-distinct-profiles}
    For any $t \in \{1, 2\}$ and two different indices $x, y \in I_t$,
      vertices $a_x$ and $a_y$ have different distance-$d$ profiles on $B$.
    \begin{claimproof}
    \cqed
    Assume on the contrary that $a_x$ and $a_y$ have identical distance-$d$ profiles on $B$.
    Without loss of generality, assume that $x < y$.
    Because vertices $a_x, a_y, b_x, b_y$ form a distance-$d$ semi-ladder of order $2$,
      it follows that
    $$
      \dist(a_y, b_x) \leq d, \quad\text{ but }\quad \dist(a_x, b_x) > d.
    $$
    \begin{figure}[h]
      \centering
      \input{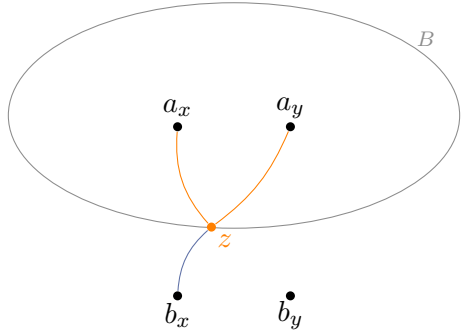}
      \caption{The configuration in Claim \ref{large-neighbor-distinct-profiles}.}
      \label{separating-profiles-fig}
    \end{figure}
    Since $x, y \in I_k$, we have that $a_y$ is separated from $b_x$ by $\Area(\Cc, i, j)$.
    Therefore, the shortest path
      between these two vertices must pass through a vertex of $\partial \Area(\Cc, i, j)$ ---
      that is, some vertex $z \in B$ (Figure \ref{separating-profiles-fig}).
    Thus,
    $$
      \dist(a_y, b_x) = \dist(a_y, z) + \dist(z, b_x).
    $$
    It follows that $\dist(a_y, z) \leq d$ and hence $\pi_d[a_y, B](z) = \dist(a_y, z)$.
    Because distance-$d$ profiles of $a_y$ and $a_x$ on $B$ are equal, this implies that
    $\pi_d[a_x, B](z) = \pi_d[a_y, B](z) = \dist(a_y, z)$.
    Since also $\pi_d[a_x, B](z) \leq d$, we get that $\dist(a_x, z) = \pi_d[a_x, B](z) = \dist(a_y, z)$,
      so the following inequality follows from the triangle inequality:
    $$
      \dist(a_x, b_x) \leq \dist(a_x, z) + \dist(z, b_x) = \dist(a_y, z) + \dist(z, b_x) = \dist(a_y, b_x)
      \leq d
    $$
    --- a contradiction.
    Therefore, $a_x$ and $a_y$ must have distinct distance-$d$ profiles on $B$.
    \end{claimproof}
  \end{claim}
  
  A combination of both claims above provides us with a link between the notion of cages and
    the Noose Profile Lemma.
  Assume for contradiction that for an integer $t \in \{1, 2\}$, the set $I_t$ contains more than
    $64 d^3 (d+2)^4$ elements.
  By Claim \ref{large-neighbor-distinct-profiles}, we infer that
    $\left|\{\pi_d[a_x, B]\, \mid\, x \in I_k\}\right| > 64 d^3 (d + 2)^4$.
  We also remark that since $B$ is a cycle in $G$, there exists a noose $\Lc$ in the embedding
    of $G$ in the plane passing through vertices in $B$ and no other vertices
    ($\Lc$ should closely follow the cycle $B$, passing through each vertex of $B$ in the
    same order as $B$, but not touching the edges of the graph).
  $\Lc$ is now subject to the Noose Profile Lemma.
  We can see that either all vertices in $\{a_k\,\mid\, k \in I_t\}$ are enclosed by $\Lc$, or
    none of them are enclosed by $\Lc$ (depending on the value of $t$ and whether
    $\Area(\Cc, i, j)$ is bounded or unbounded).
  If these vertices are all enclosed by $\Lc$, then Noose Profile Lemma (Theorem
    \ref{noose-profile-lemma}) applies, so the number of different distance-$d$ profiles
    seen from $\Lc$ is bounded from above by
  $$ |B|^3 \cdot (d+2)^4 \leq (4d)^3 \cdot (d + 2)^4 = 64d^3 (d+2)^4 $$
  --- a contradiction.
  Analogously, if none of vertices in $\{a_k\,\mid\, k \in I_t\}$ are enclosed by $\Lc$,
    then a variant of the Noose Profile Lemma (Corollary \ref{noose-profille-lemma-rev}) applies,
    and we contradict our assumption in the same fashion.
  Therefore, $|I_t| \leq 64d^3 (d+2)^4$.
    
  Hence, $|I| = |I_1| + |I_2| \leq 128d^3 (d+2)^4$.
  \end{proof}
\end{lemma}

We will also prove that in an identity ordered cage $\Cc$, for each $i \in [2, \ell - 1]$,
  we can expect to find the vertex $a_i$ in the interior of the neighbor area $\Area(\Cc,i-1,i+1)$
  --- roughly speaking, ``near'' the path $\SPath(\Cc, i)$.

\begin{lemma}
  \label{neighbor-cage-ai-good-lemma}
  In an identity ordered cage $\Cc = 
    \IdOrderedCage((a_1, \dots, a_\ell), (b_1, \dots, b_\ell), \Pc, \Qc)$ in graph $G$,
    for every $i \in [2, \ell - 1]$,
    we have that $a_i \in \mathrm{Int}\, \Area(\Cc, i-1, i+1)$.
  \begin{proof}
    For the sake of contradiction assume that $a_i \not\in \mathrm{Int}\, \Area(\Cc, i-1, i+1)$.  
  
    For all $j \in [1, \ell]$, we choose $f_j$ as the second vertex on $\Path(\Pc, j)$ --- that is,
      the vertex on this path immediately following $\Root(\Pc)$.
    We also let $p := \Root(\Pc)$,
    $\Pc = \mathsf{Tree}(P_1, \dots, P_\ell)$,
    and $\Qc = \mathsf{Tree}(Q_1, \dots, Q_\ell)$.

    Since edges $pf_{i-1}, pf_i, pf_{i+1}$ are ordered anti-clockwise around $p$ (by the
      definition of an ordered cage), and because the boundary of $\Area(\Cc, i-1, i+1)$,
      when ordered anti-clockwise, contains edges $f_{i+1}p$ and $pf_{i-1}$ in this order,
      we have that $pf_i \subseteq \Area(\Cc, i-1, i+1)$ (Figure~\ref{ai-acw-ordering-fig}).
    Thus, $f_i \in \Area(\Cc, i-1, i+1)$.
    Since $a_i \not\in \mathrm{Int}\, \Area(\Cc, i-1, i+1)$, the path $P_i[f_i, a_i]$ must intersect
      $\partial\Area(\Cc, i-1, i+1)$.
      
    \begin{figure}[h]
      \centering
      \input{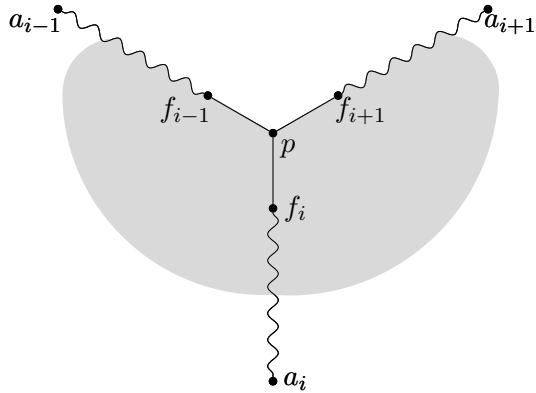}
      \caption{A simple geodesic tree $\Pc$ (a part of the identity ordered cage $\Cc$)
        with important vertices in Lemma \ref{neighbor-cage-ai-good-lemma}.
        $\Area(\Cc, i-1, i+1)$ is shaded.}
      \label{ai-acw-ordering-fig}
    \end{figure}
    
    By the definition of neighbor areas, $\partial\Area(\Cc, i-1, i+1) \subseteq
      P_{i-1} \cup P_{i+1} \cup Q_{i-1} \cup Q_{i+1}$.
    However, $P_i[f_i, a_i]$ is vertex-disjoint with $P_j$ for each $j \neq i$ (since
      by the definition of a quasi-cage, paths $P_i$ and $P_j$ are vertex-disjoint apart
      from $p$).
    Similarly, $P_i$ is vertex-disjoint with $Q_j$ for each $j \neq i$ (by the
      definition of a cage).
    Hence, $P_i[f_i, a_i]$ is vertex-disjoint with $\partial\Area(\Cc, i-1, i+1)$ --- a contradiction.
  \end{proof}
\end{lemma}

\subsection{Neighbor cages}
\label{neighbor-cages-section}
We are now ready to use the tools developed in Section \ref{neighbor-area-section}
  to define neighbor cages and explore their properties.

\begin{definition}
\label{neighbor-cage-def}
An identity ordered cage $\Cc = \IdOrderedCage((a_i)_{i=1}^\ell, (b_i)_{i=1}^\ell, \Pc, \Qc)$
  of order $\ell$ is also a \textbf{neighbor cage}
  if for each index $i \in \{2, 3, \dots, \ell - 1\}$, the vertex $b_i$ belongs to
  $\Area(\Cc, i - 1, i + 1)$.
  
We denote neighbor cages in the following way:
  $$ \NeighborCage((a_i)_{i=1}^\ell, (b_i)_{i=1}^\ell, \Pc, \Qc). $$
\end{definition}

\begin{figure}[h]
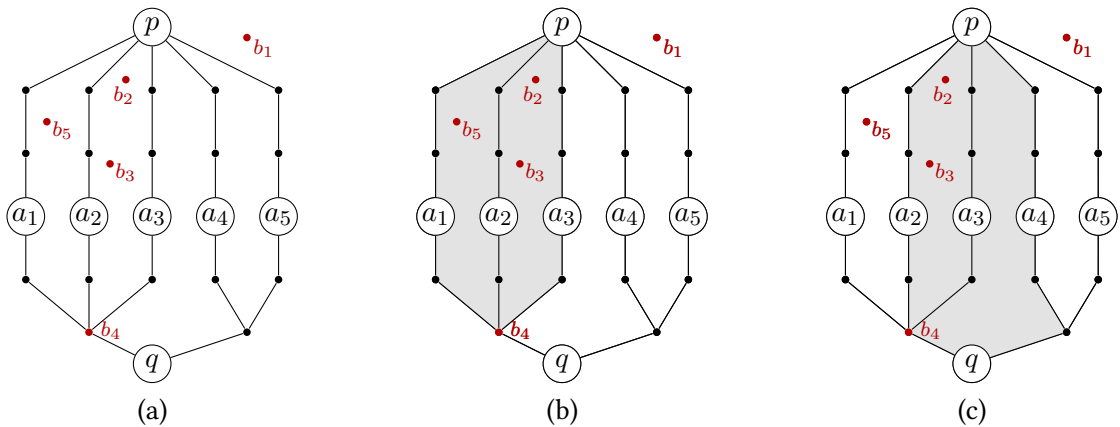

  \centering
  \begin{minipage}[b]{0.33\textwidth}
      \centering
      \input{figures/neighbor-a.tex}
      (a)
    \end{minipage}\begin{minipage}[b]{0.33\textwidth}
      \centering
      \input{figures/neighbor-b.tex}
      (b)
    \end{minipage}\begin{minipage}[b]{0.33\textwidth}
      \centering
      \input{figures/neighbor-c.tex}
      (c)
    \end{minipage}
  \caption{(a) --- An example neighbor cage of order $5$, with edges connecting
    $b_1, \dots, b_5$ with the remaining part of the graph removed. Note that $b_1$ and $b_5$ may
    be anywhere in the graph. \\
    (b) --- The neighbor area which must contain $b_2$ is shaded. \\
    (c) --- The neighbor area which must contain $b_3$ is shaded.}
  \label{neighbor-examples-fig}
\end{figure}

Figure \ref{neighbor-examples-fig} shows an example of a neighbor cage.
Unnecessary parts of the graph (including the edges connecting the vertices $b_i$
  with the remaining part of the graph) have been omitted.
  
We will now employ the toolchain we have developed to prove that large identity ordered cages
  contain large neighbor cages as subsets.

\begin{lemma}
\label{large-neighbor-lemma}
Each identity ordered cage of order $[256d^3(d+2)^4 + 2] \ell + 2$ contains a neighbor cage
  of order $\ell$ as a subset.
  
  \begin{proof}
  Let $L := [256d^3(d+2)^4 + 1] \ell + 2$,
  and fix an identity ordered cage of order $L$:
  $$\Cc = \IdOrderedCage((a_i)_{i=1}^L, (b_i)_{i=1}^L, \Pc, \Qc).$$
  We build a permutation $\pi$ of vertices of the semi-ladder $\{a_1, a_2, \dots, a_L,
    b_1, b_2, \dots, b_L\}$ in the following way:
    
  \begin{itemize}
    \item $\pi$ contains vertices $a_1, a_2, \dots, a_L$ in this order; immediately after each $a_i$
      (for $i \in [1, L]$), we add a (possibly empty) sequence $\Gamma_i$ which is a
      permuted subset of $\{b_1, b_2, \dots, b_L\}$. That is,
      $$ \pi = (a_1) \cdot \Gamma_1 \cdot (a_2) \cdot \Gamma_2 \cdot \dots \cdot (a_L) \cdot
         \Gamma_L, $$
       where $\cdot$ denotes the concatenation of sequences.
    \item For each index $i \in [1, L - 1]$, $\Gamma_i$ contains those vertices $b_1, \dots, b_L$
       which belong to $\Area(\Cc, i, i+1)$,
       but do not belong to any previous sequence $\Gamma_1, \dots, \Gamma_{i-1}$.
       Formally, the set of vertices in $\Gamma_i$ is
       $$ \left[\{b_1, b_2, \dots, b_L\} \cap \Area(\Cc, i, i+1)\right] \setminus (\Gamma_1 \cup
         \dots \cup \Gamma_{i-1}). $$
    \item $\Gamma_L$ contains all vertices in $\{b_1, b_2, \dots, b_L\}$ which do not belong
      to any previous sequence $\Gamma_1, \dots, \Gamma_{L-1}$.
    \item Each sequence $\Gamma_i$ (for $i \in [1, L]$) is ordered arbitrarily.
  \end{itemize}
  
  It can be clearly seen that the procedure above produces a permutation $\pi$.
  
  \vspace{1em}

  We will now find a large subsequence of $\pi$ with the following property:
    if $a_i$ belongs to the subsequence for some $i \in [1,L]$, then $b_i$ is a neighbor of $a_i$
    in this subsequence; and conversely, if $b_i$ belongs to this subsequence, then
    $a_i$ neighbors $b_i$ in this subsequence.

  For all integers $i \in [1, L]$, we define the values $\alpha(i)$ and $\beta(i)$
    so that $a_i = \pi_{\alpha(i)}$ and $b_i = \pi_{\beta(i)}$.
  Note that by the definition of $\pi$, we have $\alpha(1) < \alpha(2) < \dots < \alpha(L)$.

  This definition allows us to define the \textbf{conflict interval} $X_i$ for each $i \in [2, L-1]$ ---
    a closed interval whose endpoints are at $\alpha(i)$ and $\beta(i)$.
  Each conflict interval is a subinterval of the interval $[1, 2L]$.
  We stress that we do not define $X_1$ or $X_L$.
  
  We note the following property of the conflict intervals:
  \begin{claim}
    \label{large-neighbor-separation-claim}
    For every two different integers $s, t$ ($s \in [2, L-1]$, $t \in [2, L]$),
      $a_s$ is separated from $b_s$ by $\Area(\Cc, 1, t)$
      if and only if $\alpha(t) \in X_s$.
    \begin{claimproof}
    \cqed
    Let us first remark that
    \begin{equation}
    \label{prefix-area-sum-eq}
    \Area(\Cc, 1, t) = \Area(\Cc, 1, 2) \cup \Area(\Cc, 2, 3) \cup \dots \cup \Area(\Cc, t-1, t)
    \end{equation}
    (Lemma \ref{area-sum-lemma}).
    Therefore, if $\beta(s) < \alpha(t)$, then $b_s$ was added to $\pi$ before $a_t$;
      it follows that $b_s$ belongs to one of the subsequences $\Gamma_1, \Gamma_2, \dots,
      \Gamma_{t-1}$.
    Hence, $b_s$ belongs to $\Area(\Cc, i, i+1)$ for some $i < t$.
    By (\ref{prefix-area-sum-eq}), we conclude that $b_s \in \Area(\Cc, 1, t)$.
    Conversely, if $b_s \in \Area(\Cc, 1, t)$, then $b_s \in \Area(\Cc, i, i+1)$ for some $i < t$.
    If we fix $i$ to be the smallest index for which $b_s \in \Area(\Cc, i, i+1)$, we will
      get that $b_s \in \Gamma_i$, where $i < t$, so also $\beta(s) < \alpha(t)$.
    We conclude that
    \begin{equation}
    \label{large-neighbor-separation-bs}
    b_s \in \Area(\Cc, 1, t)\ \Leftrightarrow\ \beta(s) < \alpha(t).
    \end{equation}

    Now, let us compare $\alpha(s)$ with $\alpha(t)$.
    Naturally, we have $\alpha(s) < \alpha(t)$ if and only if $s < t$.
    Also, by Lemma \ref{neighbor-cage-ai-good-lemma}
      we have that $a_s \in \mathrm{Int}\,\Area(\Cc, s-1, s+1)$.
    If $s < t$, we immediately get that $a_s \in \Area(\Cc, 1, t)$.
    However, if $s > t$, then we infer the following relation from Corollary~\ref{area-int-lemma}:
    $$ \Area(\Cc, 1, s-1) \cap \mathrm{Int}\,\Area(\Cc, s-1, s+1) = \varnothing. $$
    Hence, $a_s \not\in \Area(\Cc, 1, s - 1)$, so also $a_s \not\in \Area(\Cc, 1, t)$.
    We now proved that
    \begin{equation}
    \label{large-neighbor-separation-as}
    a_s \in \Area(\Cc, 1, t)\ \Leftrightarrow\ \alpha(s) < \alpha(t).
    \end{equation}

    It is obvious that $\alpha(t) \in X_s$ if and only if exactly one of the conditions
      $\alpha(s) < \alpha(t)$, $\beta(s) < \alpha(t)$ holds.
    From (\ref{large-neighbor-separation-bs}) and (\ref{large-neighbor-separation-as})
      it follows that this is equivalent to exactly one of the conditions
      $a_s \in \Area(\Cc, 1, t)$, $b_s \in \Area(\Cc, 1, t)$ being satisfied.
    A verification with Definition \ref{separated-from-def} concludes the proof.
    \end{claimproof}
  \end{claim}

  This claim enables us to prove that the conflict intervals do not intersect \emph{too frequently}:
  
  \begin{claim}
    \label{large-neighbor-sparse-conflicts}
    For each $x \in \mathbb{R}$, at most $256d^3(d+2)^4 + 2$ conflict intervals contain $x$.
    \begin{claimproof}
    \cqed
    Fix $x \in \mathbb{R}$, and consider the set $I = \{i \in [2, L-1]\,\mid\, x \in X_i\}$ of all indices
      for which the corresponding conflict interval contains $x$.
    We partition $I$ into two subsets, depending on whether $\alpha(i)$ is the left or the right
      endpoint of the conflict interval $X_i$:
    \begin{equation*}
    \begin{split}
      I_{\alpha\beta} &= \{i \in I\,\mid\, \alpha(i) < \beta(i)\}, \\
      I_{\beta\alpha} &= \{i \in I\,\mid\, \alpha(i) > \beta(i)\}.
    \end{split}
    \end{equation*}
    Let $m = |I_{\alpha\beta}|$ and $I_{\alpha\beta} = \{i_1, i_2, \dots, i_m\}$ where
      $2 \leq i_1 < i_2 < \dots < i_m \leq L-1$.
    Since $\alpha(i_1) < \alpha(i_2) < \dots < \alpha(i_m)$ and each interval $X_{i_j}$ for
      $j \in [1, m]$ has its left endpoint at $\alpha(i_j)$, we conclude that
      $\alpha(i_m)$ also belongs to each conflict interval $X_{i_1}, X_{i_2}, \dots, X_{i_m}$.
    Hence, by Claim \ref{large-neighbor-separation-claim}, each pair of vertices
      $a_{i_j}, b_{i_j}$ for $j \in [1, m - 1]$ is separated by $\Area(\Cc, 1, i_m)$.
    By Lemma \ref{neighbor-separation-lemma} we conclude that $m - 1 \leq
      128d^3(d+2)^4$, so $|I_{\alpha\beta}| \leq 128d^3(d+2)^4 + 1$.
      
    We can prove in a similar fashion that $|I_{\beta\alpha}| \leq 128d^3(d+2)^4 + 1$ ---
      if $I_{\beta\alpha} = \{j_1, j_2, \dots, j_n\}$ for some $n \in \mathbb{N}$
      and $2 \leq j_1 < j_2 < \dots < j_n \leq L-1$,
      then $\alpha(j_1)$ belongs to each conflict interval $X_{j_1}, X_{j_2}, \dots, X_{j_n}$,
      and the remaining part of the argument follows.
      
    Therefore, $|I| = |I_{\alpha\beta}| + |I_{\beta\alpha}| \leq 256d^3(d+2)^4 + 2$.
    \end{claimproof}
  \end{claim}
  
  We now consider an auxiliary undirected graph $H$ where $V(H) = \{2, 3, \dots, L-1\}$ and
    vertices $i$ and $j$ are connected by an edge if the conflict intervals $X_i$ and $X_j$
    intersect.
  This is a well-known definition of an interval graph.
  Every interval graph is also a perfect graph \cite{west2001introduction},
    so $H$ has a proper coloring where the number of colors is equal
    to the maximum size of a clique in $H$.
  In interval graphs, the maximum size of a clique is equal to the maximum number of intervals
    which intersect at a single point.
  This number is bounded by Claim \ref{large-neighbor-sparse-conflicts} by
    $256d^3(d+2)^4 + 2$,
    and therefore $H$ can be colored using that many colors.
  Since $L = [256d^3(d+2)^4 + 2] \ell + 2$, we have
    $|V(H)| = [256d^3(d+2)^4 + 2] \ell$ and there exists a color that has been used for at least
    $\ell$ vertices of $H$.
  Therefore, $H$ contains an independent set of size $\ell$, or in other words,
    there exists a set of $\ell$ pairwise disjoint conflict intervals
    $X_{d_1}, X_{d_2}, \dots, X_{d_\ell}$ ($2 \leq d_1 < d_2 < \dots < d_\ell \leq L-1$).
  
  \vspace{0.5em}
  We create a modified identity ordered cage $\Cc'$ which is a subset of $\Cc$
    given by indices $d_1, d_2, \dots, d_\ell$.
    
  \begin{claim}
    \label{large-neighbor-found}
    $\Cc'$ is a neighbor cage.
    \begin{claimproof}
    \cqed
  We fix an index $i \in \{2, 3, \dots, \ell - 1\}$, and we prove that
    $b_{d_i}$ belongs to $\Area(\Cc, d_{i - 1}, d_{i + 1})$.
  Naturally, $\Area(\Cc, d_{i-1}, d_{i+1}) = \Area(\Cc', i-1, i+1)$, so this is enough to
    complete the proof of the claim.

  Lemma \ref{neighbor-cage-ai-good-lemma} asserts that $a_{d_i} \in \mathrm{Int}\,\Area(
    \Cc, d_i-1, d_i+1) \subseteq \mathrm{Int}\,\Area(\Cc, d_{i-1}, d_{i+1})$.
  Hence, $a_i \in \Int\Area(\Cc, 1, d_{i+1})$.
  Moreover, from Corollary \ref{area-int-lemma} we deduce that
    $\Area(\Cc, 1, d_{i-1}) \cap \Int\Area(\Cc, d_{i-1}, d_{i+1}) = \varnothing$.
  Hence, $a_i \not\in \Area(\Cc, 1, d_{i-1})$.
  
  Now, we apply Claim \ref{large-neighbor-separation-claim}.
  Since the conflict intervals $X_{d_1}, X_{d_2}, \dots, X_{d_\ell}$ are pairwise disjoint,
    we get that $\alpha(d_i) \not\in X_{d_{i-1}}$ and $\alpha(d_i) \not\in X_{d_{i+1}}$, and therefore
    $a_{d_i}$ is separated from $b_{d_i}$ by neither $\Area(\Cc, 1, d_{i-1})$ nor
    $\Area(\Cc, 1, d_{i+1})$.
  Due to the fact that $a_{d_i} \in \Int\Area(\Cc, d_{i-1}, d_{i+1})$, we get that
  $b_{d_i} \not\in \Area(\Cc, 1, d_{i-1})$ and $b_i \in \Area(\Cc, 1, d_{i+1})$.
  Now, we observe that
    $$\Area(\Cc, 1, d_{i + 1}) = \Area(\Cc, 1, d_{i-1}) \cup \Area(\Cc, d_{i-1}, d_{i+1})$$
    (Lemma \ref{area-sum-lemma}).
  Since $b_{d_i} \in \Area(\Cc, 1, d_{i+1})$ and $b_{d_i} \not\in \Area(\Cc', 1, d_{i-1})$, we
    infer that $b_{d_i} \in \Area(\Cc, d_{i-1}, d_{i+1})$.
  Hence, $\Cc'$ is a neighbor cage of order $\ell$.
    \end{claimproof}
  \end{claim}

  Claim \ref{large-neighbor-found} completes the proof.
  \end{proof}
\end{lemma}

Let us now verify a few basic properties of neighbor cages.

\begin{lemma}
  \label{neighbor-cage-spread-lemma}
  Fix a neighbor cage $\Cc = \NeighborCage((a_i)_{i=1}^\ell, (b_i)_{i=1}^\ell, \Pc, \Qc)$ in
    a graph $G$ embedded in the plane, and consider an edge $uv \in E(G)$.
  Assume that for some $L, R$ ($2 \leq L < R < \ell$), we have that
    $u \in \Area(\Cc, L, R)$, but $u \not\in V(\mathcal{Q}) \cup \{p\}$.
  Then $v \in \Area(\Cc, L-1, R+1)$.
  \begin{proof}
  Let $\Pc = \Tree(P_1, P_2, \dots, P_\ell)$.  
  
  Since $\Area(\Cc, L, R) \subseteq \Area(\Cc, L-1, R+1)$ (a consequence of Lemma
    \ref{area-sum-lemma}),
    we have that $u \in \Area(\Cc, L-1, R+1)$.
  We now want to prove that $u \in \Int\Area(\Cc, L-1, R+1)$.

  Assume on the contrary that $u \in \partial\Area(\Cc, L-1, R+1)$.
  As $\Area(\Cc, L, R) \subseteq \Area(\Cc, L-1, R+1)$ and $u \in \Area(\Cc, L, R)$,
    it follows that $u \in \partial\Area(\Cc, L, R)$.
  Hence,
  $$ u \in \left[\partial\Area(\Cc, L-1, R+1) \cap \partial\Area(\Cc, L, R)\right] \setminus \left[
    V(\Qc) \cup \{p\}\right]. $$
  
  Since $\partial\Area(\Cc, L, R)\,\subseteq\,V(\Pc) \cup V(\Qc)$, we infer from the definition
    of a neighbor area that
    $$u \in \partial\Area(\Cc, L, R) \setminus V(\Qc)\quad\Rightarrow\quad
      u \in \partial\Area(\Cc, L, R) \cap V(\Pc)\quad\Rightarrow\quad
      u \in P_L \cup P_R. $$
  Analogously, $u \in P_{L-1} \cup P_{R+1}$.
  
  However, as $\Pc$ is a simple geodesic tree,
    every pair of different paths $P_i, P_j$ ($1 \leq i, j \leq \ell$) intersects only at $p$.
  Hence, $u = p$ --- a contradiction, since we assumed that $u \neq p$.
  \end{proof}
\end{lemma}

Therefore, short paths can only reach very local parts of the neighbor cage without
  intersecting $\mathcal{Q}$ or passing through $p$.
This fact is explored in the following lemmas:

\begin{lemma}
\label{semiladder-path-no-root}
  Fix a neighbor cage $\Cc = \NeighborCage((a_i)_{i=1}^\ell, (b_i)_{i=1}^\ell, \Pc, \Qc)$
    and two indices $i, j$ ($1 \leq i < j \leq \ell$).
  No shortest path between $b_i$ and $a_j$ can pass through $\Root(\Pc)$.
  \begin{proof}
    If any shortest path connecting $b_i$ and $a_j$ contains $\Root(\Pc)$ as a vertex, then
    \begin{equation*}
    \begin{split}
    d &< \dist(b_i, a_i) \leq \dist(b_i, \Root(\Pc)) + \dist(\Root(\Pc), a_i) = \\
      &= \dist(b_i, \Root(\Pc)) + \dist(\Root(\Pc), a_j) = \dist(b_i, a_j) \leq d
    \end{split}
    \end{equation*}
    --- a contradiction.  
  \end{proof}
\end{lemma}

\begin{lemma}
\label{path-first-intersection-lemma}
  Fix a neighbor cage $\Cc = \NeighborCage((a_i)_{i=1}^\ell, (b_i)_{i=1}^\ell, \Pc, \Qc)$ in $G$
    and two indices $i, j$ ($i < j$, $i \in [d+2, \ell-d-1]$).
  Take an oriented shortest path $R$ from $b_i$ to $a_j$.
  Note that this path intersects $\mathcal{Q}$ because $a_j \in V(\mathcal{Q})$.
  The first intersection of $R$ with $\mathcal{Q}$ belongs
    to $\Area(\Cc, i-d-1, i+d+1)$.
  \begin{proof}
    From the definition of a neighbor cage, we have that $b_i \in \Area(\Cc, i-1, i+1)$.
    Let $R = (v_0, v_1, \dots, v_\delta)$ be the shortest path from $b_i$ to $a_j$,
      where $v_0 = b_i$, $v_\delta = a_j$, and $\delta \leq d$ (as required by the definition
      of a semi-ladder).
    Let also $v_k$ be the first vertex of this path belonging to $V(\mathcal{Q})$.
    Since $p \not\in (v_0, v_1, \dots, v_\delta)$ (Lemma \ref{semiladder-path-no-root}),
      by applying Lemma~\ref{neighbor-cage-spread-lemma}
      we inductively deduce that $v_t \in \Area(\Cc, i-t-1, i+t+1)$ for each $t \in \{0, 1, \dots, k\}$.
    Therefore, $v_k \in \Area(\Cc, i-k-1, i+k+1)$.
    As $k \leq \delta \leq d$, we infer that $v_k \in \Area(\Cc, i-d-1, i+d+1)$.
  \end{proof}
\end{lemma}

Now, for a vertex $v$ belonging to a geodesic tree $\Tc$, we define
  the \textbf{depth} $\mu_\Tc(v)$ of $v$ as the distance from $v$ to $\Root(\Tc)$ in $\Tc$.
Note that it is implied by the definition of a geodesic tree that $\mu_\Tc(v)$ is also
  equal to the distance between $v$ and $\Root(\Tc)$ in $G$.

We remark a few simple facts about $\mu_\Tc$:
\begin{itemize}
  \item $\mu_\Tc(\Root(\Tc)) = 0$.
  \item $\mu_\Tc$ is constant on the set of leaves of $\Tc$.
  \item If a vertex $u$ is an ancestor of another vertex $v$ in $\mathcal{T}$, then
    $\dist(u, v) = \mu_\Tc(v) - \mu_\Tc(u)$.
  \item More generally, for any two vertices $u, v \in V(\mathcal{T})$,
    we have $\dist(u, v) \geq |\mu_\Tc(u) - \mu_\Tc(v)|$.
  \item On any oriented simple path in $\Tc$, the depths of vertices are first decreasing,
    and then increasing. In other words, there do not exist three vertices $x, y, z$ lying
    on an oriented simple path in $\Tc$ in this order,
    for which $\mu_\Tc(x) < \mu_\Tc(y) > \mu_\Tc(z)$.
\end{itemize}

\vspace{0.5em}
We already know from Lemmas \ref{semiladder-path-no-root} and
  \ref{path-first-intersection-lemma} that each shortest path $R$ connecting the vertices
  $b_i$ and $a_j$ of the semi-ladder ($i < j$) must avoid $\Root(\Pc)$ and the first intersection
  of $R$ with $\Qc$ occurs inside of $\Area(\Cc, i-d-1, i+d+1)$.
It turns out that this first intersection imposes serious restrictions on each following intersection
  of $R$ with $\Qc$.

\begin{lemma}
\label{intersection-depth-lower-bound-lemma}
Fix a neighbor cage $\Cc = \NeighborCage((a_i)_{i=1}^\ell, (b_i)_{i=1}^\ell, \Pc, \Qc)$
  and two indices $i, j$ ($1 \le i < j \le \ell$).
Consider any shortest path from $b_i$ to $a_j$: $R = (v_0, v_1, v_2, \dots, v_\delta)$ for some
  $\delta \leq d$, $v_0 = b_i$, $v_\delta = a_j$.
Suppose that for some indices $0 \leq x \leq y \leq \delta$ we have $v_x, v_y \in V(\mathcal{Q})$.
Then $\mu_\Qc(v_y) > \mu_\Qc(\mathrm{lca}_\Qc(v_x, a_i))$.
  \begin{proof}
  We set $\mu := \mu_\Qc$, $\lca := \lca_\Qc$.
  
  For the sake of contradiction let us fix a shortest path $R$
    from $b_i$ to $a_j$ and two indices $x, y$ as above such that
    $\mu(v_y) \leq \mu(\lca(v_x, a_i))$.
  We consider a path $G$ from $b_i$ to $a_i$ which first follows $R$ from $b_i$ to $v_x$,
    then goes up the tree $\Qc$ from $v_x$ to $\lca(v_x, a_i)$, and then goes down the tree
    $\Qc$ from $\lca(v_x, a_i)$ to $a_i$ (Figure \ref{intersection-depth-setup-fig}).
    
  \begin{figure}[h]
    \centering
    \begin{minipage}[b]{0.49\textwidth}
        \centering
        \input{figures/depth-lca-a.tex}
        (a)
      \end{minipage}\begin{minipage}[b]{0.49\textwidth}
        \centering
        \input{figures/depth-lca-b.tex}
        (b)
      \end{minipage}
      \caption{The setup in Lemma \ref{intersection-depth-lower-bound-lemma}. \\
      (a) --- The shortest path $R$ from $b_i$ to $a_j$ (red) crosses $\Qc$ twice.
        We are going to prove
        that there exists a path $G$ from $b_i$ to $a_i$ (green) which is not longer than $R$. \\
      (b) --- Both $R$ and $G$ have been split into three parts: red, blue and green. $R$ and
        $G$ share the red fragment. (\ref{intersection-depth-eq-b}) proves that the blue segment
        in $G$ is not longer than the blue segment in $R$. (\ref{intersection-depth-eq-c}) proves
        that the green segment in $G$ is not longer than the green segment in $R$.}
      \label{intersection-depth-setup-fig}
    \end{figure}

  We obviously have that
  \begin{equation}
  \label{intersection-depth-main-eq}
  \dist(b_i, a_i) \leq \len(G) = \len(G[b_i, v_x]) + \len(G[v_x, \lca(v_x, a_i)]) + 
    \len(G[\lca(v_x, a_i), a_i]).
  \end{equation}

  As $G$ follows $R$ on the fragment between $b_i$ and $v_x$, we have that
  \begin{equation}
  \label{intersection-depth-eq-a}
  \len(G[b_i, v_x]) = \len(R[b_i, v_x]).
  \end{equation}

  Next, $\lca(v_x, a_i)$ is an ancestor of $v_x$ in $\Qc$ and thus
  \begin{equation}
  \label{intersection-depth-eq-b}
  \begin{split}
  \len(G[v_x, \mathrm{lca}(v_x, a_i)]) & = \mu(v_x) - \mu(\mathrm{lca}(v_x, a_i)) \leq
    \mu(v_x) - \mu(v_y) \leq \dist(v_x, v_y) = \\ &= \mathrm{len}(R[v_x, v_y]).
  \end{split}
  \end{equation}

  Analogously, $\mathrm{lca}(v_x, a_i)$ is an ancestor of $a_i$ in $\mathcal{Q}$.
  Therefore,
  \begin{equation}
  \label{intersection-depth-eq-c}
  \begin{split}
  \len(G[\mathrm{lca}(v_x, a_i), a_i]) & = \mu(a_i) - \mu(\mathrm{lca}(v_x, a_i)) \leq
    \mu(a_i) - \mu(v_y) = \mu(a_j) - \mu(v_y) \leq \\
    & \leq \dist(v_y, a_j) = \mathrm{len}(R[v_y, a_j]).
  \end{split}
  \end{equation}
  
  By substituting (\ref{intersection-depth-eq-a}), (\ref{intersection-depth-eq-b}),
    (\ref{intersection-depth-eq-c}) into (\ref{intersection-depth-main-eq}), we conclude that
  \begin{equation*}
  \begin{split}
  \dist(b_i, a_i) & \leq \mathrm{len}(R[b_i, v_x]) + \mathrm{len}(R[v_x, v_y]) +
    \mathrm{len}(R[v_y, a_j]) = \mathrm{len}(R) = \dist(b_i, a_j) = \\ &= \delta \leq d
  \end{split}
  \end{equation*}
  --- a contradiction since the semi-ladder requires that $\dist(b_i, a_i) > d$.
  \end{proof}
\end{lemma}

A less formal way of thinking about Lemma \ref{intersection-depth-lower-bound-lemma}
  is that as soon as the shortest path from $b_i$ to $a_j$ touches any vertex of $\mathcal{Q}$,
  all subsequent intersections of the path with $\mathcal{Q}$ need to happen
  sufficiently deep in the tree.
There are also some limitations on where this path may intersect $\mathcal{Q}$ at all ---
  by letting $x = y$ in the lemma above, we infer that each intersection $v_x$
  of a shortest path from $b_i$ to $a_j$ must satisfy $v_x \neq \lca_\Qc(v_x, a_i)$.
That is, no intersection of this path with $\Qc$ can be an ancestor of $a_i$ in $\Qc$.

The most important conclusion will become clear in a moment:
  since the first intersection of every oriented shortest path from $b_i$ to $a_j$ ($i < j$)
  with $\mathcal{Q}$ must happen topologically relatively close to $b_i$
  (Lemma \ref{path-first-intersection-lemma}),
  we should be able to infer a common lower bound on the depth of
  every such intersection, no matter the path or value of $j$.
This is indeed possible, albeit with slightly stronger assumptions,
  as proved in Lemma \ref{ladder-path-depth-bound-lemma} below.

In order to proceed, we need a small definition.
\begin{definition}
In a neighbor cage $\Cc = \NeighborCage((a_i)_{i=1}^\ell, (b_i)_{i=1}^\ell, \Pc, \Qc)$,
  a neighbor area $\Area(\Cc, i, j)$ is \textbf{rootless} if $\Root(\Qc) \not\in \Int\Area(\Cc, i, j)$.
  
Note that not all neighbor areas must be rootless, as seen in Figure \ref{rootless-def-figure}.
    \begin{figure}[h]
    \centering
    \begin{minipage}[b]{0.4\textwidth}
        \centering
        \input{figures/rootless-def-a.tex}
        (a)
      \end{minipage}\begin{minipage}[b]{0.4\textwidth}
        \centering
        \input{figures/rootless-def-b.tex}
        (b)
      \end{minipage}
    \caption{(a) --- An example neighbor cage $\Cc$ of order $5$. $\Area(\Cc, 3, 4)$ (green) is
      rootless, but $\Area(\Cc, 1, 3)$ (red) is not. Vertices $b_1, \dots, b_5$ were omitted
      for clarity. \\
      (b) --- Another neighbor cage $\Dc$. Now, both $\Area(\Dc, 1, 3)$ and $\Area(\Dc, 3, 4)$ are
      rootless.}
      \label{rootless-def-figure}
    \end{figure}
\end{definition}

We are now ready to introduce the promised lemma.

\begin{lemma}
\label{ladder-path-depth-bound-lemma}
Fix a neighbor cage $\Cc = \NeighborCage((a_i)_{i=1}^\ell, (b_i)_{i=1}^\ell, \Pc, \Qc)$
  and an index $i$ ($i \in [d + 2, \ell - d - 1]$). If $\Area(\Cc, i-d-1, i+d+1)$ is rootless, then
  each intersection $v$ of an oriented shortest path from $b_i$ to $a_j$ ($j > i$)
  with $\Qc$ satisfies
$$ \mu_\Qc(v) > \mu_\Qc(\lca_\Qc(a_{i-d-1}, a_{i+d+1})). $$
\begin{proof}
  We set $\mu := \mu_\Qc$, $\lca := \lca_\Qc$.
  We also let $v_{\lca} := \lca(a_{i-d-1}, a_{i+d+1})$ and $\Qc = \Tree(Q_1, \dots, Q_\ell)$.
  
  Basing on the definition of neighbor areas, we infer that $v_{\lca}$ is the shallowest vertex
    in $\Qc$ belonging to $\partial\Area(\Cc, i-d-1, i+d+1)$.
  We will now prove that $v_{\lca}$ is also the shallowest vertex in $\Qc$ belonging to
    $\Area(\Cc, i-d-1, i+d+1)$.
  Assume for contradiction that there exists a vertex $x \in \Area(\Cc, i-d-1, i+d+1) \cap V(\Qc)$
    for which $\mu(x) < \mu(v_{\lca})$.
  Since $\Root(\Qc) \not\in \Int\Area(\Cc, i-d-1, i+d+1)$,
    the unique simple path in $\Qc$ connecting $x$ with $\Root(\Qc)$ must intersect
    $\partial\Area(\Cc, i-d-1, i+d+1)$ at some vertex $s \in V(\Qc)$.
  Obviously, $\mu(s) \geq \mu(v_{\lca})$.
  However, it means that the simple path in $\Qc$ connecting $x$ with $\Root(\Qc)$
    contains three vertices in the following order: $x, s, \Root(\Qc)$. But
  $$ \mu(x) < \mu(v_{\lca}) \leq \mu(s) > \mu(\Root(\Qc)) = 0 \qquad
    \Rightarrow \qquad \mu(x) < \mu(s) > \mu(\Root(\Qc)) $$
  --- a contradiction.
  Hence, $v_{\lca}$ is the shallowest vertex of $\Area(\Cc, i-d-1, i+d+1) \cap V(\Qc)$.
  
  \vspace{0.5em}
  We fix a shortest path $R$ from $b_i$ to $a_j$. Let $v_1$ be the first intersection of $R$
    with $V(\mathcal{Q})$.
  By Lemma \ref{path-first-intersection-lemma}, we have that
    $v_1 \in \Area(\Cc, i-d-1, i+d+1)$.
    
  \begin{figure}[h]
    \centering
    \input{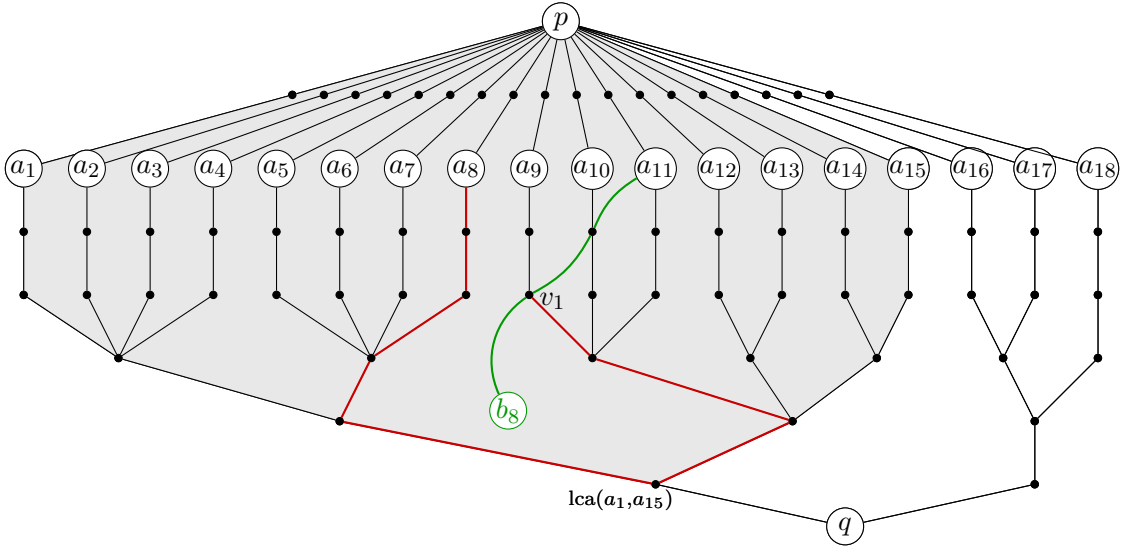}
    \caption{An example of the setup in Lemma \ref{ladder-path-depth-bound-lemma} ---
      a cage $\Cc$, $d = 6$ and $i = 8$.
      We consider $v_1$ ---
        the first intersection of the shortest path between $b_8$ and $a_{11}$ (green) with $\Qc$.
      We prove that the simple path $S$ in $\Qc$ connecting $v_1$ with $a_8$ (red) fully belongs
        to $\Area(\Cc, 1, 15)$ (shaded).
      It will follow that $\mu(\lca(v_1, a_8)) \geq \mu(v_{\lca})$, and Lemma
        \ref{intersection-depth-lower-bound-lemma} will conclude the proof.}
    \label{common-bound-setup-fig}
  \end{figure}

  We now prove (Figure \ref{common-bound-setup-fig})
    that the simple path in $\Qc$ between $v_1$ and $a_i$ --- call it $S$ --- is fully contained
    within $\Area(\Cc, i-d-1, i+d+1)$.
  Firstly, from the properties of neighbor cages, it follows that
    $a_i \in \Area(\Cc, i-1, i+1)$; thus, $a_i \in \Area(\Cc, i-d-1, i+d+1)$.
  Secondly, by the definition of neighbor areas, the intersection of
    $\partial\Area(\Cc, i-d-1, i+d+1)$ with $\Qc$ is the simple path in $\Qc$
    connecting $\Split(\Cc, i-d-1)$ with $\Split(\Cc, i+d+1)$.
  Let us call this path $B$.
  Therefore, if $S$ was not contained
    within $\Area(\Cc, i-d-1, i+d+1)$, it would have to intersect $R$ at least twice ---
    once when leaving the neighbor area, and once when reentering it.
  This would, however, mean that $\Qc$ contains a cycle --- which is impossible.

  Hence, $S$ fully belongs to $\Area(\Cc, i-d-1, i+d+1)$.
  It means in particular that
  \begin{equation*}
  \begin{split}
  \mu(\lca(v_1, a_i)) &=
    \min \{\mu(w)\,\mid\, w \in S\} \geq \\
    &\geq \min\{\mu(w)\,\mid\, w \in \Area(\Cc, i-d-1, i+d+1)\,\cap\,V(\Qc)\} = \mu(v_{\lca}).
  \end{split}
  \end{equation*}

  Therefore, by Lemma \ref{intersection-depth-lower-bound-lemma},
    each intersection $v$ of $R$ with $\Qc$ --- either $v_1$, or any
    later intersection --- must satisfy
  $$ \mu(v) > \mu(\lca(v_1, a_i)) \geq \mu(v_{\lca}). $$
  \end{proof}
\end{lemma}

\subsection{Separating cages}
\label{separating-cages-lemma}

Lemma \ref{ladder-path-depth-bound-lemma}
  allows us to introduce the final variant of a cage, which will eventually lead us
  to a polynomial bound on the maximum order of any neighbor cage.
This in turn will produce an upper bound on the maximum order of any distance-$d$ semi-ladder.

\begin{definition}
  A neighbor cage $\Cc = \NeighborCage((a_i)_{i=1}^\ell, (b_i)_{i=1}^\ell, \Pc, \Qc)$ in a graph $G$
    is a \textbf{separating cage} if there exists an integer $\lambda \in \{0, 1, 2, \dots, d - 1\}$,
    called \textbf{threshold}, such that the following properties hold
    (Figure~\ref{separating-def-fig}):
  \begin{itemize}
    \item for every two different indices $i, j$ ($1 \leq i, j \leq \ell$),
      paths $\Path(\Qc, i)$ and $\Path(\Qc, j)$ do not share any vertices
      which are at depth larger than $\lambda$ in $\Qc$;
    \item for every pair of indices $i, j$ ($1 \leq i < j \leq \ell$),
      no shortest path from $b_i$ to $a_j$ contains any vertex belonging to $\Qc$ which is at
      depth at most $\lambda$ in $\Qc$.
  \end{itemize}
  
  We denote a separating cage in the following way:
  $$ \SeparatingCage((a_i)_{i=1}^\ell, (b_i)_{i=1}^\ell, \Pc, \Qc, \lambda). $$
   \begin{figure}[h]
      \centering
      \begin{minipage}[b]{0.495\textwidth}
        \centering
        \input{figures/separating-def-a.tex}
        (a)
      \end{minipage}\begin{minipage}[b]{0.495\textwidth}
        \centering
        \input{figures/separating-def-b.tex}
        (b)
      \end{minipage}
      \caption{An example separating cage $\Cc$ of order $9$ and threshold $\lambda = 2$.
       The threshold is a blue dotted
        line; no two paths connecting $\Root(\Qc)$ with a leaf of $\Qc$ can intersect above the
        line, and no path connecting $b_i$ with $a_j$ (for $i < j$)
        can intersect $\Qc$ at this line or below. \\
        (a) --- Three valid paths connecting $b_i$ with $a_j$ for $i < j$ (green). \\
        (b) --- Two invalid paths connecting $b_i$ with $a_j$ for $i < j$ (red).
          These cannot occur in the separating cage because
          each of them passes through a vertex of $\Qc$ at depth $\lambda$ or less.}
      \label{separating-def-fig}
    \end{figure}
\end{definition}

We will now prove that graphs with large neighbor cages also contain large separating cages.

\begin{lemma}
  \label{large-separating-cage-lemma}
  Each neighbor cage of order $2(2d+3)((\ell-1)d + 1)$ contains a separating cage
    of order $\ell$ as a subset.
  \begin{proof}
  Let $L := 2(2d+3)((\ell-1)d + 1)$ and
  fix a cage $\Cc = \NeighborCage((a_i)_{i=1}^L,
    (b_i)_{i=1}^L, \Pc, \Qc)$ of order $L$.
  We let $\Qc = \Tree(Q_1, Q_2, \dots, Q_L)$.
  
  Let also $M := 2((\ell - 1)d + 1)$. We define a sequence of $M$ neighbor areas
    $(A_1, A_2, \dots, A_M)$:
  $$ A_i = \Area(\Cc,\, (2d + 3)(i - 1) + 1,\, (2d + 3)i) \qquad \text{for each index }i \in [1, M]. $$
  We remark that because the set $\{[(2d + 3)(i-1) + 1,\, (2d+3)i]\,\mid\, i \in \{1, \dots, M\}\}$
    contains only pairwise disjoint segments,
    we infer from Lemma \ref{area-int-lemma} that the interiors of the neighbor areas
    $A_1, A_2, \dots, A_M$ are pairwise disjoint.
  Hence, at least $M - 1 = 2(\ell - 1)d + 1$ out of neighbor areas $A_1, A_2, \dots, A_M$
    are rootless.
  
  For each neighbor area $A_i$, we define its depth $\mu(A_i) = \mu_\Qc(\lca_\Qc(
    a_{(2d + 3)(i - 1) + 1},\, a_{(2d+3)i}))$.
  Note that $\mu(A_i)$ is chosen specifically to match the conditions of Lemma
    \ref{ladder-path-depth-bound-lemma}:
    provided that $A_i$ is rootless, each intersection $v$ of a shortest path from
    $b_{(2d+3)(i-1) + (d+2)}$ to $a_j$ (for $j > (2d+3)(i-1) + (d+2)$) with $\Qc$ satisfies
    $\mu_\Qc(v) > \mu(A_i)$.
  We also remark that $\mu(A_i) \in \{0, 1, 2, \dots, d-1\}$ for each $i \in [1, M]$.
    
  \newcommand{\Block}{\mathsf{Block}}

  The pigeonhole principle allows us to select $2\ell-1$ out of $(2\ell-2)d + 1$ rootless
    neighbor areas so that the selected neighbor areas have the same depth $\mu$.
  We fix this subset of rootless neighbor areas: $A_{t_1}, A_{t_2}, \dots, A_{t_{2\ell-1}}$ where
    $1 \leq t_1 < t_2 < \dots < t_{2\ell-1} \leq M$. Note that
  $$ \mu(A_{t_1}) = \mu(A_{t_2}) = \dots = \mu(A_{t_{2\ell-1}}) =: \lambda, \qquad
    \Root(\Qc) \not\in \Int\Area(A_{t_j})\text{ for }j\in[1,2\ell-1]. $$
  We also define $m_j$ for $j \in [1, 2\ell-1]$ as the index of the central vertex of the semi-ladder
    in the neighbor area $A_{t_j}$: $m_j = (2d + 3)(t_j - 1) + (d + 2)$.
  Naturally, the sequence of indices $m_j$ is increasing: $m_1 < m_2 < \dots < m_{2\ell-1}$.
  Also, for every pair of indices $i, j \in [1, 2\ell - 1]$ ($i < j$), each intersection $v$ of
    a shortest path from $b_{m_i}$ to $a_{m_j}$ with $\Qc$ satisfies $\mu_\Qc(v) > \lambda$.

  \vspace{0.5em}
  Let $\Cc'$ be the neighbor cage of order $\ell$ which is a subset of $\Cc$ given by
    the sequence of indices $m_1, m_3, m_5, \dots, m_{2\ell-1}$:
  $$ \Cc' = \NeighborCage((a_{m_1}, a_{m_3}, \dots, a_{m_{2\ell-1}}),
    (b_{m_1}, b_{m_3}, \dots, b_{m_{2\ell-1}}), \Pc', \Qc'), \quad \Pc' \subseteq \Pc,\ 
    \Qc' \subseteq \Qc. $$
  We verify that $\Cc'$ is a separating cage with threshold $\lambda$.
  
  Firstly, we see that $\lambda = \mu(B_{t_1}) \in \{0, 1, 2, \dots, d-1\}$.

  Secondly, we want to prove that no shortest path whose shortness is required
    by the semi-ladder in $\Cc'$ contains any vertex belonging to $\Qc$ that is at
    depth at most $\lambda$ in $\Qc$.
  This is quickly resolved by Lemma \ref{ladder-path-depth-bound-lemma}.
  Fix any two vertices of the semi-ladder in $\Cc'$ which require a short path:
    $b_{m_i}, a_{m_j}$ ($i < j$, both $i, j$ odd), and fix a shortest path $R$ between
    these two vertices.
  As mentioned before, each intersection $v$ of $R$ with $\Qc$ satisfies $\mu_\Qc(v) > \lambda$.
  Since $\Qc' \subseteq \Qc$, we obviously have $\mu_{\Qc'}(v) > \lambda$ as well.
  
  Finally, we want to prove that no two different oriented paths connecting $\Root(\Qc')$
    with the leaves of $\Qc'$ intersect at any vertex at depth larger than $\lambda$.
  We again fix two vertices of the semi-ladder in $\Cc'$: $a_{m_i}, a_{m_j}$
    ($i < j$, both $i, j$ odd), and consider paths $Q_{m_i}$ and $Q_{m_j}$, connecting
    $\Root(\Qc)$ with $a_{m_i}$ and $a_{m_j}$, respectively.
  Assume for contradiction that these paths intersect at a vertex $v$ such that
    $\mu(v) > \lambda$.
  This means that $v \in Q_{m_i}$ and $v \in Q_{m_j}$.
  By the definition of identity ordered cages, we have that $v \in Q_k$ for every
    $k \in [m_i, m_j]$.

  However, since both $i$ and $j$ are odd, we have that $i + 1 < j$.
  This means that $$[m_{i+1}-d-1, m_{i+1}+d+1] \subseteq [m_i, m_j].$$
  Hence, $v \in Q_{m_{i+1}-d-1}$ and $v \in Q_{m_{i+1}+d+1}$ since $\Cc$ is an identity
    ordered cage.
  However, we verify that
  $$ \mu_\Qc\left(\lca_\Qc(a_{m_{i+1}-d-1}, a_{m_{i+1}+d+1})\right) =
    \mu(A_{t_{i+1}}) = \lambda. $$
  This means that paths $Q_{m_{i+1}-d-1}$ and $Q_{m_{i+1}+d+1}$ do not intersect
    at any vertex at depth greater than $\lambda$, yet $v$ is an intersection of these paths
    at depth greater than $\lambda$ --- a contradiction.
  Therefore, paths $Q_{m_i}, Q_{m_j}$ do not share any vertex at depth larger than $\lambda$.
  Since these paths were chosen arbitrarily from $\Qc'$,
    all the paths connecting $\Root(\Qc)$ with a leaf of $\Qc'$ are vertex-disjoint,
    as long as vertices at depth larger than $\lambda$ are concerned.

  Since all required properties of separating cages have been verified, $\Cc'$ together
    with threshold $\lambda$ is a separating cage.
  \end{proof}
\end{lemma}

In separating cages, we can prove a much more powerful variant of Lemma
  \ref{neighbor-cage-spread-lemma}:

\begin{lemma}
  \label{separating-cage-spread-lemma}
  Consider a separating cage
    $\Cc = \SeparatingCage((a_i)_{i=1}^\ell, (b_i)_{i=1}^\ell, \Pc, \Qc, \lambda)$.
  For every pair of indices $i, j$ ($1 \leq i < j \leq \ell$) and every edge $uv$ on a shortest
    path between $b_i$ and $a_j$, if $u \in \Area(\Cc, L, R)$ for some integers $2 \leq L < R < \ell$,
    then $v \in \Area(\Cc, L-1, R+1)$.

  \begin{proof}
    In a similar vein to the proof of Lemma \ref{neighbor-cage-spread-lemma}, we will prove that
      $u \not\in \partial\Area(\Cc, L-1, R+1)$; the correctness of the statement of the lemma
      will follow.

    We remark that $u \in \Area(\Cc, L, R)$, and thus $u \in \Area(\Cc, L-1, R+1)$.
    We assume for contradiction that $u \in \partial\Area(\Cc, L-1, R+1)$.
    Since $\Area(\Cc, L, R) \subseteq \Area(\Cc, L-1, R+1)$, we also have that
      $u \in \partial\Area(\Cc, L, R)$.
    Therefore, both conditions below must be satisfied:
    \begin{equation*}
    \begin{split}
      u \in \SPath(\Cc, L-1)\ \ &\vee\ \ u \in \SPath(\Cc, R+1), \\
      u \in \SPath(\Cc, L)\ \ &\vee\ \ u \in \SPath(\Cc, R).
    \end{split}
    \end{equation*}
    
    Hence, $v$ belongs to the intersection of $\SPath(\Cc, i)$ and $\SPath(\Cc, j)$
      for $i \neq j$.
    Therefore, by Lemma~\ref{splitting-path-intersection}, these paths intersect exactly at $p$
      and their common suffix belonging to $\Qc$.
    However, since $u$ lies on the shortest path whose shortness is required by the
      semi-ladder, Lemma \ref{semiladder-path-no-root} applies
      and asserts that $u \neq p$.
    It means that $u$ belongs to the intersection of paths $\Path(\Qc, i)$ and $\Path(\Qc, j)$.
    But the definition of a separating cage requires that each such intersection must occur
      at depth at most $\lambda$.
    This means that the considered shortest path intersects $\Qc$ at a vertex $u$
      located at depth at most $\lambda$, which is explicitly forbidden by the definition
      of a~separating cage --- a contradiction.
      
    As $u \in \Area(\Cc, L-1, R+1)$, $u \not\in \partial\Area(\Cc, L-1, R+1)$, and $v$ is
      connected to $u$ by an edge, we have that $v \in \Area(\Cc, L-1, R+1)$.
  \end{proof}
\end{lemma}

The lemma above leads to a straightforward linear bound on the maximum order of any
  separating cage:

\begin{lemma}
\label{no-large-separating-cage}
  Every separating cage has order smaller than $2d+5$.
  \begin{proof}
    Assume for contradiction
      that $\Cc = \SeparatingCage((a_i)_{i=1}^\ell, (b_i)_{i=1}^\ell, \Pc, \Qc, \lambda)$
      is a separating cage of order $\ell \geq 2d + 5$.
    We consider a shortest path from $b_{d+2}$ to $a_{2d+4}$: $(u_0, u_1, u_2, \dots, u_\delta)$
      where $u_0 = b_{d+2}$, $u_\delta = a_{2d+4}$ and $\delta \leq d$
      (as required by the semi-ladder).
      
    By the properties of neighbor cages, we have that $u_0 \in \Area(\Cc, d + 1, d + 3)$.
    Using Lemma \ref{separating-cage-spread-lemma},
      we prove inductively that $u_i \in \Area(\Cc, d - i + 1, d + i + 3)$ for each
      $i \in \{0, 1, \dots, \delta\}$.
    In particular, we infer that $a_{2d+4} = u_\delta \in \Area(\Cc, d - \delta + 1, d + \delta + 3)$.
    Since $\delta \leq d$, we get that $a_{2d+4} \in \Area(\Cc, 1, 2d + 3)$.
    
    However, Lemma \ref{neighbor-cage-ai-good-lemma} asserts
      that $a_{2d+4} \in \mathrm{Int}\,\Area(\Cc, 2d+3, 2d+5)$.
    Since
      $$\Area(\Cc, 1, 2d+3) \cap \mathrm{Int}\,\Area(\Cc, 2d+3, 2d+5) = \varnothing$$
    (Corollary \ref{area-int-lemma}), we have a contradiction.
  \end{proof}
\end{lemma}

We can now finalize the proof of Theorem \ref{planar-upper-bound} bounding the maximum
  distance-$d$ semi-ladder order in planar graphs:

\begin{corollary}
\label{planar-upper-bound-cor}
  For $d \geq 1$,
  all distance-$d$ semi-ladders in planar graphs have order smaller than
  $$ d \cdot \left\{ d \cdot (2d - 1) \cdot
    \left(\left\{[256d^3(d+2)^4 + 2] \cdot 2(2d+3)[(2d+4)d + 1]
    + 1\right\}^2 + 1\right)
    + 2 \right\}^d + 1. $$
  \begin{proof}
    We define the following polynomial functions:
    \begin{equation*}
    \begin{split}
      \chi_1(d) &= 2d + 5, \\
      \chi_2(d) &= 2(2d+3)[(\chi_1(d) - 1)d + 1], \\
      \chi_3(d) &= [256d^3(d+2)^4 + 2]\chi_2(d) + 2, \\
      \chi_4(d) &= (\chi_3(d) - 1)^2 + 1, \\
      \chi_5(d) &= (2d-1)\chi_4(d).
    \end{split}
    \end{equation*}
    It can be verified that the formula in the statement of the corollary simplifies to
      $$d \left\{ d \cdot \chi_5(d) + 2\right\}^d + 1 =: M_d.$$
    Assume that a distance-$d$ semi-ladder of order at least $M_d$ exists.
    We take its subset which is a semi-ladder of order $M_d$, and call it $\Cc_1$.
    
    Using Lemma \ref{quasi-cage-exists}, we find a quasi-cage $\Cc_2$ of order
      $\chi_5(d)$ in the graph.
      
    Using Lemma \ref{cage-exists}, we find a cage $\Cc_3$ of order $\chi_4(d)$ as
      a subset of $\Cc_3$.
      
    Using Lemma \ref{order-lemma}, we associate the cage $\Cc_3$ with an order and
      create an ordered cage $\Cc_4$ underlying the same semi-ladder as $\Cc_3$.
      
    Using Lemma \ref{identity-ordered-cage-exists}, we find an identity ordered cage $\Cc_5$ of
      order $\chi_3(d)$ as a subset of $\Cc_4$.
    Note that this may require altering the embedding of the graph in the plane.
      
    Using Lemma \ref{large-neighbor-lemma}, we find a neighbor cage $\Cc_6$ of
      order $\chi_2(d)$ as a subset of $\Cc_5$.
      
    Using Lemma \ref{large-separating-cage-lemma}, we find a separating cage $\Cc_7$ of
      order $\chi_1(d)=2d+5$ as a subset of $\Cc_6$.

    Finally, Lemma \ref{no-large-separating-cage} contradicts our assumption since
      it asserts that $\Cc_7$ cannot exist.
      
    Hence, all distance-$d$ semi-ladders in the class of planar graphs must have
      order smaller than $M_d$.
  \end{proof}
\end{corollary}

We remark that this upper bound is of the form $d \cdot \rho(d)^d + 1$ where
  $\rho(d)$ is a polynomial of degree $22$,
  so the upper bound is of the form $d^{O(d)}$.
This concludes the proof of Theorem \ref{planar-upper-bound}.

\section{Upper bounds in other classes of graphs}\label{other-upper-bounds-chapter}
In this section we present upper bounds on the maximum semi-ladder
  orders in other classes of graphs: graphs with bounded maximum degree,
  graphs with bounded pathwidth or treewidth, and graphs excluding the complete graph
  $K_t$ as a minor.
  
In Subsection \ref{degree-up-section}, we prove the $\Delta^d + 1$ upper bound for the class
  of graphs with the maximum degree bounded by $\Delta$.
The proof is straightforward, since in these graphs, no vertex can have too many
  other vertices in its distance-$d$ neighborhood.

In Subsection \ref{pw-up-section}, we bound the maximum semi-ladder order from above in
  graphs with bounded pathwidth.
As an important part of the proof, we generalize the sunflower lemma
  in a way that allows us to assign labels to the elements of each set of the sunflower.

In Subsection \ref{kt-up-section}, we prove the $d^{O(d^{t-1})}$ upper bound
  for the class of graphs excluding $K_t$ as a minor.
This proof will utilize the elements of the standard toolchain of Sparsity --- weak
  coloring numbers and uniform quasi-wideness --- and can be generalized to any class
  of graphs with bounded expansion.
However, in our proof, we also apply a new trick specific to the considered class of graphs
  (Lemma \ref{kt-reduce-uqw}), which allows us to prove a much tighter upper bound.
The same section also concludes with a similar upper bound in the class of graphs with
  bounded treewidth, facilitated by Theorem \ref{tw-minor-free}.

\subsection{Graphs of bounded degree}
\label{degree-up-section}

The following simple theorem presents an upper bound on the maximum semi-ladder order
  in graphs with maximum degree bounded by $\Delta$.
The main idea of the proof is that in a distance-$d$ semi-ladder $a_1, a_2, \dots, a_\ell,
  b_1, b_2, \dots, b_\ell$, vertex $b_1$ must remain at distance at most $d$
  from vertices $a_2, a_3, \dots, a_\ell$.

\begin{theorem}
\label{degree-upper-bound}
For each $\Delta \geq 2$, $d \geq 1$,
  the order of every distance-$d$ semi-ladder in a graph with maximum degree bounded by
  $\Delta$ is bounded from above by $\Delta^d + 1$.
  
  \begin{proof}
  Fix a graph $G$ in which the maximum degree of a vertex is bounded from above by $\Delta$,
    and let $a_1, a_2, \dots, a_\ell, b_1, b_2, \dots, b_\ell$ be a distance-$d$ semi-ladder in $G$
    of order $\ell \geq 2$.
  By the properties of semi-ladders, we infer that for each $i \in \{2, 3, \dots, \ell\}$, we have that
    $\dist(b_1, a_i) \leq d$ --- that is, there are $\ell - 1$ vertices (different than $b_1$)
    in the distance-$d$ neighborhood of $b_1$.
  Since each vertex $a_2, a_3, \dots, a_\ell$ is different than $b_1$
    (by the definition of a distance-$d$ semi-ladder in Section \ref{preliminaries-chapter}),
    it means that each of them is at distance at least $1$ and at most $d$ from $b_1$.

  Obviously, there are at most $\Delta(\Delta - 1)^k$ vertices at distance exactly $k$
    ($k \geq 1$) from $b_1$.
  Summing this over all $k \in \{1, 2, \dots, d\}$, we get that
  
  \begin{equation*}
  \begin{split}
  \ell - 1 & \leq \sum_{k = 1}^d \Delta (\Delta - 1)^{k-1} = \Delta + (\Delta-1) \sum_{k=2}^d
    \Delta(\Delta-1)^{k-2} \leq \Delta + (\Delta-1) \sum_{k=1}^{d-1} \Delta^k = \\
    & = 1 + (\Delta - 1) \sum_{k=0}^{d-1} \Delta^k = 1 + \frac{\Delta^d - 1}{\Delta - 1}(\Delta - 1) =
    \Delta^d.
  \end{split}
  \end{equation*}
  \end{proof}
\end{theorem}

We contrast the $\Delta^d + 1$ upper bound with the $\left\lfloor\frac{\Delta}{2}\right\rfloor^{
  \left\lceil \frac{d}{2} \right\rceil}$ lower bound (Theorem \ref{degree-lower-bound})
  proved in Section \ref{lower-bounds-chapter}.

\subsection{Graphs with bounded pathwidth}
\label{pw-up-section}

We begin this section by reminding about a classic combinatorial result, first proved
  by Erd\H{o}s and Rado \cite{sunflower_first}: the \textbf{Sunflower Lemma}.
  
\begin{definition}
\label{sunflower-def}
Suppose $\Omega$ is a universe and $\Fc = (F_1, F_2, \dots, F_n)$ is a family of
  (not necessarily different) subsets of $\Omega$.
Then $\Fc$ is a \textbf{sunflower} if there exists a set $C$ (which we call the \textbf{core} of
  the sunflower) such that for every two different indices $i, j \in [1, n]$, we have
  $F_i \cap F_j = C$.
The \textbf{order} of a sunflower is the number of sets in it, that is, $|\Fc| = n$.
\end{definition}

\begin{theorem}[Sunflower Lemma]
  \label{sunflower-lemma}
  Let $\Gc$ be a family of $b!a^{b+1}$ (not necessarily different) subsets of $\Omega$
    in which each subset has at most
  $b$ elements. Then $\Gc$ contains a subfamily $\Fc \subseteq \Gc$ which
  is a sunflower of order $a$.
\end{theorem}

The ``textbook proof'' of this lemma can be found e.g. in the book by Cygan et al.
  \cite{ParamAlgo}.
  
\medskip

In order to prove an upper bound on the maximum semi-ladder order in the class of graphs with
  bounded pathwidth, we first need to generalize the Sunflower Lemma.

\subsubsection*{Labeled Sunflower Lemma}

We wish to generalize the Sunflower Lemma (Theorem \ref{sunflower-lemma})
  by allowing the sets in the family to assign \textbf{labels} to the elements of the set.
We first need to describe formally what this labeling means.

\begin{definition}
  For a finite set of labels $\Sigma$ and a universe $\Omega$,
  a \textbf{labeled subset} $A \subseteq \Omega$ over $\Sigma$
    is a subset of $\Omega$ in which each element is also assigned a
    label from $\Sigma$.
  The \textbf{cardinality} of $A$ is the number of elements belonging to $A$.
  
  In other words, a labeled subset is a partial function $A : \Omega \rightharpoonup \Sigma$,
  and the cardinality of $A$ is the size of its domain.
\end{definition}

The definition of a sunflower now naturally generalizes to the families of labeled subsets:
  
\begin{definition}
  For a finite set of labels $\Sigma$ and a universe $\Omega$, a family $\Fc$ of
    (not necessarily different) labeled subsets of $\Omega$ over $\Sigma$ is a
    \textbf{labeled sunflower} if:
  \begin{itemize}
  \item $\Fc$ is a sunflower (as in Definition \ref{sunflower-def}),
    whose core we will call $C$, and
  \item for each element $v \in C$, each labeled subset $A \in \Fc$ assigns the same label to $v$;
    that is, for each $v \in C$, there exists a label $\rho \in \Sigma$ such that $A(v) = \rho$
    for each $A\in\Fc$.
  \end{itemize}
\end{definition}

It turns out that the textbook proof of the Sunflower Lemma \cite{ParamAlgo} can be
  generalized to produce labeled sunflowers in a straightforward manner.

\begin{lemma}[Labeled Sunflower Lemma]
  \label{labeled-sunflower-lemma}
  For a finite set of labels $\Sigma$, a universe $\Omega$ and two integers $a \geq 1$,
    $b \geq 0$, any family of $ab!(a|\Sigma|)^b$ (not necessarily different) labeled subsets of
    $\Omega$ over $\Sigma$, each of cardinality at most $b$,
  contains a labeled sunflower of order $a$ as a subset.
  \begin{proof}
  We apply an induction on $b$.
  If $b = 0$, then the family contains at least $ab!(a|\Sigma|)^b = a$ empty sets and already
    is a labeled sunflower.

  Now, assume that $b \geq 1$, and let $\Fc$ be the family as in the statement of the lemma.
  We consider two cases:
  \begin{itemize}
    \item There exists an element $v \in \Omega$ belonging to at least $(b-1)!(a|\Sigma|)^b$
      labeled sets in~$\Fc$.
      Hence, there exists a label $\rho \in \Sigma$ and a subfamily $\Fc' \subseteq \Fc$
        such that we have $|\Fc'| \geq a(b-1)!(a|\Sigma|)^{b-1}$, and each set $A \in \Fc'$
        contains $v$ and satisfies $A(v) = \rho$.
        
      We create a modified universe $\Omega' := \Omega \setminus \{v\}$, and form a family
        $\Fc''$ from $\Fc'$ by removing $v$ from each of the sets within $\Fc'$.
      Then, Labeled Sunflower Lemma applies inductively to $\Fc''$ and parameters $a$ and $b-1$.
      Hence, we can find a labeled sunflower $S$ of order $a$ as a subset of $\Fc''$.
      By reintroducing $v$ to each of the sets in $S$, we find a labeled sunflower $S'$ of order $a$
        as a subset of $\Fc'$.
        
    \item Each element $v \in \Omega$ occurs in fewer than $(b-1)!(a|\Sigma|)^b$
      labeled sets in $\Fc$.
      Since $(b-1)!(a|\Sigma|)^b \leq \frac{1}{ab} |\Fc|$, we can find a family $\Gc$ of $a$ pairwise
        disjoint labeled sets in $\Fc$ greedily.
      In each of $a$ steps, we extend $\Gc$ by any labeled set $A \in \Fc$, and remove
        from $\Fc$ all labeled sets intersecting $A$.
      As each step removes at most $b!(a|\Sigma|)^b \leq \frac1a |\Fc|$ elements from $\Fc$,
        we can see that this is a valid construction of $\Gc$.
      Since the labeled sets in $\Gc$ are pairwise disjoint, they form a labeled sunflower.
  \end{itemize}
  In both cases, we construct a labeled sunflower of order $a$ as a subset of $\Fc$, so the
    proof is complete.
  \end{proof}
\end{lemma}

We remark that the lemma requires a relatively small multiplicative overhead on the size of
  the initial family of subsets, compared to the original statement of the lemma (Theorem
  \ref{sunflower-lemma}):
  $$ \frac{ab!(a|\Sigma|)^b}{b!a^{b+1}} = |\Sigma|^b. $$
As a sidenote, there exists a simpler proof of a looser upper bound:
  find a sunflower of order $a|\Sigma|^b$ using the sunflower lemma,
  and then find a subset of this sunflower of order $a$ which is a~labeled sunflower
  (which is an easy task since there are at most $|\Sigma|^b$ labelings of the core of the
  sunflower).
This proof, however, requires the multiplicative overhead equal to $|\Sigma|^{b(b+1)}$
  compared to the Sunflower Lemma.
If we used this version of the Labeled Sunflower Lemma later instead of the version
  we have just proved in Lemma \ref{labeled-sunflower-lemma}, we would ultimately derive
  the $d^{O(p^2)}p^{O(p)}$ upper bound on the maximum distance-$d$ semi-ladder size
  in graphs with pathwidth bounded by $p$, which is unfortunately too high for our needs.

\subsubsection*{Upper bound proof}

In the first part of our proof, we will define a structure whose occurrence in a graph
  proves that the graph both contains a large semi-ladder and has low pathwidth.

\begin{definition}
In a graph $G$, for $p, d, k, \ell \geq 1$,
a path decomposition $(W_1, W_2, \dots, W_k)$ of width $p$, together with $2\ell$
  different vertices $a_1, a_2, \dots, a_\ell, b_1, b_2, \dots, b_\ell$,
  and $\ell$ different indices $t_1, t_2, \dots, t_\ell \in [1, k]$,
  is a \textbf{distance-$d$ alignment} of order $\ell$ if:
\begin{itemize}
  \item vertices $a_1, a_2, \dots, a_\ell, b_1, b_2, \dots, b_\ell$ form a distance-$d$
    semi-ladder in $G$;
  \item vertices $a_1, a_2, \dots, a_\ell$ belong to the bags $W_{t_1}, W_{t_2}, \dots,
    W_{t_\ell}$, respectively.
\end{itemize}

We denote a distance-$d$ alignment in the following way:

$$ \mathsf{Alignment}_d\left( (W_i)_{i=1}^k, (a_i)_{i=1}^\ell, (b_i)_{i=1}^\ell, (t_i)_{i=1}^\ell\right).$$
\end{definition}

We note that given a distance-$d$ alignment
  $\mathsf{Alignment}_d\left( (W_i)_{i=1}^k, (a_i)_{i=1}^\ell, (b_i)_{i=1}^\ell, (t_i)_{i=1}^\ell\right)$,
  one can easily take a subset of the underlying semi-ladder, given by indices
  $i_1, i_2, \dots, i_m$ such that $1 \leq i_1 < i_2 < \dots < i_m \leq \ell$,
  and form a distance-$d$ alignment with order $m$:
  $\mathsf{Alignment}_d(\allowbreak
   (W_i)_{i=1}^k, (a_{i_j})_{j=1}^m, (b_{i_j})_{j=1}^m, (t_{i_j})_{j=1}^\ell)$.

In our proof, we will transform a graph with bounded pathwidth, and containing a large distance
  distance-$d$ semi-ladder, into a distance-$d$ alignment containing this semi-ladder.
Then, we will find a reasonably large subset of this alignment having necessary structural
  properties.
The following definition formalizes these properties:

\begin{definition}
For $p, d, k, \ell \geq 1$, a distance-$d$ alignment of width $p$ and order $\ell$:
  $$\mathsf{Alignment}_d(\allowbreak (W_i)_{i=1}^k, (a_i)_{i=1}^\ell, (b_i)_{i=1}^\ell, (t_i)_{i=1}^\ell)$$
  is also a \textbf{distance-$d$ sunflower alignment} if the following additional
  properties hold:

\begin{itemize}
  \item bags $W_{t_1}, W_{t_2}, \dots, W_{t_\ell}$ form a sunflower, whose core we will
    denote by $C$;
  \item in $G$, the distance-$d$ profiles on $C$ are identical for each vertex
    $a_1, a_2, \dots, a_\ell$.
\end{itemize}

We will refer to distance-$d$ sunflower alignment in the following way:

$$ \mathsf{SunflowerAlignment}_d\left( (W_i)_{i=1}^k, (a_i)_{i=1}^\ell, (b_i)_{i=1}^\ell,
  (t_i)_{i=1}^\ell\right).$$
  
The set $C$ is the \textbf{core} of the decomposition.
\end{definition}

We are now ready to find large distance-$d$ sunflower alignments in the graphs
  with large distance-$d$ semi-ladders.

\begin{lemma}
\label{large-sunflower-lemma}
For $p, d, \ell \geq 1$,
  every graph with pathwidth not exceeding $p$ and containing a semi-ladder of order
  $\ell(p+1)! [\ell(d+2)]^{p+1}$
  also contains a distance-$d$ sunflower path decomposition of order $\ell$.
  \begin{proof}
  Let $L := \ell(p+1)! [\ell(d+2)]^{p+1}$.
  Fix a graph $G$, together with its distance-$d$ semi-ladder $a_1, a_2, \dots, a_L$,
    $b_1, b_2, \dots, b_L$.
  Since $\pw{G} \leq p$, there exists a path decomposition $W = (W_1, W_2, \dots,
    W_k)$ of $G$ of width at most $p$; that is, we have $|W_i| \leq p + 1$ for each $i \in [1, k]$.
  Moreover, we can assume that $W$ is a \textbf{nice path decomposition}; that is, we have
    $W_1 = W_k = \varnothing$, and for
    each index $i \in \{2, 3, \dots, k\}$, either of the following conditions is satisfied:
  \begin{itemize}
    \item $W_i$ \textbf{introduces} a vertex $v \in V(G)$; that is, $W_i = W_{i-1} \cup \{v\}$;
    \item $W_i$ \textbf{forgets} a vertex $v \in V(G)$; that is, $W_i = W_{i-1} \setminus \{v\}$.
  \end{itemize}
  It can be easily shown that each path decomposition can be converted into a
    nice path decomposition of the same width \cite{ParamAlgo},
    so we can safely assume that $W$ is a nice path decomposition.

  For each vertex $a_j$ ($1 \leq j \leq L$) of the semi-ladder, we can locate the bag
    $W_{t_j}$ introducing $a_j$ to the decomposition.
  Since the path decomposition is nice, the values $t_j$ are different for all $j \in [1, L]$.
  Therefore, the following describes a distance-$d$ alignment of order $L$:
  $$ \mathsf{Alignment}_d\left((W_i)_{i=1}^k, (a_i)_{i=1}^L, (b_i)_{i=1}^L, (t_i)_{i=1}^L \right). $$
  
  We will now find a subset of this distance-$d$ alignment which is a distance-$d$
    sunflower alignment.
    
  We first create a family $\Fc$ of $L$ labeled subsets of $V(G)$ over $\{0, 1, 2, \dots, d, +\infty
    \}$: for each $i \in [1, L]$, we let the $i$-th labeled subset in the family to be the bag
    $W_{t_i}$ (containing the vertex $a_i$), where each vertex $v \in W_{t_i}$ has the label
    $\pi_d[a_i, W_{t_i}](v)$.
    
  Secondly, since each set in $\Fc$ contains at most $p + 1$ elements, the cardinality
    of the set of labels is equal to $d + 2$, and $|\Fc| = L = \ell(p+1)! [\ell(d+2)]^{p+1}$,
    we infer from the Labeled Sunflower Lemma
    (Lemma \ref{labeled-sunflower-lemma}) that there exists a subfamily $\Fc' \subseteq \Fc$,
    which is a labeled sunflower of order $\ell$.
  
  Finally, we let $\Fc' = (W_{t_{i_1}}, W_{t_{i_2}}, \dots, W_{t_{i_\ell}})$ for $i_1 < i_2 < \dots < i_\ell$.
  We can now easily verify that the following structure is a distance-$d$ sunflower alignment
    of order $\ell$:
  \[ \mathsf{SunflowerAlignment}_d\left( (W_j)_{j=1}^k, (a_{i_j})_{j=1}^\ell, (b_{i_j})_{j=1}^\ell,
    (t_{i_j})_{j=1}^\ell \right). \]
  \end{proof}
\end{lemma}

Now, it turns out that we can give an explicit upper bound on the order of every
  distance-$d$ sunflower alignment:

\begin{lemma}
\label{no-large-sunflower-lemma}
  For $d \geq 1$, no graph can contain a distance-$d$ sunflower alignment of order $2d + 3$.
  \begin{proof}
    Fix $\ell := 2d + 3$ and a graph $G$ containing a distance-$d$ sunflower alignment
      of order $\ell$ and width $p$:
    $$\mathsf{SunflowerAlignment}_d\left( (W_i)_{i=1}^k, (a_i)_{i=1}^\ell, (b_i)_{i=1}^\ell,
      (t_i)_{i=1}^\ell \right). $$
    We will prove that the shortest path connecting $b_1$ with $a_i$ for some $i \in [2, \ell]$
      has length greater than $d$, which contradicts the assumption that
      $a_1, a_2, \dots, a_\ell, b_1, b_2, \dots, b_\ell$ is a distance-$d$ semi-ladder.
      
    Let $C$ be the core of the sunflower alignment.
    Firstly, we will prove that no shortest path from $b_1$ to $a_i$ for $i \in [2, \ell]$ passes
      through $C$.
    Assume for contradiction that for some $i \in [2, \ell]$, the shortest path from $b_1$ to $a_i$
      contains a vertex $x \in C$ (Figure \ref{pw-core-intersection-fig}).
    Since $\dist(x, a_i) \leq d$ and $\pi_d[a_1, C] = \pi_d[a_i, C]$, we infer that the vertex $x$ is
      equidistant from $a_1$ and $a_i$.
    Hence,
    $$ \dist(b_1, a_1) \leq \dist(b_1, x) + \dist(x, a_1) = \dist(b_1, x) + \dist(x, a_i) = \dist(b_1, a_i)
      \leq d. $$
    We have a contradiction as $\dist(b_1, a_1) > d$ is required by the definition of a semi-ladder.
    We note that this fact means that, in particular, neither $b_1$ nor any of the vertices $a_2, a_3,
      \dots, a_\ell$ belongs to $C$.
      
    \begin{figure}[h]
      \centering
      \input{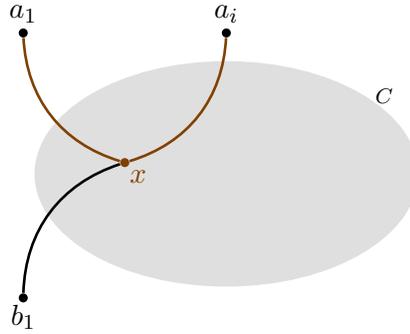}
      \caption{The setup in the proof that no shortest path from $b_1$ to $a_i$ intersects $C$.}
      \label{pw-core-intersection-fig}
    \end{figure}

    Next, we arrange the vertices $a_2, a_3, \dots, a_\ell$ in the order in which they appear in $W$.
    Formally, we define a permutation $\sigma = (\sigma_1, \sigma_2, \dots, \sigma_{\ell-1})$ of
      $\{2, 3, \dots, \ell\}$ such that $t_{\sigma_1} < t_{\sigma_2} < \dots < t_{\sigma_{\ell - 1}}$.
    This way, for each pair of indices $i, j$ ($i, j \in [1, \ell - 1]$) with $i < j$,
      the bag $W_{t_{\sigma_i}}$ containing $a_{\sigma_i}$ appears earlier in $W$ than the
      bag $W_{t_{\sigma_j}}$ containing $a_{\sigma_j}$ (as $t_{\sigma_i} < t_{\sigma_j}$).
      
    Let us find the vertex $b_1$ in the path decomposition.
    We already know that $b_1 \not\in C$.
    Since $W_{t_{\sigma_{d+1}}} \cap W_{t_{\sigma_{d+2}}} = C$, we deduce that $b_1$ may
      belong to at most one of the bags $W_{t_{\sigma_{d+1}}}$, $W_{t_{\sigma_{d+2}}}$.
    As the definition of a path decomposition requires that $b_1$ belongs to the family of the
      bags forming a subinterval of the path decomposition,
      we get that at least one of the following cases is satisfied:
    \begin{enumerate}[(a)]
      \item $b_1$ does not belong to any of the bags
        $W_1, W_2, W_3, \dots, W_{t_{\sigma_{d+1}}}$ of the decomposition;
      \item $b_1$ does not belong any of the bags $W_{t_{\sigma_{d+2}}}, \dots,
        W_{k-2}, W_{k-1}, W_k$ of the decomposition.
    \end{enumerate}
    As these cases are symmetrical, we assume without loss of generality that the case (a) holds.
    With this assumption, we consider the shortest path between $b_1$ and $a_{\sigma_1}$.
    We remark two facts about this path:
    \begin{itemize}
      \item it cannot pass through $C$, and
      \item it has to intersect each of the bags $W_{t_{\sigma_1}}, W_{t_{\sigma_2}}, \dots,
        W_{t_{\sigma_{d+1}}}$ of the decomposition
        (since $t_{\sigma_1} < t_{\sigma_2} < \dots < t_{\sigma_{d+1}}$, and
        $a_{\sigma_1} \in W_{t_{\sigma_1}}$).
    \end{itemize}
    Since each pair of the bags $W_{t_{\sigma_1}}, W_{t_{\sigma_2}}, \dots, W_{t_{\sigma_{d+1}}}$
      intersects exactly at $C$, and the path cannot pass through $C$,
      we conclude that the path has to intersect each $W_{t_{\sigma_1}},
      W_{t_{\sigma_2}}, \dots, W_{t_{\sigma_{d+1}}}$ at a~different vertex.
    By our assumption that case (a) holds, $b_1$ does not belong to any of these bags.
    We conclude that the path must contain at least $d + 2$ vertices,
      hence its length is at least $d+1$ --- a contradiction.

    This contradiction completes the proof of the lemma.
  \end{proof}
\end{lemma}

This fact immediately allows us to conclude with a concrete upper bound.

\begin{theorem}
\label{pw-upper-bound}
  For $p \geq 1$, $d \geq 1$, no graph with pathwidth bounded from above by $p$
    contains a distance-$d$ semi-ladder of order
  $$ (2d + 3)(p + 1)! [(2d + 3)(d + 2)]^{p + 1}. $$
  \begin{proof}
    If a graph with pathwidth bounded by $p$ contained a distance-$d$ semi-ladder of this
      order, it would also contain a distance-$d$ sunflower alignment of order
      $2d + 3$ (Lemma~\ref{large-sunflower-lemma}).
    However, by Lemma \ref{no-large-sunflower-lemma}, such an alignment cannot exist.
  \end{proof}
\end{theorem}

We also remark the asymptotic version of this bound.

\begin{corollary}
\label{pw-upper-bound-asym}
  There exists a constant $\delta > 0$ such that for each $p, d \geq 2$, the maximum order
    of a distance-$d$ semi-ladder in the class of graphs with pathwidth bounded by $p$
    is bounded from above by $(pd)^{\delta p}$.
\end{corollary}

We recall the construction attesting the $d^{\Omega(p)}$ lower bound in graphs
  with pathwidth bounded by $p$ (Corollary \ref{pw-lower-bound-asym}).

\subsection{Graphs with bounded treewidth and minor-free graphs}
\label{kt-up-section}

In this part, we will prove an upper bound on the maximum order of a distance-$d$
  semi-ladder in the class of graphs excluding $K_t$ as a minor; an analogous bound for graphs
  with bounded treewidth will follow.
  
We begin by introducing two concepts of Sparsity that will become useful in the proof:
  \textbf{weak coloring numbers} and \textbf{uniform quasi-wideness}.

\begin{definition}[weak coloring numbers]
\label{wcol-def}
We follow the definition by Zhu \cite{DBLP:journals/dm/Zhu09}.

Let $\Sigma_{V(G)}$ denote the set of all permutations of $V(G)$.
For a permutation $\sigma \in \Sigma_{V(G)}$ and two vertices $u, v \in V(G)$, we write
  $u <_\sigma v$ if $u$ occurs in $\sigma$ before~$v$.

For a fixed distance $d \geq 1$, we say that $v$ is \textbf{weakly $d$-reachable} from $u$
  with respect to $\sigma \in \Sigma_{V(G)}$ if $v \leq_\sigma u$ and
  there exists a path of length at most $d$
  connecting $u$ and $v$ such that $v$ is the minimum vertex on this path with respect to
  $\sigma$.
We remark that $u$ is weakly $d$-reachable from itself for every $d \in \N$.
We denote the set of vertices weakly $d$-reachable from $u$ with respect to $\sigma$
  as $\mathrm{WReach}_d[G, \sigma, u]$ (Figure \ref{weak-reach-def-fig}).

\begin{figure}[h]
\centering
\input{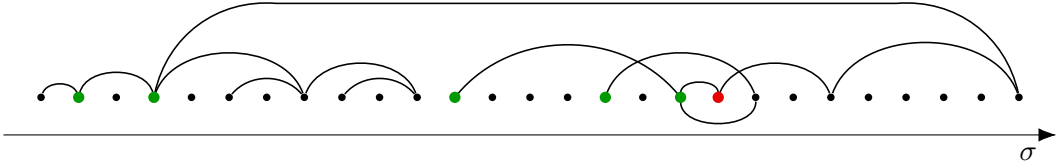}
\caption{An example setup. In a permutation $\sigma$, for a vertex $u$ (red),
  $\mathrm{WReach}_4[G, \sigma, u] \setminus \{u\}$ is marked green.}
\label{weak-reach-def-fig}
\end{figure}

Given a graph $G$ and an order $\sigma \in \Sigma_{V(G)}$ we define the weak coloring number
  of $G$ with respect to $\sigma$ as $$\wcol_d(G, \sigma) :=
  \max_{u \in V(G)} \left|\mathrm{WReach}_d[G, \sigma, u]\right|.$$
  
Finally, given a graph $G$, we define the \textbf{weak coloring number} of $G$ as
  $$\wcol_d(G) := \min_{\sigma \in \Sigma_{V(G)}} \wcol_d(G, \sigma).$$
\end{definition}

\begin{definition}[uniform quasi-wideness]
\label{uniform-quasi-wideness-def}
A class $\Cc$ of graphs is \textbf{uniformly quasi-wide} with margins $N\,:\,\N\times\N\to\N$ and
  $s\,:\,\N\to\N$ if for every $m, d \geq 1$, graph $G\in\Cc$ and a subset $A$ of vertices such that
  $|A| \geq N(m, d)$, there exists a set of vertices $S \subseteq V(G)$, $|S| \leq s(d)$ and
  a~subset $B \subseteq A-S$, $|B| \geq m$, which is a distance-$d$ independent set in
  $G - S$.
\end{definition}

It turns out that a subgraph-closed class $\Cc$ of graphs is uniformly quasi-wide for any margins
  if and only if it is nowhere dense \cite{DBLP:journals/jsyml/NesetrilM10}.
For more properties and applications of uniform quasi-wideness, we refer to the lecture notes
  from the University of Warsaw \cite{SparsityUWNotes}.

In this section, we may provide the margins $N, s$ as partial functions, defined only
  for distances $d$ from an infinite subset $D \subsetneq \N$.
It is straightforward to see that these margins can be easily extended to total functions:
  given $m, d \in N$, we set $d' := \min \{x \in D\,\mid\,x \geq d\}$, and we put
  $N(m, d) := N(m, d')$ and $s(d) := s(d')$.
  
\medskip

In the proof of the upper bound for $K_t$-minor-free graphs,
  we will use the upper bound by Fabiański et al.~\cite{fabiaski2018progressive},
  applicable to every uniformly quasi-wide class of graphs:
\begin{lemma}[{{\cite[Lemma 29]{fabiaski2018progressive}}}]
  \label{semi-ladder-from-uqw}
  In every uniformly quasi-wide class of graphs $\Cc$ with margins $N\,:\,\N \times \N \to \N$ and
    $s\,:\,\N \to \N$, each semi-ladder in every graph $G \in \Cc$ has its order bounded by
    $N(2 \cdot (d + 2)^{s(2d)} + 1,\, 2d)$.
\end{lemma}

We will now prove that the class of $K_t$-minor-free graphs is uniformly quasi-wide by
  utilizing the generalized coloring numbers.
Precisely, we will use the following upper bound on the weak coloring numbers in this class,
  proved by van den Heuvel et al. \cite{DBLP:journals/endm/HeuvelMRS15}.
\begin{lemma}[\cite{DBLP:journals/endm/HeuvelMRS15}]
  \label{kt-wcol-bound}
  For $t \geq 4$ and any graph $G$ which excludes $K_t$ as a minor, we have that
    $\wcol_d(G) \leq \binom{d+t-2}{t-2}(t - 3)(2d + 1)$.
\end{lemma}

We only need to prove that if a class has fairly low weak coloring numbers, then it is
  also uniformly quasi-wide with reasonably low margins.

\begin{lemma}
  \label{wcol-gives-uqw}
  Assume that a class of graphs $\Cc$ has a finite weak coloring number
    $wcol_d(\Cc)$ for each $d \geq 1$.
  Then $\Cc$ is also uniformly quasi-wide with margins
    $$N(m, d) = (\wcol_d(\Cc))!\,(m+1)^{\wcol_d(\Cc) + 1}, \quad
      s(d) = \wcol_d(\Cc) \quad \text{for each }d \geq 1,\, m \geq 2. $$
  \begin{proof}
  Fix integers $d \geq 1, m \geq 2$, a graph $G\in\Cc$, and a subset $A$ of the vertices of $G$
    containing at least $N(m, d)$ vertices.
  We will prove the lemma by showing that we can erase at most $\wcol_d(\Cc)$ vertices
    from $G$ in order to find a distance-$d$ independent set $B \subseteq A$
    of size $m$ in the remaining part of the graph.
  
  In $G$, we can find an ordering $\sigma$ of its vertices for which
    $\wcol_d(G, \sigma) \leq \wcol_d(\Cc)$.
  Knowing that, we define a family $\Fc$ of subsets of vertices:
  $$ \Fc := \{R_v\,\mid\, v \in S\} \qquad\text{where }R_v := \mathrm{WReach}_d[G, \sigma, v]
    \text{ for } v \in S. $$
  Since each set in the family contains at most $\wcol_d(\Cc)$ elements and
    $|\Fc| \geq N(m, d) = (\wcol_d(\Cc))!\,(m+1)^{\wcol_d(\Cc) + 1}$,
    the Sunflower Lemma (Theorem \ref{sunflower-lemma})
    lets us find a subfamily $\Fc' \subseteq \Fc$ which is a sunflower of order $m + 1$.
  We let $\Fc' = \{R_{v_1}, R_{v_2}, \dots, R_{v_{m+1}}\}$ for $m+1$ different vertices
    $v_1, v_2, \dots, v_{m+1}$. We denote the core of $\Fc'$ as $D$.

  Without loss of generality, we assume that the vertices $v_1, v_2, \dots, v_{m+1}$ are ordered
    decreasingly with respect to the order $\sigma$: $v_{m+1} <_\sigma v_m <_\sigma \dots
    <_\sigma v_1$.
  By the definition of weak reachability in Definition \ref{wcol-def},
    for each $i \in [1, m + 1]$, $v_i$ is the maximum element of $R_{v_i}$
    with respect to $\sigma$.
  Since for $i \in [1, m]$ we have $v_i >_\sigma v_{m + 1}$, and the core $D$ must be a subset
    of $R_{v_{m+1}}$, we infer that $v_i \not\in D$ for each $i \in [1, m]$.
  Therefore, the set $\{v_1, v_2, \dots, v_m\}$ of vertices is disjoint with $D$.
  We now focus on the sunflower $\Fc'' := \{R_{v_1}, R_{v_2}, \dots, R_{v_m}\}$ of order $m$,
    whose core is naturally $D$.
  Obviously, $|D| \leq \wcol_d(\Cc)$.
  
  We define the graph $G' := G - D$.
  We will now prove that $\{v_1, v_2, \dots, v_m\}$ is a distance\nobreakdash-$d$
    independent set in $G'$.
  Assume for contradiction that there exists a path of length at most $d$ connecting two
    different vertices of the set --- say, $v_i, v_j$ for two different indices $i, j \in [1, m]$.
  Let $w_{\min}$ be the vertex on this path which is minimal with respect to $\sigma$.
  This, however, implies that $w_{\min}$ belongs to both $\mathrm{WReach}_d[G, \sigma, v_i]$
    and $\mathrm{WReach}_d[G, \sigma, v_j]$, and thus $w_{\min} \in R_{v_i} \cap R_{v_j}$.
  Since $\Fc'' = \{R_{v_1}, R_{v_2}, \dots, R_{v_m}\}$, we infer that
    $w_{\min}$ belongs to $D$ --- the core of $\Fc''$.
  This is a contradiction since $D$ is disjoint with $G'$.
  \end{proof}
\end{lemma}

We remark that a combination of Lemmas \ref{semi-ladder-from-uqw}, \ref{kt-wcol-bound} and
  \ref{wcol-gives-uqw} can already provide us a concrete upper bound on the maximum
  order of a distance-$d$ semi-ladder in the class of $K_t$-minor-free graphs.
However, the produced bound will be insufficient for our purposes --- it can be verified that
  it would result in the $d^{O(d^{2t - 2})}$ upper bound for each fixed $t \geq 4$, while
  we aim at the $d^{O(d^{t - 1})}$ bound.

In order to prove a tighter upper bound, we will improve  the margins $N, s$ with which the
  class of $K_t$-minor-free graphs is uniformly quasi-wide using the lemma below.
This lemma is inspired by Lemma 3.16 from the work of Nešetřil and Ossona de Mendez
  \cite{DBLP:journals/jsyml/NesetrilM10}.
However, we significantly modified the original lemma in order to ensure a huge drop of the
  margin $s$ while keeping the margin $N$ relatively low.

\begin{lemma}
  \label{kt-reduce-uqw}
  Fix $t \geq 4$. Assume that a class of graphs $\Cc$:
  \begin{itemize}
  \item excludes the complete bipartite graph $K_{t-1, t-1}$ as a minor, and
  \item is uniformly quasi-wide with margins $N\,:\,\N\times \N\to\N$, $s\,:\,\N\to\N$.
  \end{itemize}
  Then $\Cc$ is also uniformly quasi-wide with margins $\widehat{N}\,:\,\N \times \N \to \N$
    and $\widehat{s}\,:\,\N \to \N$, where $\widehat{s}$~is a constant function,
    $\widehat{s} \equiv t - 2$, and
  $$ \widehat{N}(m, 2d) = N\left(\binom{s(2d) + 1}{t - 1}(m+t),\,2d\right) \quad
    \text{for each }m, d \geq 1.$$
    
  \begin{proof}    
    \vspace{0.5em}
    We fix a graph $G \in \Cc$, integers $m, d \geq 1$, and a set of vertices $A \subseteq V(G)$ of
      size at least $\widehat{N}(m, 2d)$.
    By our assumption, there exists a~subset $S \subseteq V(G)$ with $|S| \leq s(2d)$
      and a~subset $B \subseteq A-S$ with $|B| \geq \binom{s(2d) + 1}{t - 1}(m+t)$
      such that $B$ is a distance-$2d$ independent set in $G - S$.
      
    Similarly to \cite[Lemma 3.16]{DBLP:journals/jsyml/NesetrilM10},
      for each $v \in B$, we let $L_v$ to be a minimal subset of $S$
      such that in $G - L_v$, there is no path of length $d$ or less from $v$
      to any vertex from $S \setminus L_v$; such a subset exists since $S$ itself satisfies
      this requirement.
    Now, we will define a new set $L'_v$ in the following way:
      if $|L_v| \leq t - 2$, then we set $L'_v := L_v$; otherwise, we take $L'_v$ to be any subset of
      $L_v$ of size exactly $t - 1$.
    Since $L'_v \subseteq S$, $|L'_v| \leq t - 1$ and $|S| \leq s(2d)$, we conclude that there
      are at most
    $$\binom{s(2d)}{0} + \binom{s(2d)}{1} + \dots + \binom{s(2d)}{t - 1} =\binom{s(2d) + 1}{t - 1}$$
    different possible sets $L'_v$.
    It follows that for some set $C \subseteq S$, $|C| \leq t - 1$, there are at least $m + t$
      different vertices $v_1, v_2, \dots, v_{m + t}$ for which $L'_{v_1} = L'_{v_2} = \dots =
      L'_{v_{m+t}} = C$.
      
    We consider two cases, distinguished by the cardinality of $C$:
    \begin{itemize}
      \item If $|C| \leq t - 2$, then $L_{v_1} = L_{v_2} = \dots = L_{v_{m + t}} = C$.
      Now, assume there are two different vertices $v_i, v_j$ ($i, j \in [1, m + t]$) at distance
        at most $2d$ in $G - C$.
      Since $v_i$ and $v_j$ are at distance greater than $2d$ in $G - S$, the shortest path
        between these two vertices in $G - C$ must pass through some vertex $x \in S \setminus C$.
      This, however, means that $x$ is at distance at most $d$ from one of the vertices $v_i$,
        $v_j$ in $G - C$ --- a contradiction since $x \not\in L_{v_i}$ and $x \not\in L_{v_j}$.
        
      Therefore, $\{v_1, v_2, \dots, v_{m+t}\}$ is a distance-$2d$ independent set in $G - C$,
        where $|C| \leq t - 2$.
        
      \item If $|C| = t - 1$, then we will find $K_{m + t, t - 1}$ as a minor in $G$, which will lead
        to a contradiction with $G \in \Cc$.
        Namely, for each $x \in \{v_1, v_2, \dots, v_{m+t}\}$, we construct a tree $Y_x$ rooted at $x$
          with depth not exceeding $d$ whose set of leaves is $L_x$ and that does not contain
          any other vertices of $S$.
        Such a tree exists as the minimality of $L_x$ implies that for each $v \in L_x$,
          the distance between $x$ and $v$ in $G - (L_x \setminus \{v\})$ does not exceed $d$.

        We see that no pair of trees $Y_x$, $Y_y$ can intersect at any vertex outside of $S$
          --- otherwise, we would be able to construct a path of length at most $2d$
          between $x$ and $y$ avoiding $S$.
        This is, of course, impossible as $x, y \in A$.
        
        Finally, we can find $K_{m+t, t-1}$ as a minor of $G$ by contracting all subgraphs
          $Y_x - S$ into single vertices, and then observing
          that for each $x \in \{v_1, v_2, \dots, v_{m+t}\}$ and $y \in C$, the vertex $y$ is adjacent
          to the contracted subgraph $Y_x - S$.
        We have a contradiction since we assumed $G$ to be $K_{t-1, t-1}$-free.
    \end{itemize}
  \end{proof}
\end{lemma}

We remark that Lemma \ref{kt-reduce-uqw} also applies if $\Cc$ excludes $K_t$ as a minor
  for some fixed $t \geq 4$, for such classes of graphs exclude $K_{t-1, t-1}$ as a minor as well.

With our additional lemma in hand, we can prove the following upper bound:

\begin{theorem}
  \label{kt-upper-bound}
  There exists a polynomial $p$ such that for $t \geq 4$, $d \geq 2$,
    there is no distance-$d$ semi-ladder of order larger than
  $$ d^{p(t) \cdot (2d+1)^{t-1}} $$    
  in any $K_t$-minor-free graph.
  
  \begin{proof}
  In this proof, $\poly(t)$ will denote any polynomial of $t$ with constant coefficients independent
    on $d$ or $t$.
  However, each occurrence of this term in any equation might mean a different polynomial.
  
  Fix $t \geq 4$ and $d \geq 2$, and consider the class $\Cc$ of graphs excluding $K_t$ as
    a minor.
  We first bound the weak coloring number of $\Cc$ using Theorem \ref{kt-wcol-bound}:
  \begin{equation}
  \label{kt-wcol-ineq}
  \begin{split}
  \wcol_d(\Cc) &\leq \binom{d + t - 2}{t - 2}(t - 3)(2d + 1) =
    \frac{(d + 1)(d + 2) \dots (d+t-2)}{(t - 2)!} (t - 3)(2d + 1) = \\
    &= \frac{d+1}{1} \cdot \frac{d+2}{2} \cdot \ldots \cdot \frac{d+t-2}{t-2} \cdot
      (t - 3)(2d + 1) \leq \\ &\leq (d + 1)^{t-2}(t-3)(2d+1) \leq
    \poly(t) \cdot (d + 1)^{t - 1}.
  \end{split}
  \end{equation}

  By Lemma \ref{wcol-gives-uqw}, $\Cc$ is uniformly quasi-wide with margins
  \begin{equation}
  \label{kt-init-nmd-ineq}
  \begin{split}
  N(m, d) & = (\wcol_d(\Cc))!(m + 1)^{\wcol_d(\Cc) + 1} \leq
    [(m + 1)\wcol_d(\Cc)]^{\wcol_d(\Cc) + 1} = \\
    & \stackrel{(\ref{kt-wcol-ineq})}{=}
      [m \cdot \poly(t) \cdot (d + 1)^{t - 1}]^{\poly(t) \cdot (d + 1)^{t - 1}} \leq
    (md)^{\poly(t) \cdot (d + 1)^{t - 1}}
  \end{split}
  \end{equation}
  and
  \begin{equation}
  \label{kt-init-sd-ineq}
  s(d) = \wcol_d(\Cc) \stackrel{(\ref{kt-wcol-ineq})}{\leq} \poly(t) \cdot (d + 1)^{t - 1}.
  \end{equation}
  
  In order to apply Lemma \ref{kt-reduce-uqw}, we first bound the following expression for
    $m \geq 1$:
  \begin{equation}
  \label{kt-init-binom-ineq}
  \begin{split}
  \binom{s(2d) + 1}{t - 1}(m + t) & \leq
    [s(2d) + 1]^{t - 1} (m + t) \stackrel{(\ref{kt-init-sd-ineq})}{\leq}
    [\poly(t) \cdot (2d + 1)]^{(t - 1)^2} (m + t) \leq \\
    & \leq m \cdot [d \cdot \poly(t)]^{(t - 1)^2}.
  \end{split}
  \end{equation}
  
  Now, we get that $\Cc$ is also uniformly quasi-wide with margins $\widehat{s} \equiv t - 2$ and
  \begin{equation}
  \label{kt-later-nmd-ineq}
  \begin{split}
    \widehat{N}(m, 2d) & = N\left(\binom{s(2d) + 1}{t - 1}(m+t),\,2d\right)
      \stackrel{(\ref{kt-init-binom-ineq})}{\leq}
      N\left(m \cdot [d \cdot \poly(t)]^{(t - 1)^2},\, 2d\right) \leq \\
    & \stackrel{(\ref{kt-init-nmd-ineq})}{\leq}
      \left\{ m \cdot [d \cdot \poly(t)]^{(t - 1)^2} \cdot 2d \right\}^{\poly(t) \cdot
      (2d + 1)^{t - 1}} \leq (md)^{\poly(t) \cdot (2d + 1)^{t - 1}}.
  \end{split}
  \end{equation}
  
  Finally, we apply Lemma \ref{semi-ladder-from-uqw} and infer that the maximum order
    of any distance-$d$ semi-ladder is bounded from above by
  $$ \widehat{N}\left(2 \cdot (d + 2)^{t - 2} + 1, 2d\right) \stackrel{(\ref{kt-later-nmd-ineq})}{\leq}
    d^{\poly(t) \cdot (2d + 1)^{t - 1}}. $$
  
  Hence, there exists a polynomial $p$ for which this maximum order is bounded from
    above by $d^{p(t) \cdot (2d+1)^{t - 1}}$ for all $t \geq 4$, $d \geq 2$.
  \end{proof}
\end{theorem}

We also infer the following asymptotic version of Theorem \ref{kt-upper-bound}:
\begin{corollary}
  \label{kt-upper-bound-asym}
  For each fixed $t \geq 4$, the maximum order of a distance-$d$ semi-ladder in the class
    of graphs excluding $K_t$ as a minor is bounded from above by $d^{O(d^{t - 1})}$.
\end{corollary}
We contrast this result with the $2^{d^{\Omega(t)}}$ lower bound
  (Corollary \ref{kt-lower-bound-asym}) in the class of $K_t$-minor-free graphs.
  
\medskip

By Theorem \ref{tw-minor-free}
  (each graph $G$ with $\tw(G) \leq t$ is also $K_{t+2}$-minor-free),
  we immediately prove the following result from Corollary \ref{kt-upper-bound-asym}:
\begin{corollary}
  \label{tw-upper-bound-asym}
  For each fixed $t \geq 2$, the maximum order of a semi-ladder in the class of graphs with
    treewidth not exceeding $t$ is bounded from above by $d^{O(d^{t + 1})}$.
\end{corollary}
Similarly, we recall the $2^{d^{\Omega(t)}}$ lower bound on the maximum order
  of distance-$d$ semi-ladders in the class of graphs with treewidth bounded by $t$
  (Corollary \ref{tw-lower-bound-asym}).

\section{Conclusions}\label{conclusions-chapter}
In this work we gave asymptotic lower and upper bounds on the maximum orders
  of distance-$d$ half graphs and semi-ladders that can appear in various classes of sparse graphs:
  graphs with bounded maximum degree, planar graphs, graphs with bounded pathwidth
  or treewidth, and graphs excluding $K_t$ as a minor.
All these bounds are asymptotically almost tight and allow us to determine how large structures
  of this type may appear in graphs from these classes.
Graphs with degree bounded by $\Delta$ and planar graphs admit exponentially large
  distance-$d$ half graphs.
Graphs with pathwidth bounded by $p$ admit polynomially large half graphs, with the degree
  of the polynomial depending linearly on $p$.
Finally, the maximum order of distance-$d$ half graphs in graphs with treewidth
  bounded by $t$ and in $K_t$-minor-free graphs grows exponentially on~$d$,
  where the exponent is a polynomial of degree depending linearly on~$t$.

There are a few natural open problems regarding half graphs and semi-ladders in sparse classes of graphs and related topics:
\begin{itemize}
\item In case of the maximum order of distance-$d$ semi-ladders in planar graphs,
  it would be interesting to close the asymptotic gap between $2^{\Omega(d)}$ and $d^{O(d)}$.
\item We can easily see that in an arbitrary class~$\Cc$ with bounded expansion,
  the combination of Lemmas \ref{semi-ladder-from-uqw} and \ref{wcol-gives-uqw}
  results in a $(d \cdot \wcol_{2d}(\Cc))^{O(\wcol_{2d}(\Cc)^2)}$ upper bound on the maximum
  order of any distance-$d$ semi-ladder.
  However, we can deduce from Subsection \ref{kt-up-section} that as soon as $\Cc$~excludes
  a fixed graph as a minor, we can derive a much tighter $(d \cdot \wcol_{2d}(\Cc))^{
  O(\wcol_{2d}(\Cc))}$ upper bound thanks to Lemma \ref{kt-reduce-uqw}.
  Does the refined bound apply for all classes $\Cc$ with bounded expansion?
\item Can we characterize the classes $\Cc$ of graphs where the orders of distance-$d$
  semi-ladders are bounded polynomially in $d$? For instance, any class of
  graphs with bounded pathwidth has this property, but this is not the case for planar graphs or
  even for graphs with maximum degree bounded by $4$.
\item Theorem \ref{neighborhood-complexity-planar} proves that the neighborhood complexity
  in planar graphs is bounded by a polynomial of the size of the set $|X|$ and the distance $d$.
  The proof of Reidl et al.~\cite{DBLP:journals/ejc/ReidlVS19} provides a bound which is linear
  in $|X|$, but exponential in $d$.
  Can these bounds be reconciled?
  That is, can we bound the neighborhood complexity in planar graphs by a function which is
  linear in $|X|$ and polynomial in $d$?
\end{itemize}

\bibliographystyle{abbrv}
\bibliography{references}

\end{document}


\begin{tikzpicture}[auto, every loop/.style={},
                    top/.style={circle,minimum size=0.4cm,inner sep=0pt,draw},
                    smol/.style={circle,minimum size=0.12cm,inner sep=0pt,draw, gray}]

  \node [top] (r) at (0, 1) {$r$};
  \node [top] (b) at (0, -0.2) {$b_i$};
  \foreach \name/\x in {1/-1, 2/-0.5, 3/0, 4/0.5, 5/1} {
    \foreach \i/\y in {1/-1, 2/-1.5, 3/-2, 4/-2.5, 5/-3, 6/-3.5} {
      \node [smol] (m\name\i) at (\x, \y) {};
    };
    \foreach \i [count=\yi] in {2,...,6} {
      \draw [gray] (m\name\yi) -- (m\name\i);
    };
    \draw [gray] (b) -- (m\name 1);
  };
\end{tikzpicture}


\begin{tikzpicture}[auto, every loop/.style={},
                    top/.style={circle,minimum size=0.4cm,inner sep=0pt,draw},
                    smol/.style={circle,minimum size=0.12cm,inner sep=0pt,draw, gray}]

  \node [top] (r) at (0, 1) {$r$};
  \node [top] (b) at (0, -0.2) {$b_i$};
  \foreach \name/\x in {1/-1, 2/-0.5, 3/0, 4/0.5, 5/1} {
    \foreach \i/\y in {1/-1, 2/-1.5, 3/-2, 4/-2.5, 5/-3, 6/-3.5} {
      \node [smol] (m\name\i) at (\x, \y) {};
    };
    \foreach \i [count=\yi] in {2,...,6} {
      \draw [gray] (m\name\yi) -- (m\name\i);
    };
    \draw [gray] (b) -- (m\name 1);
  };
  \foreach \i/\y in {1/-1, 2/-1.5, 3/-2, 4/-2.5} {
    \node [smol, black] (m1\i) at (-1, \y) {};
  }
  \foreach \i [count=\yi] in {2,...,4} {
    \draw (m1\yi) -- (m1\i);
  };
  \draw (b) -- (m11);
  \draw (m14) to [bend left=20] (-1.6, -1.2) to [bend left=38] (r);
\end{tikzpicture}